\documentclass[11pt]{amsart}
\usepackage{graphpap}
\usepackage{amssymb}
\newtheorem{theoreme}{Theorem}

\newtheorem{proposition}{Proposition}
\newtheorem{lemme}[proposition]{Lemma}

\newtheorem{remarque}[proposition]{Remark}

\numberwithin{equation}{section}
\numberwithin{proposition}{section}

\def\Im{\textrm{Im}}
\def\Re{\textrm{Re}} 
\def\11{{\rm 1~\hspace{-1.4ex}l} }
\def\R{\mathbb R}
\def\C{\mathbb C}
\def\Z{\mathbb Z}
\def\N{\mathbb N}

\def\T{\mathbb T}

\begin{document}
\title[NLS and invariant measures]{Invariant measures for the defocusing NLS }
\author{N. Tzvetkov}
\address{D\'epartement de Math\'ematiques, Universit\'e Lille I, 59 655 Villeneuve d'Ascq Cedex, France}
\email{nikolay.tzvetkov@math.univ-lille1.fr}
\begin{abstract} 
We prove the existence and the invariance of a Gibbs measure associated to the defocusing sub-quintic 
Nonlinear Schr\"odinger equations on the disc of the plane $\R^2$. 
We also prove an estimate giving some intuition to what may happen in $3$ dimensions.
\\[0.4cm]
{\large R\'esum\'e.} On d\'emontre l'existence et l'invariance d'une mesure de Gibbs par 
le flot de l'\'equation de  Schr\"odinger non lin\'eaire pos\'ee sur le disque du plan $\R^2$. 
On d\'emontre \'egalement une estim\'ee qui donne une id\'ee de ce qui pourrait arriver en dimension $3$. 
\end{abstract}

\maketitle
\tableofcontents
%
\section{Introduction}
In \cite{Tz}, we constructed and proved the invariance of a Gibbs measure associated to the sub-cubic, 
focusing or defocusing Nonlinear Schr\"odinger equation (NLS) on the disc of the plane $\R^2$. 
For focusing non-linear interactions the cubic threshold is critical for the argument in \cite{Tz} 
because of measure existence obstructions.
The main goal of this paper is to show that, in the case of defocusing nonlinearities, one can 
extend the result of \cite{Tz} to the case of sub-quintic nonlinearities.
Thus we will be able to treat the relevant for the Physics case of cubic defocusing NLS. The argument 
presented here requires some significant elaborations with respect 
to \cite{Tz} both in the measure existence analysis and the Cauchy problem issues.
The main facts, proved in \cite{Tz} which will be used here without proof are
some properties of the Bessel functions and their zeros and the bilinear
Strichartz estimates of Proposition~\ref{str1} and  Proposition~\ref{str2} below.
\subsection{Presentation of the equation}
Let $V:\C\longrightarrow\R$ be a $C^\infty(\C)$ function. 
We suppose that $V$ is gauge invariant 
which means that there exists a smooth function $G:\R\longrightarrow\R$ such that $V(z)=G(|z|^2)$.
Set $F=\bar{\partial}V$, i.e. $F(z)=G'(|z|^2)z$. Consider the NLS
\begin{equation}\label{1}
(i\partial_t+\Delta)u-F(u)=0,
\end{equation}
where $u:\R\times\Theta\rightarrow \C$ is a complex valued function defined on the product 
of the real line (corresponding to the time variable) and 
$\Theta$, the unit disc of $\R^2$ corresponding to the spatial variable. More precisely
$$
\Theta=\big((x_1,x_2)\in\R^2\, :\, x_1^2+x_2^2<1\big).
$$
In this paper, we consider (\ref{1}) subject to Dirichlet boundary conditions
$u|_{\R\times \partial\Theta}=0$. It is likely that Neumann boundary
conditions are in the scope of applicability of our methods too.

We suppose that 
\begin{equation}\label{rast}
\exists\, \alpha\in ]0,4[\,:\,\forall\, (k_1,k_2),\,
\exists\, C>0\, :\, \forall\, z\in \C,\,
\big|\partial^{k_1}\bar{\partial}^{k_2}V(z)\big|\leq 
C(1+|z|)^{2+\alpha-k_1-k_2}\, .
\end{equation}
The number $\alpha$ involved in (\ref{rast}) measures the ``degree'' of the
nonlinearity. In this paper we will also suppose the defocusing assumption
\begin{equation}\label{defocus}
\exists\,\, C>0,\,\,\exists\,\, \beta\in [2,4[\,\, :\,\, 
\forall\, z\in\C,\,V(z)\geq -C(1+|z|)^{\beta}.
\end{equation}
A typical example for $V$ is 
$$
V(z)=\frac{2}{\alpha+2}\,\,(1+|z|^2)^{\frac{\alpha+2}{2}}
$$
with corresponding 
$$
F(z)=(1+|z|^2)^{\frac{\alpha}{2}}z\,.
$$
In the case $\alpha =2$ one can take $V(z)=\frac{1}{2}|z|^{4}$ which leads to a cubic defocusing 
nonlinearity $F(u)=|u|^{2}u$. Observe that $V(z)=-\frac{1}{2}|z|^{4}$, which
is the potential of the cubic 
focusing nonlinearity $F(u)=-|u|^{2}u$, does not satisfy assumption (\ref{defocus}).

We restrict our consideration only to radial solutions, i.e. we shall suppose that $u=u(t,r)$, where
$$
x_1=r\cos\phi,\quad x_2=r\sin\phi,\quad 0\leq r<1,\quad \phi\in[0,2\pi].
$$
Our goal here is to construct a Gibbs type measure, on a suitable phase space,
associated to the radial solutions 
of (\ref{1}) which is invariant under the (well-defined) global flow of (\ref{1}).
\subsection{Bessel expansion and formal Hamiltonian form}
Since we deal with radial solutions of (\ref{1}), it is natural to use Bessel function expansions.
Denote by $J_{0}(x)$ the zero order Bessel function. We have that (see \cite{Tz} and the 
references therein) $J_{0}(0)=1$ and $J(x)$ decays as 
$x^{-1/2}$ when $x\rightarrow\infty$. More precisely
$$
J_{0}(x)= \sqrt{\frac{2}{\pi}}\,\,\frac{\cos(x-\pi/4)}{\sqrt{x}}+{\mathcal O}(x^{-3/2}).
$$
Let $0<z_1<z_2<\cdots$ be the (simple) zeroes of $J_0(x)$. Then (see e.g. \cite{Tz}) $z_n\sim n$ as 
$n\rightarrow\infty$.
Each $L^2$ radial function may be expanded with respect to the Dirichlet bases formed by
$J_0(z_n r)$, $n=1,2,3,\cdots$. The functions $J_0(z_n r)$ are eigenfunctions of $-\Delta$ 
with eigenvalues $z_n^2$.
Define $e_n:\Theta\rightarrow\R$ by
$$
e_n\equiv e_n(r)=\|J_0(z_n\cdot)\|_{L^2(\Theta)}^{-1}J_0(z_nr)\,.
$$
We have (see \cite{Tz}) that $\|J_0(z_n\cdot)\|_{L^2(\Theta)}\sim n^{-1/2}$ as $n\rightarrow\infty$.
Therefore $\|e_n\|_{L^2(\Theta)}=1$ but $\|e_n\|_{L^{\infty}(\Theta)}\sim n^{1/2}$ as $n\rightarrow\infty$.
Hence we observe a significant difference between the disc and the flat torus $\T^2$, where the sup norm 
of the eigenfunctions can not grow so fast.

Let us fix from now on a real number $s$ such that
\begin{equation}\label{s}
\max\big(\frac{1}{3},1-\frac{2}{\alpha},1-\frac{2}{\beta}\big)<s<\frac{1}{2}
\end{equation}
(recall that $\alpha,\beta<4$ and thus a proper choice of the index $s$ is indeed possible).
Set $e_{n,s}=z_{n}^{-s}e_{n}$ ($H^s$ normalization) and if
$$
u(t)=\sum_{n=1}^{\infty}c_{n}(t)e_{n,s}
$$
then we need to analyze the equation
\begin{equation}\label{2}
iz_{n}^{-s}\dot{c_n}(t)-z_n^2\,z_{n}^{-s} c_n(t)-
\Pi_{n}\Big(
F\Big(
\sum_{m=1}^{\infty}c_m(t)\, e_{m,s}\Big)
\Big)=0,\quad n=1,2,\cdots
\end{equation}
where $\Pi_n$ is the projection on the mode $e_n$, i.e.
$\Pi_n(f)=\langle f, e_n\rangle.$
Equation (\ref{2}) is a Hamiltonian PDE for $c\equiv (c_n)_{n\geq 1}$ with Hamiltonian
$$
H(c,\overline{c})=\sum_{n=1}^{\infty} z_n^{2-2s}\, |c_n|^2 +2\pi\int_{0}^1V\Big(\sum_{n=1}^{\infty}c_n\, 
e_{n,s}(r)\Big)rdr\, ,
$$
and a formal Hamiltonian form
$$
ic_t=J\frac{\delta H}{\delta \overline{c}},\quad i\overline{c}_{t}=-J\frac{\delta H}{\delta c}\, ,
$$
where $J={\rm diag}(z_n^{2s})_{n\geq 1}$ is the
map inducing the symplectic form in the coordinates $(c,\overline{c})$. 
Thus the quantity $H(c,\overline{c})$ is, at least formally, conserved by the flow.
In fact we will need to use the energy conservation only for finite dimensional (Hamiltonian) 
approximations of (\ref{2}). 
Let us also observe that the $L^2$ norm of $u(t)$ expressed in terms of $c$ as
$$
\|c\|^2\equiv\sum_{n=1}^{\infty}z_{n}^{-2s}|c_n|^2
$$
is also conserved by the flow. Following Lebowitz-Rose-Speer \cite{LRS}, we will construct a 
{\bf renormalization} of the formal measure $\chi(\|c\|)\exp(-H(c,\overline{c}))dc\,d\bar{c}$
($\chi$ being a cut-off) which is invariant under the (well-defined) flow, living on a low regularity phase space
(for a finite dimensional Hamiltonian model the invariance would follow from the Liouville theorem 
for volume preservation by flows induced by divergence free vector fields). 
\subsection{The free measure}
Define the Sobolev spaces $H^\sigma_{rad}(\Theta)$, $\sigma\geq 0$ equipped with the norm
$$
\Big\|\sum_{n=1}^{\infty}c_n\, e_{n,s}\Big\|^{2}_{H^{\sigma}_{rad}(\Theta)}
\equiv\sum_{n=1}^{\infty}z_n^{2(\sigma-s)}|c_n|^{2}\,.
$$
The Sobolev spaces $H^\sigma_{rad}(\Theta)$ are related to the domains of $\sigma/2$ powers of the Dirichlet
Laplacian. In several places in the sequel, we shall denote $\|\cdot\|_{H^{\sigma}_{rad}(\Theta)}$
simply by $\|\cdot\|_{H^{\sigma}(\Theta)}$.
We can identify $l^2(\N;\C)$ with $H^s_{rad}(\Theta)$ via the map
$$
c\equiv (c_{n})_{n\geq 1}\longmapsto \sum_{n=1}^{\infty}c_n\, e_{n,s}\,.
$$
Consider the free Hamiltonian
$$
H_0(c,\overline{c})=\sum_{n=1}^{\infty} z_n^{2-2s}\, |c_n|^2 
$$
and the measure
$$
\frac{``\exp(-H_0(c,\overline{c}))dcd\bar{c}\,''}{\int\exp(-H_0(c,\overline{c}))dcd\bar{c}}=
\prod_{n=1}^{\infty}
\frac{e^{-z_n^{2-2s}|c_n|^{2}}dc_nd\bar{c}_{n}}
{\int_{\C}e^{-z_n^{2-2s}|c_n|^{2}}dc_nd\bar{c}_{n}}\equiv d\mu(c)\,.
$$
Denote by ${\mathcal B}$ the Borel sigma algebra of $H^s_{rad}(\Theta)$.
The measure $d\mu$ is first defined on cylindrical sets (see \cite{Tz}) in the natural way and
since for $s<1/2$,
$$
\sum_{n=1}^{\infty}z_{n}^{2s-2}<\infty
$$
we obtain that $d\mu$ is countably additive on the cylindrical sets and thus 
may be defined as a probability measure on $(H^s_{rad}(\Theta),{\mathcal B})$ via the map considered above. 
Let us recall that $A\subset H^s_{rad}(\Theta)$ is called cylindrical if there exist an integer $N$ and 
a Borel set of $V$ of $\C^N$ so that
$$
A=\Big\{u\in H^s_{rad}(\Theta)\, :\, \big( (u,e_{1,s}),\dots,(u,e_{N,s})\big)\in V\Big\}.
$$
In addition, the minimal sigma algebra on $H^s_{rad}(\Theta)$ containing the cylindrical sets is 
${\mathcal B}$.
\\

The measure $d\mu$ may also equivalently be defined as the distribution of the $H^s_{rad}(\Theta)$ 
valued random variable
\begin{equation*}
\varphi(\omega,r)=
\sum_{n= 1}^{\infty}\frac{g_n(\omega)}{z_n^{1-s}}e_{n,s}(r)=
\sum_{n= 1}^{\infty}\frac{g_n(\omega)}{z_n}e_{n}(r)\, ,
\end{equation*}
where $g_n(\omega)$ is a sequence of centered, normalised, independent identically distributed (i.i.d.)
complex Gaussian random variables, defined in a probability space
$(\Omega,{\mathcal F},p)$.
By normalised, we mean that 
$$
g_{n}(\omega)=\frac{1}{\sqrt{2}}\big(h_{n}(\omega)+i\, l_{n}(\omega)\big),
$$
where $h_n,l_n\in{\mathcal N}(0,1)$ are standard independent real gaussian variables.
Indeed, if we consider the sequence $(\varphi_N(\omega,r))_{N\in\N}$ defined by
\begin{equation}\label{3}
\varphi_N(\omega,r)=\sum_{n= 1}^{N}\frac{g_n(\omega)}{z_n}e_{n}(r)
\end{equation} 
then using that $s<1/2$ we obtain that $(\varphi_N(\omega,r))_{N\in\N}$ is a Cauchy sequence in 
$L^2(\Omega;H^s_{rad}(\Theta))$ and $\varphi(\omega,r)$ is, by definition, the limit of this sequence. 
Thus the map which to $\omega\in\Omega$ associates $\varphi(\omega,r)$ is measurable from
$(\Omega,{\mathcal F})$ to $(H^s_{rad}(\Theta),{\mathcal B})$.
Therefore $\varphi(\omega,r)$ 
may be seen as a $H^s_{rad}(\Theta)$ valued random variable and for every Borel set $A\in {\mathcal B}$,
$$
\mu(A)=p(\omega\,:\, \varphi(\omega,r)\in A).
$$
Moreover, if $f:H^s_{rad}(\Theta)\rightarrow \R$ is a measurable function then $f$ is integrable if and 
only if the real random variable $f\circ \varphi\,:\, \Omega\rightarrow \R$ is integrable and
$$
\int_{H^s_{rad}(\Theta)}f(u)d\mu(u)=\int_{\Omega}f(\varphi(\omega,\cdot))dp(\omega)\,.
$$ 
\subsection{Measure existence}
Following the basic idea one may expect that the measure (Gibbs measure)
\begin{equation}\label{mesdef}
d\rho(u)\equiv\chi\big(\|u\|_{L^2(\Theta)}\big)\exp\Big(-\int_{\Theta}V(u)\Big)d\mu(u)
\end{equation}
is invariant under the flow of (\ref{1}).
In (\ref{mesdef}), 
$$
\chi:\R\longrightarrow [0,\infty[
$$ 
is a non-negative continuous function with compact support.
In (\ref{mesdef}), $\exp\big(-\int_{\Theta}V(u)\big)$ is the contribution of
the nonlinearity of (\ref{1}) to the Hamiltonian, while the free Hamiltonian
(coming from the linear part of (\ref{1})) is incorporated in $d\mu(u)$.
One may wish to see $d\rho(u)$ as the image measure on $H^s_{rad}(\Theta)$ under the map
$$
\omega\longmapsto \sum_{n= 1}^{\infty}\frac{g_n(\omega)}{z_n}e_{n}(r)
$$
of the measure
$$
\chi\big(\|\varphi(\omega,\cdot)\|_{L^2(\Theta)}\big)
\exp\Big(-2\pi\int_{0}^{1}V\big(\varphi(\omega,r)\big)rdr \Big)dp(\omega)\,.
$$
A first problem (in order to ensure that $\rho$ is not trivial) is whether 
$\int_{\Theta}V(u)$ is finite $\mu$ almost surely (a.s.). Let us notice that an appeal to the 
(\ref{rast}) and the
Sobolev inequality gives
\begin{equation}\label{sobolev}
\Big|
\int_{\Theta}V(u)
\Big|
\leq
C\big(1+\|u\|_{L^{\alpha+2}(\Theta)}^{\alpha+2}\big)
\leq
C\big(1+\|u\|^{\alpha+2}_{H^{\sigma}_{rad}(\Theta)}\big),
\end{equation}
provided $\sigma\geq 2(\frac{1}{2}-\frac{1}{2+\alpha})=\frac{\alpha}{2+\alpha}$.
For $\alpha\geq 2$ (a case excluded in \cite{Tz}), inequality (\ref{sobolev}) does 
not suffice to conclude that $\int_{\Theta}V(u)$ is finite $\mu$ a.s.
Indeed, for $\alpha\geq 2$ one has $\sigma\geq \frac{1}{2}$ and,
using for instance the Fernique integrability theorem, one may show that 
$\|u\|_{H^{\sigma}_{rad}(\Theta)}=\infty$, $\mu$ a.s.
We can however resolve this problem by using a probabilistic argument 
(which ``improves'' on the Sobolev inequality).
Let us also mention the recent work \cite{AT}, where one studies $L^p$ 
properties of Gaussian random series with a particular attention to radial functions. 
Here is a precise statement.
\begin{theoreme}\label{thm1}
We have that $\int_{\Theta}V(u)\in L^{1}(d\mu(u))$
(in particular $\int_{\Theta}V(u)$ is $\mu$ a.s. finite).
\end{theoreme}
Essentially, the assertion of Theorem~\ref{thm1} follows from the considerations in \cite{AT}. 
We will however give below a proof of Theorem~\ref{thm1} using an argument slightly different from \cite{AT}.
\subsection{Finite dimensional approximations}
Let $E_N$ be the finite dimensional complex vector space spanned by $(e_{n})_{n=1}^{N}$.
We consider $E_N$ as a measured space with the measure induced by $\C^n$ under
the map from $\C^N$ to $E_N$ defined by
$$
(c_1,\cdots,c_{N})\longmapsto \sum_{n=1}^{N}c_{n}e_{n,s}\,.
$$
Following Zhidkov (cf. \cite{Zh} and the references therein), we consider the finite dimensional 
projection (an ODE) of (\ref{1})
\begin{equation}\label{N}
(i\partial_t+\Delta)u-S_{N}(F(u))=0,\quad u|_{t=0}\in E_N,
\end{equation}
where $S_N$ is the projection on $E_N$. 
Notice that $S_{N}(F(u))$ is well-defined for $u\in E_N$ since $E_N\subset
C^{\infty}(\overline{\Theta})$.
The equation (\ref{N}) is a Hamiltonian ODE for $u\in E_N$ with Hamiltonian
$$
H_{N}(u,\bar{u})=\int_{\Theta}|\nabla u|^{2}+\int_{\Theta}V(u),\quad u\in E_N\,.
$$
Thus $H_{N}(u,\bar{u})$ is conserved by the flow of (\ref{N}). One may directly check this by 
multiplying (\ref{N}) with $\bar{u}_{t}\in E_N$ and integrating over $\Theta$
(observe that the boundary terms in the integration by parts disappear).
Multiplying (\ref{N}) by $\bar{u}$ and integrating over $\Theta$, we see that the $L^2(\Theta)$ norm 
is also preserved by the flow of (\ref{N}) and thus (\ref{N}) has a well-defined global dynamics.
Denote by $\Phi_{N}(t):E_N\rightarrow E_N$, $t\in\R$ the flow of (\ref{N}).
Let $\mu_N$ be the distribution of the $E_N$ valued random variable
$\varphi_{N}(\omega,r)$ defined by (\ref{3}). Set 
$$
d\rho_N(u)\equiv
\chi\big(\|u\|_{L^2(\Theta)}\big)
\exp\big(-\int_{\Theta}V(u)\big)d\mu_N(u).
$$
One may see $\rho_{N}$ as the image measure on $E_N$ under the map
$\omega\mapsto \varphi_{N}(\omega,r)$ of the measure
$$
\chi\big(\|\varphi_{N}(\omega,\cdot)\|_{L^2(\Theta)}\big)
\exp\Big(-2\pi\int_{0}^{1}V\big(\varphi_{N}(\omega,r)\big)rdr\Big)dp(\omega)\,.
$$
From the Liouville theorem for divergence free vector fields, the measure $\rho_{N}$ is 
invariant under $\Phi_{N}(t)$. Indeed, if we write the solution of (\ref{N}) as
$$
u(t)=\sum_{n=1}^{N}c_{n}(t)e_{n,s}\,,\quad c_{n}(t)\in \C
$$
then in the coordinates $c_{n}$, the equation (\ref{N}) can be written as
\begin{equation}\label{4bis}
iz_{n}^{-s}\dot{c_n}(t)-z_n^2\,z_{n}^{-s} c_n(t)-
\int_{\Theta}S_{N}(F(u(t)))\overline{e_{n}}=0,\quad 1\leq n\leq N.
\end{equation}
Equation (\ref{4bis}) in turn can be written in a Hamiltonian format as follows
$$
\partial_{t}c_n=-iz_{n}^{2s}\frac{\partial H}{\partial \overline{c_n}},
\quad 
\partial_{t}\overline{c_n}=iz_{n}^{2s}\frac{\partial H}{\partial c_n},\quad 1\leq n\leq N,
$$
with
$$
H(c,\overline{c})=\sum_{n= 1}^{N} z_n^{2-2s}\, |c_n|^2 + 2\pi\int_{0}^1
V\Big(\sum_{n=1}^{N}c_n\, e_{n,s}(r)\Big)rdr\, ,\quad c=(c_1,\cdots,c_N).
$$
Since
$$
\sum_{n=1}^{N}
\Big(
\frac{\partial}{\partial c_n}\big(-iz_{n}^{2s}\frac{\partial H}{\partial \overline{c_n}}\big)
+
\frac{\partial}{\partial\overline{c_n}}
\big(iz_{n}^{2s}\frac{\partial H}{\partial c_n}\big)
\Big)=0,
$$
we can apply the Liouville theorem for divergence free vector fields to
conclude that the measure $dcd\overline{c}$ is invariant under the flow of
(\ref{4bis}).
On the other hand the quantities $H(c,\overline{c})$ and 
$$
\|c\|^{2}\equiv\sum_{n=1}^{N}z_{n}^{-2s}|c_n|^{2}
$$
are conserved under the flow of (\ref{4bis}).
Moreover, by definition if $A$ is a Borel set of $E_{N}$ defined by
$$
A=\Big\{
u\in E_{N}\,:\, u=\sum_{n=1}^{N}c_{n}e_{n,s},\quad (c_1,\cdots,c_N)\in A_1
\Big\},
$$
where $A_1$ is a Borel set of $\C^N$, then
$$
\rho_{N}(A)=\kappa_{N}\int_{A_1}e^{-H(c,\overline{c})}\chi(\|c\|)dcd\overline{c},
$$
with
$$
\kappa_{N}=\pi^{-N}\Big(\prod_{1\leq n\leq N}z_{n}^{2-2s}\Big).
$$
Therefore the measure $\rho_{N}$ is invariant under $\Phi_{N}(t)$,
thanks to the invariance of $dcd\overline{c}$ and the $\Phi_{N}(t)$ 
conservations of of $H(c,\overline{c})$ and $\chi(\|c\|)$.
Let us also observe that if we write (\ref{4bis}) in terms of $(\Re(c_n),\Im(c_n))$ then
we still obtain a Hamiltonian ODE and one may show the invariance of $\rho_N$
under (\ref{N}) by analyzing that ODE.
\\

One may extend $\rho_N$ to a measure $\tilde{\rho}_{N}$ on $H^s_{rad}(\Theta)$.
If $U$ is a $\rho$ measurable set then $\tilde{\rho}_{N}(U)\equiv \rho_{N}(U\cap E_{N})$.
A similar definition may be given for $\mu_N$.
The measure $\tilde{\rho}_{N}$ is well-defined since for $U\in {\mathcal B}$ one has that
$U\cap E_{N}$ is a Borel set of $E_{N}$. Indeed, this property is clear for $U$ a cylindrical set and
then we extend it to ${\mathcal B}$ by the key property of the cylindrical sets.
Observe that for $U$, a $\rho$ measurable set, one has
$$
\tilde{\rho}_{N}(U)=\int_{U_{N}}\chi\big(\|S_N(u)\|_{L^2(\Theta)}\big)
\exp\big(-\int_{\Theta}V(S_{N}(u))\big)d\mu(u),
$$
where 
$$
U_{N}=\big\{u\in H^s_{rad}(\Theta)\,:\, S_{N}(u)\in U\big\}.
$$
The following properties relating $\rho$ and $\rho_N$ will be useful in our analysis concerning the long 
time dynamics of (\ref{1}).
\begin{theoreme}\label{thm2}
One has that for every $p\in [1,\infty[$,
$$
\chi\big(\|u\|_{L^2(\Theta)}\big)
\exp\big(-\int_{\Theta}V(u)\big)\in L^p(d\mu(u)).
$$
In addition, if we fix $\sigma\in [s,1/2[$ then for every $U$ an open set of $H^{\sigma}_{rad}(\Theta)$ 
one has
\begin{equation}\label{parvo}
\rho(U)\leq\liminf_{N\rightarrow\infty}\tilde{\rho}_{N}(U)\,\, 
(= \liminf_{N\rightarrow\infty}\rho_{N}(U\cap E_{N})) \,.
\end{equation}
Moreover if $F$ is a closed set of $H^{\sigma}_{rad}(\Theta)$ then
\begin{equation}\label{vtoro}
\rho(F)\geq\limsup_{N\rightarrow\infty}\tilde{\rho}_{N}(F)\,\, 
(=\limsup_{N\rightarrow\infty}\rho_{N}(F\cap E_{N})).
\end{equation}
\end{theoreme}
The proof of Theorem~\ref{thm2} is slightly more delicate than an analogous result used in \cite{Tz}.
In contrast with \cite{Tz} we can not exploit that $\int_{\Theta}V(S_{N}u)$ converges $\mu$ 
a.s. to $\int_{\Theta}V(u)$. 
In \cite{Tz} we deal with sub-quartic growth of $V$ and by the Sobolev embedding 
we can get directly the needed $\mu$ a.s. convergence. Here we will need to use a different argument.
\subsection{Statement of the main result}
With Theorem~\ref{thm1} and Theorem~\ref{thm2} in hand we can prove our main result.
\begin{theoreme}\label{thm3}
The measure  $\rho$ is invariant under the well-defined $\rho$ a.s. global in time flow of
the NLS (\ref{1}), posed on the disc. More precisely :
\begin{itemize}
\item
There exists a $\rho$ measurable set $\Sigma$ of full $\rho$ measure such that for every $u_0\in\Sigma$ 
the NLS (\ref{1}), posed on the disc,
with initial data data $u|_{t=0}=u_0$ has a unique (in a suitable sense) 
global in time solution $u\in C(\R;H^s_{rad}(\Theta))$. 
In addition, for every $t\in\R$, $u(t)\in\Sigma$ and the map $u_0\mapsto u(t)$
is $\rho$ measurable.
\item
For every $A\subset\Sigma$, a $\rho$ measurable set, for every $t\in\R$,
$
\rho(A)=\rho(\Phi(t)(A)),
$
where $\Phi(t)$ denotes the flow defined in the previous point.
\end{itemize}
\end{theoreme}
The uniqueness statement of Theorem~\ref{thm3} is in the sense of a uniqueness for the integral equation
(\ref{venda}) in a suitable space continuously embedded in the space of continuous $H^s_{rad}(\Theta)$
valued functions. Another possibility is to impose zero boundary conditions on $\R\times \partial\Theta$ 
and then relate the solutions of (\ref{1}) to the solutions of (\ref{venda}) 
(see also Remark~\ref{zabelejka} below). 

As a consequence of Theorem~\ref{thm3} one may apply the Poincar\'e recurrence theorem to the flow $\Phi$.
For previous works proving the invariance of Gibbs measures under the flow of NLS we refer to 
\cite{Bo1,Bo2,Zh}. 
In all these works one considers periodic boundary conditions, i.e. the spatial domain is the flat torus.
We also refer to \cite{KS}, for a construction of invariant measures, supported by $H^2$,
for the defocusing NLS.

Let us also remark that the result of Theorem~\ref{thm3} implies that the
sub-quintic defocusing NLS is almost surely globally well-posed for data
$\varphi(\omega,r)$ defined by
$$
\varphi(\omega,r)=  \sum_{n= 1}^{\infty}\frac{g_n(\omega)}{z_n}e_{n}(r)\,.
$$
Because of the low regularity of $\varphi$ for typical $\omega$'s such a
result seems to be difficult to achieve by the present deterministic methods
for global well-posedness of NLS.
\subsection{Structure of the paper and notation}
Let us briefly describe the organization of the rest of the paper. 
In the next section, we prove Theorem~\ref{thm1}. Section~3, is devoted to the proof of Theorem~\ref{thm2}.
In Section~4, we recall the definition of the Bourgain spaces and we state two bilinear 
Strichartz estimates which are the main tool in the study of the local Cauchy problem. In Section~5, 
we prove nonlinear estimates in Bourgain spaces. Section~6 is devoted to the local well-posedness analysis.
In Section~7, we establish the crucial control on the dynamics of (\ref{N}). In section~8, we construct 
the set $\Sigma$ involved in the statement of Theorem~\ref{thm3}. 
In Section~9, we prove the invariance of the measure. In the last section, we prove several bounds for the 
3d NLS with random data. 
\\

In this paper, we assume that the set of the natural numbers $\N$ is $\{1,2,3,\cdots\}$.
We call dyadic integers the non-negative powers of $2$, i.e. $1,2,4,8$ etc.
\subsection{Acknowledgements.} 
I am very grateful to Nicolas Burq for several useful discussions on the problem and for
pointing out an error in a previous version of this text.
It is a pleasure to thank A.~Ayache and H.~Queff\'elec for useful 
discussions on random series. 
I am also indebted to N.~Burq and P.~G\'erard since this work (as well as \cite{Tz}) benefited from 
our collaborations on NLS on compact manifolds.
I thank the referee for pointing out several imprecisions in a previous version
of the paper.
\section{Proof of Theorem~\ref{thm1} (measure existence)}
\subsection{Large deviation estimates}
\begin{lemme}\label{lem1}
Let $(g_n(\omega))_{n\in\N}$ be a sequence of normalized i.i.d. complex Gaussian random
variables defined in a probability space $(\Omega,{\mathcal F},p)$.
There exists $\beta>0$ such that for every $\lambda >0$, every sequence $(c_n)\in l^2(\N;\C)$
of complex numbers,
\begin{equation*}
p\Big(\omega\,:\,\big|\sum_{n=1}^{\infty} c_n g_n(\omega) \big|>\lambda\Big)\leq 
4\,e^{-\frac{\beta\lambda^2}{\sum_{n}|c_n|^2}}
\end{equation*}
(the right hand-side being defined as zero if $(c_n)_{n\in\N}$ is identically zero).
\end{lemme}
\begin{proof}
By separating the real and the imaginary parts, we can assume that $g_n$ are real valued independent
standard gaussians and $c_n$ are real constants. The bound we need to prove
is thus
\begin{equation}\label{real}
\exists\,\beta>0\,:\,
\forall\, (c_n)\in l^2(\N;\R),\, \forall\, \lambda>0,\quad
p\Big(\omega\,:\,\big|\sum_{n=1}^{\infty} c_n g_n(\omega) \big|>\lambda\Big)\leq 
2\,e^{-\frac{\beta\lambda^2}{\sum_{n}c_n^2}}\, .
\end{equation}
We may of course assume that the sequence $(c_n)_{n\in\N}$ is not identically zero.
For $t>0$ to be determined later, using the independence, we obtain that
\begin{equation*}
\int_{\Omega}\, \exp\Big(t\sum_{n= 1}^{\infty}c_n g_n(\omega)\Big)dp(\omega)= 
\exp\Big((t^2/2)\sum_{n=1}^{\infty}c_n^2\Big)\, .
\end{equation*}
Therefore
$$
\exp\Big((t^2/2)\sum_{n=1}^{\infty}c_n^2\Big)\geq \exp(t\lambda)\,\,\,  
p\,\Big(\omega\,:\,\sum_{n=1}^{\infty} c_n g_n(\omega)>\lambda\Big)
$$
and thus
$$
p\,(\omega\,:\,\sum_{n=1}^{\infty} c_n g_n(\omega)>\lambda)\leq \exp\Big((t^2/2)
\sum_{n=1}^{\infty}c_n^2\Big)\,\,\,
e^{-t\lambda}\, .
$$
For $a>0$, $b>0$ the minimum of $f(t)=at^2-bt$ is $-b^2/4a$ and this minimum is
attained in the positive number $t=b/(2a)$. It is thus natural to choose the positive number $t$ as 
$$
t\equiv \lambda/\big(\sum_{n=1}^{\infty}c_n^2\big)
$$
which leads to
$$
p\,(\omega\,:\,\sum_{n=1}^{\infty} c_n g_n(\omega)>\lambda)\leq
\exp\Big(-\frac{\lambda^2}{2\sum_{n}c_n^2}\Big)\, .
$$
In the same way (replacing $c_n$ by $-c_n$), we can show that
$$
p\,(\omega\,:\,\sum_{n=1}^{\infty} c_n g_n(\omega)<-\lambda)\leq
\exp\Big(-\frac{\lambda^2}{2\sum_{n}c_n^2}\Big)
$$
which shows that (\ref{real}) holds with $\beta=1/2$.
This completes the proof of Lemma~\ref{lem1}.
\end{proof}
We next state the following consequence of Lemma~\ref{lem1}.
\begin{lemme}\label{armen}
Let $(g_n(\omega))_{n\in\N}$ be a sequence of normalized i.i.d. complex Gaussian random
variables defined in a probability space $(\Omega,{\mathcal F},p)$.
Then there exist positive numbers $c_1,c_2$ such that for every non empty
finite set of indexes $\Lambda \subset \N $, every $\lambda >0$, 
$$
p\Big(\omega\in \Omega\, : \, \sum_{n\in \Lambda} |g_n(\omega)|^2 >\lambda\Big)\leq
e^{c_1|\Lambda|-c_2\lambda}
\, ,
$$
where $|\Lambda|$ denotes the cardinality of $\Lambda$.
\end{lemme}
\begin{proof}
A proof of this lemma is given in \cite[Lemma~3.4]{Tz}.
Here we propose a different proof based on Lemma~\ref{lem1}.
The interest of this proof is that the argument might be useful in more
general situations.
Again, we can suppose that $g_n$ are real valued standard gaussians.
A simple geometric observation shows that there exists $c_1>0$ (independent of $|\Lambda|$)
and a set ${\mathcal A}$ of the unit ball of $\R^{|\Lambda|}$ of cardinality bounded
by $e^{c_1|\Lambda|}$ such that almost surely in $\omega$,
$$
\frac{1}{2}\, \Big(\sum_{n\in \Lambda}\, |g_n(\omega)|^2\Big)^{1/2}
\leq 
\sup_{c\in {\mathcal A}}\, 
\big|\sum_{n\in\Lambda} c_n g_{n}(\omega)  \big|\, 
$$ 
($c=(c_n)_{n\in\Lambda}$ with $\sum_{n}|c_n|^2=1$).
Therefore
$$
\{\omega\,: \,\sum_{n\in \Lambda} |g_n(\omega)|^2 >\lambda\}\subset\,\bigcup_{c\in {\mathcal A}}
\,\{\omega\,: \,|\sum_{n\in\Lambda} c_n g_{n}(\omega)|\geq \frac{\sqrt{\lambda}}{2}\}\, .
$$
Consequently, using Lemma~\ref{lem1}, we obtain that there exists $c_2>0$,
independent of $\Lambda$, such that for every $\lambda>0$,
$$
p\,(\omega\,: \,\sum_{n\in \Lambda} |g_n(\omega)|^2 >\lambda)
\leq 
|{\mathcal A}|\, 4e^{-c_2\lambda}\leq 4\, e^{c_1|\Lambda|-c_2\lambda}
< e^{(c_1+2)|\Lambda|-c_2\lambda}\, .
$$
This completes the proof of Lemma~\ref{armen}.
\end{proof}
\subsection{Proof of Theorem~\ref{thm1}}
Theorem~\ref{thm1} follows from the following statement.
\begin{lemme}\label{rudin}
The sequence $\int_{\Theta}V(S_{N}(u))$ converges to $\int_{\Theta}V(u)$ in $L^1(d\mu)$.
\end{lemme}
\begin{proof}
Let us first show that $(\int_{\Theta}V(S_{N}(u)))_{N\in\N}$ is a Cauchy sequence in $L^1(d\mu)$. 
From the Sobolev embedding, we have that for a fixed $N$ the map from $H^s_{rad}(\Theta)$ to $\R$ defined by
$u\mapsto \int_{\Theta}V(S_{N}(u))$ is continuous and thus measurable. 
Write, for $N<M$, using (\ref{rast})
\begin{multline*}
\Big\|\int_{\Theta}V(S_{N}(u))-\int_{\Theta}V(S_{M}(u))\Big\|_{L^1(H^s_{rad};{\mathcal B},d\mu(u))}
\\
\leq C\Big\|\int_{\Theta}|S_{N}(u)-S_{M}(u)|(1+|S_{N}(u)|^{\alpha+1}+|S_{M}(u)|^{\alpha+1})
\Big\|_{L^1(H^s_{rad};{\mathcal B},d\mu(u))}\,.
\end{multline*}
Using the H\"older inequality, we get
\begin{multline*}
\Big|\int_{\Theta}|S_{N}(u)-S_{M}(u)|(1+|S_{N}(u)|^{\alpha+1}+|S_{M}(u)|^{\alpha+1})\Big|
\\
\leq \|S_{N}(u)-S_{M}(u)\|_{L^{\alpha+2}(\Theta)}\big(C+\|S_{N}(u)\|_{L^{\alpha+2}(\Theta)}^{\alpha+1}+
\|S_{M}(u)\|_{L^{\alpha+2}(\Theta)}^{\alpha+1}\big).
\end{multline*}
Another use of the H\"older inequality, this time with respect to $d\mu$ gives
\begin{multline*}
\Big\|\int_{\Theta}V(S_{N}(u))-\int_{\Theta}V(S_{M}(u))\Big\|_{L^1(d\mu(u))}\leq
C\big\|\|S_{N}(u)-S_{M}(u)\|_{L^{\alpha+2}(\Theta)}\big\|_{L^{\alpha+2}(d\mu(u))}
\\
\times\Big(1+\big\|\|S_{N}(u)\|_{L^{\alpha+2}(\Theta)}\big\|_{L^{\alpha+2}(d\mu(u))}^{\alpha+1}+
\big\|\|S_{M}(u)\|_{L^{\alpha+2}(\Theta)}\big\|_{L^{\alpha+2}(d\mu(u))}^{\alpha+1}\Big).
\end{multline*}
Thus
\begin{multline}\label{sarbia}
\Big\|\int_{\Theta}V(S_{N}(u))-\int_{\Theta}V(S_{M}(u))\Big\|_{L^1(d\mu(u))}\leq
C\|\varphi_{N}-\varphi_{M}\|_{L^{\alpha+2}(\Theta\times\Omega)}
\\
\times\Big(1+\big\|\varphi_{N}\|_{L^{\alpha+2}(\Theta\times\Omega)}^{\alpha+1}
+\big\|\varphi_{M}\|_{L^{\alpha+2}(\Theta\times\Omega)}^{\alpha+1}\Big),
\end{multline}
where $\varphi_{N}$ is defined by (\ref{3}).
Let us now prove that there exists $C>0$ such that for every $N$,
\begin{equation}\label{ravn}
\|\varphi_{N}\|_{L^{\alpha+2}(\Omega\times\Theta)}\leq C.
\end{equation}
Using Lemma~\ref{lem1} with $c_n=z_{n}^{-1}e_{n}(r)$, $1\leq n \leq N$ and the definition of the 
$L^{\alpha+2}$ norms by the aide of the distributional function, we obtain that for a fixed $r$
\begin{eqnarray*}
\|\varphi_{N}(\omega,r)\|_{L^{\alpha+2}(\Omega)}^{\alpha+2} & = &(\alpha+2)\int_{0}^{\infty}\lambda^{\alpha+1}
p\Big(\omega\,:\,\big|\varphi_{N}(\omega,r)\big|>\lambda\Big)d\lambda
\\
& \leq & 
C\int_{0}^{\infty}\lambda^{\alpha+1}
\exp\Big(-(\beta\lambda^2)/\big(\sum_{n=1}^{N}z_{n}^{-2}|e_{n}(r)|^{2}\big)\Big)d\lambda
\\
&= &
C\big(\int_{0}^{\infty}\lambda^{\alpha+1}e^{-\beta\lambda^2}d\lambda\big)
\big(\sum_{n= 1}^{N}z_{n}^{-2}|e_{n}(r)|^{2}\big)^{\frac{\alpha+2}{2}}\,.
\end{eqnarray*}
Therefore
$$
\|\varphi_{N}(\omega,r)\|_{L^{\alpha+2}(\Omega)}\leq 
C\big(\sum_{n=1}^{N}z_{n}^{-2}|e_{n}(r)|^{2}\big)^{\frac{1}{2}}\,.
$$
Squaring, taking the $L^{\frac{\alpha+2}{2}}(\Theta)$ norm and using the triangle inequality, we get
$$
\|\varphi_{N}\|_{L^{\alpha+2}(\Omega\times\Theta)}^{2}\leq 
\sum_{n=1}^{N}z_n^{-2}\|e_{n}\|_{L^{\alpha+2}(\Theta)}^{2}\,.
$$
On the other hand, it is shown in \cite{Tz} that for $\alpha<2$ one has that
$\|e_{n}\|_{L^{\alpha+2}(\Theta)}$ is uniformly bounded (with respect to $n$),
for $\alpha=2$, $\|e_{n}\|_{L^{\alpha+2}(\Theta)}\leq C\log(1+z_n)^{1/4}$
and for $\alpha>2$, $\|e_{n}\|_{L^{\alpha+2}(\Theta)}\leq Cz_n^{1/2-2/(\alpha+2)}$.
Since $z_n\sim n$, we obtain that there exists $C$ such that for every $N\in \N$,
$\|\varphi_{N}\|_{L^{\alpha+2}(\Omega\times\Theta)}\leq C.$
Therefore (\ref{ravn}) holds. 
\\

Similarly, we may obtain that
\begin{equation}\label{ravn2}
\|\varphi_{N}-\varphi_{M}\|_{L^{\alpha+2}(\Omega\times\Theta)}^{2}\leq 
\sum_{n=N+1}^{M}z_n^{-2}\|e_{n}\|_{L^{\alpha+2}(\Theta)}^{2}
\end{equation}
which tends to zero as $N\rightarrow\infty$ thanks to the bounds on the growth of 
$\|e_{n}\|_{L^{\alpha+2}(\Theta)}$.
Moreover, we have that
\begin{equation}\label{noel-pak}
\lim_{N\rightarrow\infty}\varphi_{N}=\varphi\quad{\rm in}\quad L^{\alpha+2}(\Theta\times\Omega)
\end{equation}
(we can identify the limit thanks to the $L^{2}(\Theta\times\Omega)$
convergence of $\varphi_{N}$ to $\varphi$ and the fact that
$L^{\alpha+2}(\Theta\times\Omega)$ convergence implies
$L^{2}(\Theta\times\Omega)$ convergence).

On the other hand thanks to (\ref{defocus}), we can write 
$
V(u)=V_1(u)+V_2(u),
$
where $V_1\geq 0$ and 
$
|V_2(u))|\leq C\big(1+|u|^{\beta}\big).
$
Thanks to the Sobolev embedding and (\ref{s}), we obtain that $\int_{\Theta}V_{2}(u)$ 
is continuous on $H^s_{rad}(\Theta)$.
Therefore the map $u\mapsto \int_{\Theta}V_{2}(u)$ is a $\mu$ measurable real valued function.
Let us next show that the map $u\mapsto \int_{\Theta}V_{1}(u)$ is $\mu$ measurable.
For that purpose, it is sufficient to show that the map 
\begin{equation}\label{ccc}
c\equiv (c_n)_{n\in\N}\longmapsto \int_{\Theta}V_{1}\Big(\sum_{n\in\N}c_{n}e_{n,s}\Big)
\end{equation}
is measurable from $l^2(\N)$ to $\R$. Indeed, we have that the map
$$
(c,r)\longmapsto \sum_{n\in\N}c_{n}e_{n,s}(r)
$$
is measurable from from $l^2(\N)\times\Theta$ to $\R$ since we can see 
$\sum_{n\in\N}c_{n}e_{n,s}(r)$ as the limit of $\sum_{n=1}^{N}c_{n}e_{n,s}(r)$
as $N\rightarrow\infty$ in $L^2( l^2(\N)\times \Theta)$ where $l^2(\N)$ is
equipped with the measure $d\mu(c)$ introduced in the
introduction. Therefore $V_1\Big(  \sum_{n\in\N}c_{n}e_{n,s}  \Big)$ is a
measurable map from $l^2(\N)\times\Theta$ to $\R$.
Since $V_1\geq 0$, using for instance the Fubini theorem,
we obtain that the map (\ref{ccc}) is indeed measurable.
This in turn implies the measurability of the map $u\mapsto \int_{\Theta}V(u)$.
Next, similarly to the proof of (\ref{sarbia}), we get 
\begin{multline*}
\Big\|\int_{\Theta}V(S_{N}(u))-\int_{\Theta}V(u)\Big\|_{L^1(d\mu(u))}\leq
C\|\varphi-\varphi_{N}\|_{L^{\alpha+2}(\Theta\times\Omega)}
\\
\times\Big(1+\big\|\varphi_{N}\|_{L^{\alpha+2}(\Theta\times\Omega)}^{\alpha+1}
+\big\|\varphi\|_{L^{\alpha+2}(\Theta\times\Omega)}^{\alpha+1}\Big)\,.
\end{multline*}
Therefore
$$
\lim_{N\rightarrow\infty}\Big\|\int_{\Theta}V(S_{N}(u))-
\int_{\Theta}V(u)\Big\|_{L^1(H^s_{rad};{\mathcal B},d\mu(u))}=0\,.
$$
This completes the proof of Lemma~\ref{rudin}.
\end{proof}
Using  Lemma~\ref{rudin}, we have that $\int_{\Theta}V(u)\in L^1(d\mu(u))$ and thus 
$\int_{\Theta}V(u)$ is finite $\mu$ a.s.
This proves that $d\rho$ is indeed a nontrivial measure.
This completes the proof of Theorem~\ref{thm1}.
\qed
\subsection{The necessity of the probabilistic argument}
In this section we make a slight digression by showing that for $\alpha\geq 2$ an argument based only on 
the Sobolev embedding may not conclude to the fact that $\int_{\Theta}V(u)$ is finite 
$\mu$ a.s. More precisely we know that for every $\sigma<1/2$, $\|u\|_{H^{\sigma}(\Theta)}$
is finite $\mu$ a.s. Therefore the deterministic inequality
\begin{equation}\label{wrong}
\exists\,\, \sigma<1/2,\,\,\exists\,\, C>0,\,\, \forall\,\, u\in H^{\sigma}_{rad}(\Theta),\quad
\|u\|_{L^{\alpha+2}(\Theta)}\leq C\|u\|_{H^{\sigma}_{rad}(\Theta)}
\end{equation}
would suffice to conclude that $\int_{\Theta}V(u)$ is finite 
$\mu$ a.s. We have however the following statement.
\begin{lemme}\label{scaling}
For $\alpha\geq 2$, estimate (\ref{wrong}) fails.
\end{lemme}
\begin{proof}
We shall give the proof for $\alpha=2$. The construction for $\alpha>2$ is similar.
Suppose that (\ref{wrong}) holds for some $\sigma<1/2$. Using the
Cauchy-Schwarz inequality, we obtain that there exists $\theta\in ]0,1/2]$ 
such that 
$$
\exists\, C>0\,:\, \forall\, u\in H^1_{rad}(\Theta),\quad
\|u\|_{H^{\sigma}(\Theta)}\leq C
\|u\|_{L^{2}(\Theta)}^{\frac{1}{2}+\theta}\|u\|_{H^{1}_{rad}(\Theta)}^{\frac{1}{2}-\theta}
$$
(observe that $H^1_{rad}(\Theta)$ may be seen as the completion of $C_{0}^{\infty}(\Theta)$
radial functions with respect to the $H^1(\Theta)$ norm).
Thus by applying (\ref{wrong}) to $H^1_{rad}(\Theta)$ functions, we obtain that 
\begin{equation}\label{wrong-bis}
\exists\, C>0\,:\,\forall\,\, u\in H^{1}_{rad}(\Theta),\quad \|u\|_{L^{4}(\Theta)}\leq C
\|u\|_{L^{2}(\Theta)}^{\frac{1}{2}+\theta}\|u\|_{H^{1}(\Theta)}^{\frac{1}{2}-\theta}\,.
\end{equation}
We now show that (\ref{wrong-bis}) fails. Let $v\in C_{0}^{\infty}(\Theta)$ be a radial bump function,
not identically zero. 
We can naturally see $v$ as a $C_{0}^{\infty}(\R^2)$ function. For $\lambda\geq 1$, we set
$$
v_{\lambda}(x_1,x_2)\equiv v(\lambda x_1,\lambda x_2)\,.
$$
Thus $v_{\lambda}\in C_{0}^{\infty}(\Theta)$ and $v_{\lambda}$ is still radial. 
We can therefore substitute $v_{\lambda}$ in (\ref{wrong-bis}) and obtain a contradiction in the limit
$\lambda\rightarrow \infty$. More precisely, one may directly check that for $\lambda\gg 1$,
$$
\|v_{\lambda}\|_{L^{4}(\Theta)}\sim \lambda^{-\frac{1}{2}},\quad 
\|v_{\lambda}\|_{L^{2}(\Theta)}\sim \lambda^{-1},\quad\|v_{\lambda}\|_{H^{1}(\Theta)}\sim 1.
$$
This completes the proof of Lemma~\ref{scaling}.
\end{proof}
\section{Proof of Theorem~\ref{thm2} (integrability and convergence properties)}
\subsection{Convergence in measure}
Let us define the $\mu$ measurable functions $f$ and $f_{N}$ by
$$
f(u)\equiv\chi\big(\|u\|_{L^2(\Theta)}\big) \exp\Big(-\int_{\Theta}V(u)\Big)
$$
and
$$ 
f_{N}(u)\equiv\chi\big(\|S_{N}(u)\|_{L^2(\Theta)}\big) \exp\Big(-\int_{\Theta}V(S_{N}(u))\Big)\,.
$$
We start by the following convergence property.
\begin{lemme}\label{lem2}
The sequence $(f_{N}(u))_{N\in\N}$ converges in measure as $N$ tends to infinity, 
with respect to the measure $\mu$, to
$f(u)$.
\end{lemme}
\begin{proof}
Since $\chi$ and the exponential are continuous functions, it suffices to show that
the sequence $\|S_{N}u\|_{L^2(\Theta)}$ converges in measure as $N$ tends to infinity, 
with respect to the measure $\mu$, to
$\|u\|_{L^2(\Theta)}$ and that
the sequence $\int_{\Theta}V(S_{N}(u))$ converges in measure as $N$ tends to infinity, 
with respect to the measure $\mu$, to $\int_{\Theta}V(u)$.
Thanks to the Chebishev inequality, it therefore suffices to prove that  
$\|S_{N}u\|_{L^2(\Theta)}$ converges in $L^2(d\mu(u))$ to $\|u\|_{L^2(\Theta)}$ and that
$\int_{\Theta}V(S_{N}(u))$ 
converges in $L^1(d\mu(u))$ to $\int_{\Theta}V(u)$.
The first assertion is trivial and the second one follows from Lemma~\ref{rudin}.
This completes the proof of Lemma~\ref{lem2}.
\end{proof}
\subsection{A gaussian estimate}
We now state a property of the measure $\mu$ resulting from its 
gaussian nature.
\begin{lemme}\label{gauss}
Let $\sigma\in [s,1/2[$. There exist $C>0$ and $c>0$ such that for every integers $M\geq N\geq 0$ 
(with the convention that $S_{0}\equiv 0$), every real number $\lambda\geq 1$,
$$
\mu
\Big(
u\in H^s_{rad}(\Theta)\, :\, 
\big\|S_{M}(u)-S_{N}(u)\big\|_{H^{\sigma}(\Theta)}>\lambda\Big)
\leq Ce^{-c\lambda^{2}(1+N)^{2(1-\sigma)}}\,.
$$
\end{lemme}
\begin{proof}
We follow the argument given in \cite[Proposition~3.3]{Tz}.
It suffices to prove that $p(A_{N,M})\leq C\exp\big(-c\lambda^{2}(1+N)^{2(1-\sigma)}\big)$, where
$$
A_{N,M}\equiv \Big(\omega\in\Omega\, :\, 
\big\|\varphi_{M}(\omega,\cdot)-\varphi_{N}(\omega,\cdot)\big\|_{H^{\sigma}(\Theta)}>\lambda\Big)\,.
$$
Let $\theta>0$ be such that $2\theta<1-2\sigma$. Notice that a proper choice of $\theta$ is possible 
thanks to the assumption $\sigma<1/2$.
For $0\leq N_1\leq N_2$ two integers and $\kappa>0$, 
we consider the set $A_{N_1,N_2,\kappa}$, defined by 
$$
A_{N_1,N_2,\kappa}\equiv \Big(\omega\in\Omega\, :\, 
\big\|\varphi_{N_2}(\omega,\cdot)-\varphi_{N_1}(\omega,\cdot)\big\|_{H^{\sigma}(\Theta)}>\kappa\lambda
\big((1+N_2)^{-\theta}+
\Big(\frac{1+N}{1+N_2}\Big)^{1-\sigma}\big)\Big)\,.
$$
Let $L_1$, $L_2$ be two dyadic integers such that 
$$L_1/2<1+N\leq L_1,\quad L_2\leq M< 2L_2.$$
We will only analyse the case $L_1\leq L_{2}/2$. If $L_{1}> L_{2}/2$ then the
analysis is simpler. Indeed, if $L_{1}> L_{2}/2$ then
$L_1\geq L_2$ which implies
$$
L_1/2<1+N\leq 1+M<1+2L_2<4L_1
$$
and the analysis of the case $L_1\leq L_{2}/2$ below (see (\ref{lidle2}),
(\ref{lidle3})) can be performed to this case by writing
$$
\varphi_{M}-\varphi_{N}=(\varphi_{L_1}-\varphi_{N})+(\varphi_{M}-\varphi_{L_1})
$$
(without the summation issue).
We thus assume that $L_1\leq L_{2}/2$.
Write 
$$\varphi_{M}-\varphi_{N}=(\varphi_{L_1}-\varphi_{N})+\Big(\sum_{\stackrel{L_1\leq L\leq  L_2/2}
{ L-{\rm dyadic }}}(\varphi_{2L}-\varphi_{L})\Big)+(\varphi_{M}-\varphi_{L_2}).$$
Using the triangle inequality and summing-up geometric series, we obtain that
there exists a sufficiently small $\kappa>0$ 
depending on $\sigma$ but independent of $\lambda$, $N$ and $M$ such that
\begin{equation}\label{union}
A_{N,M}\subset\,A_{N,L_1,\kappa}\bigcup\,\Big(\bigcup_{\stackrel{L_1\leq L\leq  L_2/2}{ L-{\rm dyadic }}} \,
A_{L,2L,\kappa}\Big)\,\bigcup A_{L_2,M,\kappa}\,.
\end{equation}
Since $z_n \sim n$, for $\omega \in A_{L,2L,\kappa}$,
$$\sum_{n=L+1}^{2L}|g_{n}(\omega)|^{2}\geq c\lambda^{2}L^{2-2\sigma}
\big(L^{-2\theta}+(L^{-1}(1+N))^{2-2\sigma}\big).
$$
Therefore using Lemma~\ref{armen} and that $2-2\sigma-2\theta>1$, we obtain that for $\lambda\geq 1$,
\begin{equation}\label{lidle1}
p(A_{L,2L,\kappa})\leq e^{c_1L-c_{2}\lambda^{2}(L^{2-2\sigma-2\theta}+(1+N)^{2-2\sigma})}
\leq
e^{-c\lambda^2 (1+N)^{2-2\sigma}}e^{-c\lambda^{2}L^{2-2\sigma-2\theta}},
\end{equation}
where the constant $c>0$ is independent of $L,N,M$ and $\lambda$. 
Similarly 
\begin{equation}\label{lidle2}
p(A_{N,L_1,\kappa})\leq e^{-c\lambda^2 (1+N)^{2-2\sigma}}e^{-c\lambda^{2}L_1^{2-2\sigma-2\theta}}
\leq e^{-c\lambda^2 (1+N)^{2-2\sigma}}
\end{equation}
and
\begin{equation}\label{lidle3}
p(A_{L_2,M,\kappa})\leq e^{-c\lambda^2 (1+N)^{2-2\sigma}}e^{-c\lambda^{2}L_2^{2-2\sigma-2\theta}}
\leq e^{-c\lambda^2 (1+N)^{2-2\sigma}}\,.
\end{equation}
Collecting estimates (\ref{lidle1}), (\ref{lidle2}), (\ref{lidle3}), coming back to (\ref{union})
and summing an obviously convergent series in $L$ completes the proof of Lemma~\ref{gauss}.
\end{proof}
\subsection{Uniform integrability}
We next prove the crucial uniform integrability property of $f_{N}$.
\begin{lemme}\label{lem5/2}
Let us fix $p\in[1,\infty[$. Then there exists $C>0$ such that for every $M\in\N$,
$$
\int_{H^s_{rad}(\Theta)}|f_{M}(u)|^{p}d\mu(u)\leq C\,.
$$
\end{lemme}
\begin{proof}
Using (\ref{defocus}), we observe that it suffices to prove that
$$
\exists\,\, C>0,\,\,\forall M\in\N,\,\,\,
\int_{\Omega}\chi^{p}\big(\|\varphi_{M}(\omega,\cdot)\|_{L^2(\Theta)}\big)
\exp\big(Cp\|\varphi_{M}(\omega,\cdot)\|^{\beta}_{L^{\beta}(\Theta)}\big)dp(\omega)\leq C.
$$
Using the Sobolev inequality, we infer that
$$
\|\varphi_{M}(\omega,\cdot)\|_{L^{\beta}(\Theta)}\leq C\|\varphi_{M}(\omega,\cdot)\|_{H^{\sigma}(\Theta)},
$$
provided
\begin{equation}\label{sob}
\sigma\geq 2\big(\frac{1}{2}-\frac{1}{\beta}\big)\,.
\end{equation}
Observe that since $\beta<4$ there exists $\sigma\in [s,1/2[$ satisfying (\ref{sob}).
Let us fix such a value of $\sigma$ for the sequel of the proof.
Since $\chi$ is with compact support, we need to study the convergence of the integral
$$
\int_{\lambda_0}^{\infty}h_{M}(\lambda)d\lambda,
$$
with
$$
h_{M}(\lambda)\equiv p\Big(\omega\in\Omega\,:\, \|\varphi_{M}(\omega,\cdot)\|_{H^{\sigma}(\Theta)}\geq
c(\log(\lambda))^{\frac{1}{\beta}},\quad 
\|\varphi_{M}(\omega,\cdot)\|_{L^{2}(\Theta)}\leq C
\Big),
$$
where $c$ and $C$ are independent of $\lambda$ and $M$ ($C$ is depending on
the support of $\chi$) and $\lambda_0$ is a large constant, independent of $M$, to be fixed later.

Since for $N\leq M$,
\begin{equation}\label{empty}
\|\varphi_{N}(\omega,\cdot)\|_{H^{\sigma}(\Theta)}\leq 
C N^{\sigma}\|\varphi_{N}(\omega,\cdot)\|_{L^{2}(\Theta)}
\leq 
C N^{\sigma}\|\varphi_{M}(\omega,\cdot)\|_{L^{2}(\Theta)}
\end{equation}
we obtain that there exists $\alpha>0$, independent of $M$ and $\lambda$ such that if $M$ satisfies
$
M \leq 
\alpha (\log(\lambda))^{\frac{1}{\sigma\beta}}
$
then $h_{M}(\lambda)=0$ (use (\ref{empty}) with $M=N$).
We can therefore assume that $M>\alpha
(\log(\lambda))^{\frac{1}{\sigma\beta}}$. 

Let us fix $\lambda\geq \lambda_0$. Define $N$ as the integer part of
$
\alpha(\log(\lambda))^{\frac{1}{\sigma\beta}-\delta},
$
where $\delta$ is such that
\begin{equation}\label{viena}
0<\delta<\frac{2-\sigma \beta}{2\sigma\beta(1-\sigma)}\,.
\end{equation}
Let us notice that a proper choice of $\delta$ is possible since $\beta<4$ and $\sigma<1/2$.
Observe also that for $\lambda_0\gg 1$, depending only on $\alpha$, we have
$N\geq 1$ and $N\leq M$.
Using (\ref{empty}), we obtain that the event
$$
\Big(\omega\in\Omega\,:\, \|\varphi_{N}(\omega,\cdot)\|_{H^{\sigma}(\Theta)}\geq
\frac{c}{2}(\log(\lambda))^{\frac{1}{\beta}},\quad 
\|\varphi_{M}(\omega,\cdot)\|_{L^{2}(\Theta)}\leq C
\Big)
$$
is of probability zero for $\lambda\geq\lambda_0$, where $\lambda_0$ is a large constant independent of $M$. 
At this place we fix the value of $\lambda_0$.
Using the triangle inequality, we obtain that for $\lambda\geq\lambda_0$,
$$
h_{M}(\lambda)\leq 
p\Big(\omega\in\Omega\,:\, \|\varphi_{M}(\omega,\cdot)-\varphi_{N}(\omega,\cdot)\|_{H^{\sigma}(\Theta)}\geq
\frac{c}{2}(\log(\lambda))^{\frac{1}{\beta}}
\Big).
$$
Using Lemma~\ref{gauss}, we arrive at
\begin{eqnarray*}
h_{M}(\lambda) & \leq & 
Ce^{-c(\log(\lambda))^{\frac{2}{\beta}}(1+N)^{2(1-\sigma)}}
\\
& \leq &
Ce^{-c(\log(\lambda))^{\frac{2}{\beta}}
(\log(\lambda))^{\frac{2(1-\sigma)}{\sigma\beta}-2\delta(1-\sigma)}}
\\
& = &
Ce^{-c(\log(\lambda))^{\frac{2}{\sigma\beta}-2\delta(1-\sigma)}}
\,.
\end{eqnarray*}
Thanks to (\ref{viena}), we have that $\frac{2}{\sigma\beta}-2\delta(1-\sigma)>1$ 
and therefore $h_{M}(\lambda)$ is 
integrable on the interval $[\lambda_0,\infty[$. 
The integrability on $[0,\lambda_0]$ is direct since $0\leq h_{M}(\lambda)\leq
1$. This completes the proof of Lemma~\ref{lem5/2}.
\end{proof}
\begin{remarque}
The exponent $\beta=4$ appears as critical in the above argument, a fact which reflects
the critical nature of the cubic non-linearity for the $2d$ NLS. This fact may be related to 
a blow-up for the cubic focusing NLS for data of positive $\mu$ measure. This is however an open problem
(see the final section of \cite{Tz}).
\end{remarque}
Using Lemma~\ref{lem5/2}, we readily arrive at the following statement.
\begin{lemme}\label{gauss-bis}
Let $\sigma\in [s,1/2[$. There exist $C>0$ and $c>0$ such that for every integer $M\geq 1$, 
every real number $\lambda\geq 1$,
$$
\tilde{\rho}_{M}
\Big(
u\in H^s_{rad}(\Theta)\, :\, \|S_{M}(u)\|_{H^{\sigma}(\Theta)}>\lambda\Big)
\leq Ce^{-c\lambda^{2}}\,.
$$
\end{lemme}
\begin{proof}
It suffices to use the Cauchy-Schwarz inequality, Lemma~\ref{gauss} and Lemma~\ref{lem5/2}.
\end{proof}
Another consequence of  Lemma~\ref{lem5/2} is the integrability of $f(u)$.
\begin{lemme}\label{integrability}
For every $p\in[1,\infty[$, $f(u)\in L^p(H^s_{rad};{\mathcal B},d\mu(u))$.
\end{lemme} 
\begin{proof}
Using Lemma~\ref{lem2}, we obtain that there is a sub-sequence $N_k$ such that the sequence 
$(f_{N_k}(u))_{k\in\N}$ converges to 
$f(u)$, $\mu$ almost surely. Thanks to Lemma~\ref{lem5/2} $(f_{N_k}(u))_{k\in\N}$ is uniformly bounded
in $L^p(H^s_{rad}(\Theta),{\mathcal B}, d\mu)$. Using Fatou's lemma we deduce that $f(u)$ belongs to 
$L^p(H^s_{rad}(\Theta),{\mathcal B}, d\mu)$ with a norm bounded by the liminf of the norms of $f_{N_k}(u)$'s.
This completes the proof of Lemma~\ref{integrability}.
\end{proof}
\subsection{End of the proof of Theorem~\ref{thm2}}
We have the following convergence property which yields the assertion of Theorem~\ref{thm2} 
in the particular case $U=F=H^s_{rad}(\Theta)$.
\begin{lemme}\label{lem3}
Let us fix $p\in[1,\infty[$.
The following holds true :
$$
\lim_{N\rightarrow\infty}\int_{H^s_{rad}(\Theta)}|f_{N}(u)-f(u)|^{p}d\mu(u)=0\,.
$$
\end{lemme}
\begin{proof}
Let us fix $\varepsilon >0$. Consider the set
$$
A_{N,\varepsilon}\equiv \big(
u\in H^{s}_{rad}(\Theta)\,:\, |f_{N}(u)-f(u)|\leq \varepsilon
\big).
$$
Denote by $A_{N,\varepsilon}^{c}$ the complementary set in
$H^{s}_{rad}(\Theta)$ of $A_{N,\varepsilon}$. 
Observe that $f$ and $f_{N}$ belong to $L^{2p}(d\mu)$ with norms bounded uniformly in $N$.
Then, using the H\"older inequality, we get
$$
\Big|
\int_{A_{N,\varepsilon}^{c}}|f_{N}(u)-f(u)|^{p}d\mu(u)
\Big|^{\frac{1}{p}}
\leq
\|f_N-f\|_{L^{2p}(d\mu)}[\mu(A_{N,\varepsilon}^{c})]^{\frac{1}{2p}}\leq
C[\mu(A_{N,\varepsilon}^{c})]^{\frac{1}{2p}}\,.
$$
On the other hand
$$
\int_{A_{N,\varepsilon}}|f_{N}(u)-f(u)|^{p}d\mu(u)\leq\varepsilon^{p}
$$
and thus we have the needed assertion since the convergence in measure of $f_N$ to $f$
implies that for a fixed $\varepsilon$,
$
\lim_{N\rightarrow \infty}\mu(A_{N,\varepsilon}^{c})=0.
$
This completes the proof of Lemma~\ref{lem3}.
\end{proof}
We can now turn to the proof of Theorem~\ref{thm2}. We follow the arguments of \cite[Lemma~3.8]{Tz}.
If we set
$$
U_{N}\equiv\big\{u\in H^s_{rad}(\Theta)\,:\, S_{N}(u)\in U\big\}
$$
then
$$
U\subset \liminf_{N}(U_{N}),
$$
where
$$
\liminf_{N}(U_{N})\equiv \bigcup_{N=1}^{\infty}\bigcap_{N_1= N}^{\infty}U_{N_1}\,.
$$
Indeed, we have that for every $u\in H^{\sigma}_{rad}(\Theta)$,
$S_{N}(u)$ converges to $u$ in $H^{\sigma}_{rad}(\Theta)$, as $N$ tends to $\infty$.
Therefore, using that $U$ is an open set, we conclude that for every $u\in U$
there exists $N_{0}\geq 1$ such that for $N\geq N_0$ one has $u\in
U_{N}$. Hence we have 
$
U\subset \liminf_{N}(U_{N}).
$
If $A$ is a $\rho$-measurable set, we denote by $\11_{A}$ the characteristic
function of $A$. 
Notice that thanks to the property $U\subset \liminf_{N}(U_{N})$,
$$
\liminf_{N\rightarrow\infty}\11_{U_{N}}\geq\11_{U}.
$$
Recall that
$$
\tilde{\rho}_{N}(U)=\rho_{N}(U\cap E_{N})=\int_{H^s_{rad}(\Theta)}\11_{U_N}(u)f_{N}(u)d\mu(u)\,.
$$
Using Lemma~\ref{lem3}, we observe that
$$
\lim_{N\rightarrow\infty}\Big(\int_{H^s_{rad}(\Theta)}\11_{U_N}(u)f_{N}(u)d\mu(u)-
\int_{H^s_{rad}(\Theta)}\11_{U_N}(u)f(u)d\mu(u)\Big)=0\,.
$$
Next, using the Fatou lemma, we get
\begin{eqnarray*}
\liminf_{N\rightarrow\infty}\rho_{N}(U\cap E_{N}) & = &
\liminf_{N\rightarrow\infty}
\int_{H^s_{rad}(\Theta)}\11_{U_N}(u)f(u)d\mu(u)
\\
& \geq &
\int_{H^s_{rad}(\Theta)}\11_{U}(u)f(u)d\mu(u)
\\
& = &\int_{U}f(u)d\mu(u)=\rho(U)\,.
\end{eqnarray*}
This proves (\ref{parvo}). 
Observe that  Lemma~\ref{lem3} implies that
$$
\lim_{N\rightarrow\infty}\rho_{N}(E_{N})=\rho(H^s_{rad}(\Theta))\,.
$$
Therefore to prove (\ref{vtoro}), it suffices to use (\ref{parvo}) by passing to complementary sets
(as in \cite{Tz}, we could give a direct proof of (\ref{vtoro})).
This completes the proof of Theorem~\ref{thm2}.
\qed
\begin{remarque}
Let us observe that the reasoning in the proof of Theorem~\ref{thm2} is of quite general nature.
It suffices to know that : 
\begin{itemize}
\item
$(f_{N})$ is bounded uniformly with respect to $N$ in $L^{p}(d\mu)$ for some $p>1$.
\item
$(f_{N})$ converges to $f$ in measure.
\end{itemize}
\end{remarque}
\subsection{A corollary of Theorem~\ref{thm2}}
Combining Lemma~\ref{gauss-bis} and Theorem~\ref{thm2}, we arrive at the following statement.
\begin{lemme}\label{gauss-tris}
Let $\sigma\in [s,1/2[$. There exist $C>0$ and $c>0$ such that for every real number $\lambda\geq 1$,
$$
\rho\Big(u\in H^s_{rad}(\Theta)\, :\, \|u\|_{H^{\sigma}(\Theta)}\in ]\lambda,\infty[\Big)\leq Ce^{-c\lambda^{2}}\,.
$$
\end{lemme}
\begin{proof}
It suffices to apply Theorem~\ref{thm2} to the open set of $H^{\sigma}_{rad}(\Theta)$,
$$
U=\Big(u\in H^s_{rad}(\Theta)\, :\, \|u\|_{H^{\sigma}(\Theta)}\in]\lambda,\infty[\Big)
$$
and to observe that $\tilde{\rho}_{N}(U)=\tilde{\rho}_{N}(U_{N})$, where 
\begin{equation*}
U_{N}=\Big(u\in H^s_{rad}(\Theta)\, :\, \|S_{N}(u)\|_{H^{\sigma}(\Theta)}\in]\lambda,\infty[\Big).
\end{equation*}
Thus by Lemma~\ref{gauss-bis}, $\tilde{\rho}_{N}(U)\leq C\exp(-c\lambda^2)$ which,
combined with Theorem~\ref{thm2}, completes the proof of
Lemma~\ref{gauss-tris}.
\end{proof}
\begin{remarque}\label{rem}
As a consequence of Lemma~\ref{gauss-tris} one obtains that for $\sigma\in [s,1/2[$ one has 
$\rho(H^{\sigma}_{rad}(\Theta))=\rho(H^{s}_{rad}(\Theta))$.
Moreover for every $\rho$ measurable set $A$,
$$
\rho\Big(u\in A\, :\, \|u\|_{H^{\sigma}(\Theta)}\in]\lambda,\infty[\Big)\leq Ce^{-c\lambda^{2}}\,.
$$
and thus $A$ may be approximated by bounded sets of $H^{\sigma}_{rad}(\Theta)$
(the intersections of $A$ and the balls of radius $\lambda\gg 1$ centered at the origin of 
$H^{\sigma}_{rad}(\Theta)$).
\end{remarque}
\section{Bourgain spaces and bilinear estimates}
The following two statements play a crucial role in the analysis of the local Cauchy problem for (\ref{1}).
\begin{proposition}\label{str1}
For every $\varepsilon>0$, there exists $\beta<1/2$, there exists $C>0$ such that for
every $N_1,N_2\geq 1$, every $L_1,L_2\geq 1$, every $u_1$, $u_2$ two functions on
$\R\times\Theta$ of the form
$$
u_{j}(t,r)=\sum_{N_j\leq \langle z_n\rangle < 2N_j}\,c_j(n,t)\, e_{n}(r),\quad j=1,2
$$
where the Fourier transform of $c_j(n,t)$ with respect to $t$ satisfies 
$$
{\rm supp}\, \widehat{c_j}(n,\tau)\subset \{
\tau\in\R\,:\, L_{j}\leq \langle\tau+z_n^2\rangle\leq 2L_j\},\quad j=1,2
$$
one has the bound
$$
\|u_1 u_2\|_{L^2(\R\times \Theta)}\leq C(\min(N_1, N_2))^{\varepsilon}(L_1L_2)^{\beta}
\|u_1\|_{L^2(\R\times \Theta)}\|u_2\|_{L^2(\R\times \Theta)}\,.
$$
\end{proposition}
\begin{proposition}\label{str2}
For every $\varepsilon>0$, there exists $\beta<1/2$, there exists $C>0$ such that for
every $N_1,N_2\geq 1$, every $L_1,L_2\geq 1$, every $u_1$, $u_2$ two functions on
$\R\times\Theta$ of the form
$$
u_{1}(t,r)=\sum_{N_1\leq \langle z_n\rangle < 2N_1}\,c_1(n,t)\, e_{n}(r)
$$
and
$$
u_{2}(t,r)=\sum_{N_2\leq \langle z_n\rangle < 2N_2}\,c_2(n,t)\, e'_{n}(r)
$$
where the Fourier transform of $c_j(n,t)$ with respect to $t$ satisfies 
$$
{\rm supp}\, \widehat{c_j}(n,\tau)\subset \{\tau\in\R\,:\, 
L_{j}\leq \langle\tau+z_n^2\rangle\leq 2L_j\},\quad j=1,2
$$
one has the bound
$$
\|u_1 u_2\|_{L^2(\R\times \Theta)}\leq C(\min(N_1, N_2))^{\varepsilon}(L_1L_2)^{\beta}
\|u_1\|_{L^2(\R\times \Theta)}\|u_2\|_{L^2(\R\times \Theta)}\,.
$$
\end{proposition}
For the proof of Propositions~\ref{str1} and \ref{str2} we refer to \cite[Proposition~4.1]{Tz}
and \cite[Proposition~4.3]{Tz} respectively.
The results of Propositions~\ref{str1} and \ref{str2} can be injected in the framework of the Bourgain spaces 
associated to the Schr\"odinger equation on the disc, in order to get local existence results for (\ref{1}).
Following \cite{Tz}, we define the Bourgain spaces $X^{\sigma,b}_{rad}(\R\times\Theta)$ of
functions on $\R\times\Theta$ which are radial with respect to the second
argument, equipped with the norm
$$
\|u\|_{X^{\sigma,b}_{rad}(\R\times\Theta)}^{2}=\sum_{n= 1}^{\infty} z_{n}^{2\sigma}\|
\langle\tau+z_n^2\rangle^{b}\widehat{\langle u(t),e_{n}\rangle}(\tau)\|_{L^{2}(\R_{\tau})}^{2}\,,
$$
where $\langle\cdot,\cdot\rangle$ stays for the $L^2(\Theta)$ pairing and $\widehat{\cdot}$ denotes 
the Fourier transform on $\R$.
We can express the norm in
$X^{\sigma,b}_{rad}(\R\times\Theta)$ in terms of the localisation
operators $\Delta_{N,L}$. More precisely, if for 
$N,L$ positive integers, we define $\Delta_{N,L}$ by
\begin{equation*}
\Delta_{N,L}(u)=\frac{1}{2\pi}\sum_{n\,:\, N\leq\langle z_n\rangle< 2N}
\Big(\int_{L\leq \langle\tau+z_n^2\rangle\leq 2L}
\widehat{\langle u(t),e_{n}\rangle}(\tau)e^{it\tau}d\tau\Big)e_{n},
\end{equation*}
then we can write 
\begin{equation*}
\|u\|_{X^{\sigma,b}_{rad}(\R\times\Theta)}^{2} \approx_{\sigma,b}
\sum_{L,N-{\rm dyadic }}L^{2b}N^{2\sigma} \|\Delta_{N,L}(u)\|_{L^2(\R\times\Theta)}^{2}\,,
\end{equation*}
where $\approx_{\sigma,b}$ means that the implicit constant may depend on $\sigma$ and $b$.
Using that (see \cite{Tz}),
$$
\exists\, C>0\,\, : \,\,\forall n\in \N,\,\, \|e_{n}\|_{L^{\infty}(\Theta)}\leq Cn^{\frac{1}{2}}
$$
and the Cauchy-Schwarz inequality in the $\tau$ integration and in the $n$ summation, we arrive at the bound
\begin{equation}\label{infty}
\|\Delta_{N,L}(u)\|_{L^{\infty}(\R\times \Theta)}\leq 
L^{\frac{1}{2}}N\|\Delta_{N,L}(u)\|_{L^{2}(\R\times \Theta)}\,.
\end{equation}
Let us next analyse $\partial_{r}(\Delta_{N,L}(u))$. We can write 
\begin{equation*}
\Delta_{N,L}(u)=\sum_{N\leq \langle z_{n}\rangle < 2N}\,c(n,t)\, e_{n}(r),\quad
{\rm supp}\, \widehat{c}(n,\tau)\subset \{\tau\in\R\,:\, L\leq \langle\tau+z_{n}^2\rangle\leq 2L\}
\end{equation*}
and thus
\begin{equation}\label{partial_r}
\partial_{r}\big(\Delta_{N,L}(u)\big)=\sum_{N\leq \langle z_{n}\rangle < 2N}\,c(n,t)\, e'_{n}(r).
\end{equation}
Recall (see \cite{Tz}) that for $m\neq n$, $e'_{m}$ and $e'_{n}$ are orthogonal in $L^2(\Theta)$ and
$
\|e'_{n}\|_{L^2(\Theta)}\approx n
$.
Therefore 
\begin{equation*}
\big\|\partial_{r}\big(\Delta_{N,L}(u)\big)\big\|_{L^2(\R\times\Theta)}^{2}=
c\sum_{N\leq \langle z_{n}\rangle < 2N}\,\|e'_{n}\|_{L^2(\Theta)}^{2}
\int_{-\infty}^{\infty}|\widehat{c}(n,\tau)|^{2}d\tau
\end{equation*}
and thus
\begin{equation}\label{jm}
\big\|\partial_{r}\big(\Delta_{N,L}(u)\big)\big\|_{L^2(\R\times\Theta)}\approx 
N\big\|\Delta_{N,L}(u)\big\|_{L^2(\R\times\Theta)}\,.
\end{equation}
Using \cite[Lemma~2.1]{Tz}, $\|\partial_{r}e_{n}\|_{L^{\infty}(\Theta)}\leq Cn^{3/2}$ and thus coming back to
(\ref{partial_r}), after writing $c(n,t)$ in terms of its Fourier transform and 
applying the Cauchy-Schwarz inequality in the $\tau$ (the dual of $t$ variable) integration, we obtain that
\begin{equation}\label{infty-bis}
\big\|\partial_{r}\big(\Delta_{N,L}(u)\big)\big\|_{L^{\infty}(\R\times\Theta)}\leq C L^{\frac{1}{2}}
N^{2}\big\|\Delta_{N,L}(u)\big\|_{L^2(\R\times\Theta)}\,.
\end{equation}
Let us next define two other projectors involved in the well-posedness analysis of (\ref{1}).
The projector $\Delta_{N}$ is defined by
$$
\Delta_{N}(u)=\sum_{n\,:\, N\leq\langle z_n\rangle< 2N}\, \langle u, e_n\rangle e_n\,.
$$
For $N\geq 2$ a dyadic integer, we define the projector $\tilde{S}_N$ by
$$
\tilde{S}_{N}= 
\sum_{\stackrel{N_1\leq N/2}{ N_1-{\rm dyadic }}}\Delta_{N_1}\, .
$$
For a notational convenience, we assume that $\tilde{S}_{1}$ is zero.
\section{Nonlinear estimates in Bourgain spaces}
In this section, we shall derive nonlinear estimates related to the problems (\ref{1}) and (\ref{N}).
We start by a lemma which improves on the Sobolev embedding.
\begin{lemme}\label{Lp}
Let us fix $p\geq 4$. Then for every $b>\frac{1}{2}$, $\sigma>1-\frac{4}{p}$ there exists $C>0$ such that 
for all $u\in X^{\sigma,b}_{rad}(\R\times\Theta)$ one has
\begin{equation}\label{che}
\|u\|_{L^p(\R\times\Theta)}\leq C\|u\|_{X^{\sigma,b}_{rad}(\R\times\Theta)}\,.
\end{equation}
\end{lemme}
\begin{proof}
It suffices to prove the assertion for $p=4$ and $p=\infty$.
Let us first consider the case $p=4$. 
Observe that $\Delta_{N,L}(u)$ fits in the scope of applicability of Proposition~\ref{str1}.
Using Proposition~\ref{str1} with $\varepsilon=\sigma/2>0$, we obtain that
$$
\|\Delta_{N,L}(u)\|_{L^4(\R\times\Theta)}\leq 
C\|\Delta_{N,L}(u)\|_{X^{\sigma/2,\beta}_{rad}(\R\times\Theta)}\,.
$$
Therefore, by writing $u=\sum_{L,N}\Delta_{N,L}(u)$, where the summation runs over the dyadic values of $L$, 
$N$, by summing geometric series in $N$ and $L$ , 
we obtain that (\ref{che}) holds true for $p=4$ (observe that we use
Proposition~\ref{str1} with $\varepsilon=\sigma/2$
instead of $\sigma$ in order to get small negative powers of $N$ and $L$ after
applying the triangle inequality to $\sum_{L,N}\Delta_{N,L}(u)$.
Let us next consider the case $p=\infty$.
In this case, the assertion holds true thanks to (\ref{infty}) and another summation of geometric series.
This completes the proof of Lemma~\ref{Lp}.
\end{proof}
The next lemma gives sense of $F(u)$, in the scale of Bourgain's spaces, for $u$ of low regularity.
\begin{lemme}\label{F(u)}
Let $(b,\sigma)$ be such that $\max(1/3,1-2/\alpha)<\sigma<1/2$, $b>1/2$ and let 
$u\in X^{\sigma,b}_{rad}(\R\times\Theta)$.
Then $F(u)\in X^{-\sigma,-b}_{rad}(\R\times\Theta)$. Moreover
$$
\lim_{N\rightarrow \infty}\|F(u)-F(\tilde{S}_{N}(u))\|_{X^{-\sigma,-b}_{rad}(\R\times\Theta)}=0\,.
$$
\end{lemme}
\begin{proof}
For $v\in X^{\sigma,b}_{rad}(\R\times\Theta)$, we write
\begin{eqnarray*}
\int_{\R\times\Theta}|F(u)v| 
& \leq & C\big(\int_{\R\times \Theta}|uv|+\int_{\R\times \Theta}|u|^{\alpha+1}|v|\big)
\\
& \leq & C\|u\|_{X^{\sigma,b}_{rad}(\R\times\Theta)}\|v\|_{X^{\sigma,b}_{rad}(\R\times\Theta)}
+C\|u\|^{\alpha+1}_{L^{\alpha+2}(\R\times\Theta)}\|v\|_{L^{\alpha+2}(\R\times\Theta)}
\end{eqnarray*}
Now, using Lemma~\ref{Lp}, we get
$$
\|u\|_{L^{\alpha+2}(\R\times\Theta)}\leq C\|u\|_{X^{\sigma_1,b}_{rad}(\R\times\Theta)}\,,
$$
where $\sigma_1>0$, when $\alpha\leq 2$ and $\sigma_1>1-4/(\alpha+2)$ when $\alpha\in]2,4[$.
Observing that for $\alpha\geq 2$, $\max(1/3,1-2/\alpha)\geq 1-4/(\alpha+2)$ shows that 
$$
\int_{\R\times\Theta}|F(u)v|\leq C\|u\|_{X^{\sigma,b}_{rad}(\R\times\Theta)}\|v\|_{X^{\sigma,b}_{rad}(\R\times\Theta)}
\big(1+\|u\|_{X^{\sigma,b}_{rad}(\R\times\Theta)}^{\alpha}\big)
$$
and thus $F(u)\in X^{-\sigma,-b}_{rad}(\R\times\Theta)$.
Similarly one shows that
$$
\int_{\R\times\Theta}|(F(u)-F(\tilde{S}_{N}(u)))v|\leq C\|u-\tilde{S}_{N}(u)\|_{X^{\sigma,b}_{rad}}
\|v\|_{X^{\sigma,b}_{rad}}\big(1+\|u\|_{X^{\sigma,b}_{rad}}^{\alpha}\big)
$$
which yields the needed convergence. This completes the proof of Lemma~\ref{F(u)}.
\end{proof}
One may prove a statement similar to Lemma~\ref{F(u)} with $\tilde{S}_{N}$ replaced by $\tilde{S}_{N,L}$
where $\tilde{S}_{N,L}$ is defined similarly to $\tilde{S}_{N}$ with $\Delta_{N_1}$ replaced by 
$\Delta_{N_1,L_1}$, $L_1\leq L$. This observation allows to consider only finite sums over dyadic integers 
in the proof of the next proposition (one can also apply a similar approximation argument to $v$ involved in
(\ref{vanbis})).

In fact a much stronger statement then Lemma~\ref{F(u)} holds true.
It turns out that under the assumptions of Lemma~\ref{F(u)} one has 
$F(u)\in X^{\sigma,-b}_{rad}(\R\times\Theta)$.
\begin{proposition}\label{main}
Let 
$\max(1/3,1-2/\alpha)<\sigma_1\leq\sigma<1/2$.
Then there exist two positive numbers $b,b'$ such that $b+b'<1$, $b'<1/2<b$, there
exists $C>0$ such that for every $u,v\in X^{\sigma,b}_{rad}(\R\times\Theta)$,
\begin{equation}\label{main1}
\|F(u)\|_{X^{\sigma,-b'}_{rad}(\R\times\Theta)}\leq
C\Big(1+\|u\|^{\max(\alpha,2)}_{X^{\sigma_{1},b}_{rad}(\R\times\Theta)}\Big)
\|u\|_{X^{\sigma,b}_{rad}(\R\times\Theta)}
\end{equation}
and
\begin{equation}\label{main2} 
\|F(u)-F(v)\|_{X^{\sigma,-b'}_{rad}}\leq C\Big(1+\|u\|^{\max(\alpha,2)}_{X^{\sigma,b}_{rad}}+
\|v\|^{\max(\alpha,2)}_{X^{\sigma,b}_{rad}}\Big)\|u-v\|_{X^{\sigma,b}_{rad}}\,.
\end{equation}
\end{proposition}
\begin{proof}
The proof of this proposition for $\alpha<2$ is given in \cite{Tz}.
We therefore may assume that $\alpha\geq 2$ in the sequel of the proof.
The proof will follow closely \cite{Tz} by incorporating an argument already appeared in \cite{BGT}.
Let us observe that in order to prove (\ref{main1}), it suffices to prove that
\begin{equation}\label{main1-bis}
\|F(\tilde{S}_{M}(u))\|_{X^{\sigma,-b'}_{rad}(\R\times\Theta)}\leq
C\Big(1+\|u\|^{\alpha}_{X^{\sigma_{1},b}_{rad}(\R\times\Theta)}\Big)
\|u\|_{X^{\sigma,b}_{rad}(\R\times\Theta)},
\end{equation}
uniformly in $M\in\N$. Indeed, if we can prove (\ref{main1-bis}) then $(F(\tilde{S}_{M}(u)))_{M\in\N}$ 
is a bounded sequence in $X^{\sigma,-b'}_{rad}(\R\times\Theta)$ 
(and thus also in $X^{-\sigma,-b}_{rad}(\R\times\Theta)$) 
and therefore it converges (up to a sub-sequence) 
weakly to some limit which satisfies the needed bound. 
In order to identify this limit with $F(u)$ it suffices make appeal to Lemma~\ref{F(u)}.
Thanks to the gauge invariance of the nonlinearity $F(u)$, we observe that
$F(u)-(\partial F)(0)u$ is vanishing at order $2$ at $u=0$ and thus in the proof of (\ref{main1-bis}), 
we may assume that 
\begin{equation}\label{vanishing}
\partial^{k_1}\bar{\partial}^{k_2}(F)(0)=0,\quad \forall\,\, k_1+k_2\leq 2.
\end{equation}
Observe that (\ref{main1-bis}) follows from the estimate
\begin{equation}\label{vanbis}
\Big|\int_{\R\times\Theta}F(\tilde{S}_{M}(u))\bar{v}\Big|
\leq C\|v\|_{X^{-\sigma,b'}_{rad}(\R\times\Theta)}
\Big(1+\|u\|^{\alpha}_{X^{\sigma_{1},b}_{rad}(\R\times\Theta)}\Big)
\|u\|_{X^{\sigma,b}_{rad}(\R\times\Theta)}
\end{equation}
(let us remark that if $v$ contains only very high frequencies with respect to
the $\Delta_{N,L}$ decomposition then the right hand-side of (\ref{vanbis}) is small).
Using that $\Delta_{N}=\tilde{S}_{2N}-\tilde{S}_{N}$ and (\ref{vanishing}), we may write
\begin{eqnarray*}
F(\tilde{S}_{M}(u)) & = &\sum_{\stackrel{2\leq N_1\leq M}{N_1-{\rm dyadic }}}
(F(\tilde{S}_{N_1}(u))-F(\tilde{S}_{N_1/2}(u))
\\
& = &
\sum_{\stackrel{N_1\leq M/2}{N_1-{\rm dyadic }}}
\big(\Delta_{N_1}(u)G_{1}(\Delta_{N_1}(u),\tilde{S}_{N_1}(u))+
\overline{\Delta_{N_1}(u)}G_{2}(\Delta_{N_1}(u),\tilde{S}_{N_1}(u))\big)
\\
& \equiv &  F_{1}(u)+F_{2}(u),
\end{eqnarray*}
where $G_{1}(z_1,z_2)$ and $G_{2}(z_1,z_2)$ are smooth functions with a
control on their growth at infinity coming from (\ref{rast})
(similar bounds to $F$ with $\alpha$ replaced $\alpha-1$).
Moreover, thanks to (\ref{vanishing}), $G_{1}(0,0)=\partial(F)(0)=0$ and $G_{2}(0,0)=\bar{\partial}(F)(0)=0$.
We will only estimate the contribution of $F_{1}(u)$ to the right hand-side of (\ref{vanbis}),
the argument for the contribution of $F_{2}(u)$ being completely analogous.
Next, we set
$$
I=\Big|\int_{\R\times\Theta}F_{1}(u)\bar{v}\Big|,\quad
I(N_0,N_1)=\Big|\int_{\R\times\Theta}\Delta_{N_1}(u)\overline{\Delta_{N_0}(v)}
G_{1}(\Delta_{N_1}(u),\tilde{S}_{N_1}(u))\Big|.
$$
Then $I\leq I_1+I_2$, where
$$
I_{1}=
\sum_{\stackrel{N_0\leq N_1\leq M/2}{N_0,N_1-{\rm dyadic }}}I(N_0,N_1),\quad
I_{2}=
\sum_{\stackrel{N_1\leq \min(N_0, M/2)}{N_0,N_1-{\rm dyadic }}}I(N_0,N_1).
$$
We first analyse $I_1$. Using (\ref{vanishing}) with $(k_1,k_2)=(1,0)$,
we decompose $G_{1}(\Delta_{N_1}(u),\tilde{S}_{N_1}(u))$ as
$$
\sum_{\stackrel{N_2\leq N_1}{N_2-{\rm dyadic }}}
\Big(G_{1}(\tilde{S}_{2N_2}\Delta_{N_1}(u),\tilde{S}_{2N_2}\tilde{S}_{N_1}(u))-
G_{1}(\tilde{S}_{N_2}\Delta_{N_1}(u),\tilde{S}_{N_2}\tilde{S}_{N_1}(u))\Big).
$$
Using that $\Delta_{N_1}\Delta_{N_2}=\Delta_{N_1}$, if $N_1=N_2$ and zero elsewhere, we obtain 
\begin{multline}\label{G1}
G_{1}(\Delta_{N_1}(u),\tilde{S}_{N_1}(u))=
\sum_{\stackrel{N_2\leq N_1}{N_2-{\rm dyadic }}}
\Delta_{N_2}(u)G_{11}^{N_2}(\Delta_{N_2}(u),\tilde{S}_{N_2}(u))
+
\\
\sum_{\stackrel{N_2\leq N_1}{N_2-{\rm dyadic }}}
\overline{\Delta_{N_2}(u)}G_{12}^{N_2}(\Delta_{N_2}(u),\tilde{S}_{N_2}(u)),
\end{multline}
where $G_{11}^{N_2}(z_1,z_2)$ and $G_{12}^{N_2}(z_1,z_2)$ are smooth functions with a
control on their growth at infinity coming from (\ref{rast}).
Moreover thanks to (\ref{vanishing}), applied with $(k_1,k_2)=(2,0)$ and $(k_1,k_2)=(1,1)$, 
we get $G_{11}^{N_2}(0,0)=0$ and $G_{12}^{N_2}(0,0)=0$. Therefore, we can expand for $j=1,2$,
\begin{multline}\label{G1j}
G_{1j}^{N_2}(\Delta_{N_2}(u),\tilde{S}_{N_2}(u))=
\sum_{\stackrel{N_3\leq N_2}{N_3-{\rm dyadic }}}
\Delta_{N_3}(u)G_{1j1}^{N_3}(\Delta_{N_3}(u),\tilde{S}_{N_3}(u))
+
\\
\sum_{\stackrel{N_3\leq N_2}{N_3-{\rm dyadic }}}
\overline{\Delta_{N_3}(u)}G_{1j2}^{N_3}(\Delta_{N_3}(u),\tilde{S}_{N_3}(u)),
\end{multline}
where, thanks to the growth assumption on the nonlinearity $F(u)$, we have
that the functions $G_{1 j_1 j_2}^{N_3}(z_1,z_2)$, $j_1,j_2\in\{1,2\}$ satisfy
\begin{equation}\label{novo}
|G_{1 j_1 j_2}^{N_3}(z_1,z_2)|\leq C(1+|z_1|+|z_2|)^{\alpha-2}.
\end{equation}
We therefore have the bound
\begin{multline*}
I_1\leq C
\sum_{\stackrel{N_0\leq N_1}{N_0,N_1-{\rm dyadic }}}
\sum_{\stackrel{N_1\geq N_2\geq N_3}{N_2,N_3-{\rm dyadic }}}
\\
\int_{\R\times\Theta}
|\Delta_{N_0}(v)\Delta_{N_1}(u)\Delta_{N_2}(u)\Delta_{N_3}(u)|
(
1+|\Delta_{N_3}(u)|+|\tilde{S}_{N_3}(u)|)^{\alpha-2}.
\end{multline*}
By splitting
$$
\Delta_{N}=\sum_{L-{\rm dyadic}}\Delta_{N,L},
$$
we may write for $b>1/2$, $0<\sigma_1<1/2$,
by using (\ref{infty}) and the Cauchy-Schwarz inequality in the $L$ summation
\begin{eqnarray*}
\|\Delta_{N_3}(u)\|_{L^{\infty}(\R\times\Theta)}
& \leq &
\sum_{L-{\rm dyadic}}\|\Delta_{N_3,L}(u)\|_{L^{\infty}(\R\times\Theta)}
\\
& \leq &
C\sum_{L-{\rm dyadic}}N_3L^{\frac{1}{2}}\|\Delta_{N_3,L}(u)\|_{L^{2}(\R\times\Theta)}
\\
& \leq & 
CN_{3}^{1-\sigma_1}\|u\|_{X^{\sigma_1,b}_{rad}(\R\times\Theta)}\,,
\end{eqnarray*}
where $C$ depends on $b$ and $\sigma_1$. Similarly
\begin{eqnarray*}
\|\tilde{S}_{N_3}(u)\|_{L^{\infty}(\R\times\Theta)}
& \leq &  
\sum_{\stackrel{N\leq N_3/2}{N-{\rm dyadic }}}
\|\Delta_{N}(u)\|_{L^{\infty}(\R\times\Theta)}
\\
& \leq &
\sum_{\stackrel{N\leq N_3/2}{L,N-{\rm dyadic }}}\|\Delta_{N,L}(u)\|_{L^{\infty}(\R\times\Theta)}
\\
& \leq &
\sum_{\stackrel{N\leq N_3/2}{L,N-{\rm dyadic }}}
CNL^{\frac{1}{2}}\|\Delta_{N,L}(u)\|_{L^{2}(\R\times\Theta)}
\\
& \leq &
C\big(\sum_{\stackrel{N\leq N_3/2}{L,N-{\rm dyadic }}}
L^{1-2b}
N^{2(1-\sigma_1)}
\big)^{\frac{1}{2}}
\big(
\sum_{\stackrel{N\leq N_3/2}{L,N-{\rm dyadic }}}
L^{2b}N^{2\sigma_1}\|\Delta_{N,L}(u)\|_{L^{2}}^{2}
\big)^{\frac{1}{2}}
\\
& \leq &
CN_{3}^{1-\sigma_1}\|u\|_{X^{\sigma_1,b}_{rad}(\R\times\Theta)}\,.
\end{eqnarray*}
Therefore
\begin{multline*}
I_1\leq C
(1+\|u\|_{X^{\sigma_1,b}(\R\times\Theta)}^{\alpha-2})
\sum_{\stackrel{N_0\leq N_1}{N_0,N_1-{\rm dyadic }}}
\sum_{\stackrel{N_1\geq N_2\geq N_3}{N_2,N_3-{\rm dyadic }}}
\\
N_{3}^{(1-\sigma_1)(\alpha-2)}
\Big(\int_{\R\times\Theta}
|\Delta_{N_0}(v)\Delta_{N_1}(u)\Delta_{N_2}(u)\Delta_{N_3}(u)|\Big).
\end{multline*}
Using Proposition~\ref{str1} and the Cauchy-Schwarz inequality, we obtain
that for every $\varepsilon>0$ there exist $\beta<1/2$ and $C_{\varepsilon}$
such that
\begin{multline*}
\int_{\R\times\Theta}
|\Delta_{N_0,L_0}(v)\Delta_{N_1,L_1}(u)\Delta_{N_2,L_2}(u)\Delta_{N_3,L_3}(u)|\leq
\\
\leq
\|\Delta_{N_0,L_0}(v)\Delta_{N_2,L_2}(u)\|_{L^2(\R\times\Theta)}
\|\Delta_{N_1,L_1}(u)\Delta_{N_3,L_3}(u)\|_{L^2(\R\times\Theta)}
\leq
\\
\leq C_{\varepsilon}(N_2 N_3)^{\varepsilon}(L_0L_1L_2L_3)^{\beta}
\|\Delta_{N_0,L_0}(v)\|_{L^2(\R\times\Theta)}
\prod_{j=1}^{3}\|\Delta_{N_j,L_j}(u)\|_{L^2(\R\times\Theta)}.
\end{multline*}
Therefore, if we set
\begin{multline}\label{Q}
Q\equiv Q(N_0,N_1,N_2,N_3,L_0,L_1,L_2,L_3)= 
CN_{0}^{-\sigma}N_1^{\sigma}(N_2 N_3)^{\sigma_{1}}L_{0}^{b'}(L_1 L_2 L_3)^{b}
\\
\times
(1+\|u\|_{X^{\sigma_1,b}(\R\times\Theta)}^{\alpha-2})
\|\Delta_{N_0,L_0}(v)\|_{L^2(\R\times\Theta)}
\prod_{j=1}^{3}\|\Delta_{N_j,L_j}(u)\|_{L^2(\R\times\Theta)},
\end{multline}
we can write
$$
I_1\leq \sum_{L_0,L_1,L_2,L_3-{\rm dyadic }}\,
\sum_{\stackrel{N_1\geq N_2\geq N_3, N_1\geq N_0}{N_0,N_1,N_2,N_3-{\rm dyadic }}}
L_{0}^{\beta-b'}(L_1L_2L_3)^{\beta-b}\Big(\frac{N_0}{N_1}\Big)^{\sigma}\,\,
\frac{N_{3}^{(1-\sigma_1)(\alpha-2)}}{(N_2N_3)^{\sigma_{1}-\varepsilon}}Q.
$$
Let us take $\varepsilon>0$ such that 
$$
\varepsilon<1-\frac{\alpha(1-\sigma_1)}{2}\,.
$$
A proper choice of $\varepsilon$ is possible thanks to the assumption $\sigma_1>1-2/\alpha$.
With this choice of $\varepsilon$ we have that $2(\sigma_1-\varepsilon)>(1-\sigma_1)(\alpha-2)$.
The choice of $\varepsilon$ fixes $\beta$ via the application of Proposition~\ref{str1}.
Then we choose $b'$ such that $\beta<b'<1/2$. We finally choose $b>1/2$ such that
$b+b'<1$. With this choice of the parameters, coming back to the definition of the projectors $\Delta_{N,L}$
and after summing geometric series in $L_0$, $L_1$, $L_2$, $L_3$, $N_2$, $N_3$, we can write that
$$
I_{1}\leq
C
(1+\|u\|_{X^{\sigma_1,b}(\R\times\Theta)}^{\alpha-2})
\|u\|^{2}_{X^{\sigma_{1},b}_{rad}(\R\times\Theta)}
\sum_{\stackrel{ N_0\leq N_1}{N_0,N_1-{\rm dyadic }}}
\Big(\frac{N_0}{N_1}\Big)^{\sigma}c(N_0)d(N_1),
$$
where
\begin{equation}\label{cd}
c(N_0)=N_{0}^{-\sigma}\|\Delta_{N_0}(v)\|_{X^{0,b'}_{rad}(\R\times\Theta)},\quad
d(N_1)=N_{1}^{\sigma}\|\Delta_{N_1}(u)\|_{X^{0,b}_{rad}(\R\times\Theta)}\,.
\end{equation}
Finally, using \cite[Lemma~6.2]{Tz}, we arrive at the bound
$$
I_1\leq
C\|v\|_{X^{-\sigma,b'}_{rad}(\R\times\Theta)}
(1+\|u\|_{X^{\sigma_1,b}(\R\times\Theta)}^{\alpha-2})
\|u\|^{2}_{X^{\sigma_{1},b}_{rad}(\R\times\Theta)}
\|u\|_{X^{\sigma,b}_{rad}(\R\times\Theta)}
$$
which ends the analysis for $I_1$.
\\

Let us now turn to the analysis of $I_2$.
The main observation is that after in integration by parts the roles of $N_0$ and $N_1$ are exchanged.
We have that
$$
I_{2}\leq \sum_{\stackrel{N_1\leq \min(N_0, M/2)}{L_0,N_0,N_1-{\rm dyadic }}}I(L_0,N_0,N_1),
$$
where
$$
I(L_0,N_0,N_1)=\Big|\int_{\R\times\Theta}\Delta_{N_1}(u)\overline{\Delta_{N_0,L_0}(v)}
G_{1}(\Delta_{N_1}(u),\tilde{S}_{N_1}(u))\Big|.
$$
Write
\begin{equation}\label{write1}
\Delta_{N_0,L_0}(v)=\sum_{N_0\leq \langle z_{n_0}\rangle < 2N_0}\,c(n_0,t)\, e_{n_0}(r),
\end{equation}
where
$$
{\rm supp}\, \widehat{c}(n_0,\tau)\subset \{\tau\in\R\,:\, L_{0}\leq \langle\tau+z_{n_0}^2\rangle\leq 2L_0\}.
$$
Define $\widetilde{\Delta}_{N_0,L_0}$ as
\begin{equation*}
\widetilde{\Delta}_{N_0,L_0}(v)\equiv
\sum_{N_0\leq \langle z_{n_0}\rangle \leq 2N_0}\,
\frac{c(n_0,t)}{z_{n_0}^2}\, e'_{n_0}(r).
\end{equation*}
Since $\| e'_{n_0}\|_{L^2(\Theta)}\approx n_0$ (see \cite{Tz}), we have
\begin{equation}\label{jm2}
\|\widetilde{\Delta}_{N_0,L_0}(v)\|_{L^2(\R\times\Theta)}
\approx
N_{0}^{-1}
\|\Delta_{N_0,L_0}(v)\|_{L^2(\R\times\Theta)}.
\end{equation}
Since $e_n$ vanishes on the boundary, using that
$$
e_{n}(r)=-\frac{1}{z_n^2}\frac{1}{r}\partial_{r}(r\partial_{r}e_{n}(r)),
$$
an integration by parts gives
$$
I(L_0,N_0,N_1)=\Big|\int_{\R\times\Theta}\overline{\widetilde{\Delta}_{N_0,L_0}(v)}
\partial_{r}\Big(\Delta_{N_1}(u)G_{1}(\Delta_{N_1}(u),\tilde{S}_{N_1}(u))\Big)\Big|.
$$
Recall that equality (\ref{G1}) shows that $G_{1}(\Delta_{N_1}(u),\tilde{S}_{N_1}(u))$ can be expanded as a
sum of two terms and then each term can be expanded according to (\ref{G1j}).
Therefore
$$
I(L_0,N_0,N_1)\leq I_1(L_0,N_0,N_1)+I_2(L_0,N_0,N_1)+I_3(L_0,N_0,N_1)+I_4(L_0,N_0,N_1),
$$
where
\begin{multline*}
I_1(L_0,N_0,N_1)=\sum_{j_1=1}^{2}\sum_{j_2=1}^{2}
\sum_{\stackrel{N_3\leq N_2\leq N_1}{N_2,N_3-{\rm dyadic }}}
\\
\int_{\R\times\Theta}
|\widetilde{\Delta}_{N_0,L_0}(v)\partial_{r}\big(\Delta_{N_1}(u)\big)\Delta_{N_2}(u)\Delta_{N_3}(u)
G_{1j_1j_2}^{N_3}(\Delta_{N_3}(u),\tilde{S}_{N_3}(u))|,
\end{multline*}
\begin{multline*}
I_2(L_0,N_0,N_1)=\sum_{j_1=1}^{2}\sum_{j_2=1}^{2}
\sum_{\stackrel{N_3\leq N_2\leq N_1}{N_2,N_3-{\rm dyadic }}}
\\
\int_{\R\times\Theta}
|\widetilde{\Delta}_{N_0,L_0}(v)\Delta_{N_1}(u)\partial_{r}\big(\Delta_{N_2}(u)\big)\Delta_{N_3}(u)
G_{1j_1j_2}^{N_3}(\Delta_{N_3}(u),\tilde{S}_{N_3}(u))|,
\end{multline*}
\begin{multline*}
I_3(L_0,N_0,N_1)=\sum_{j_1=1}^{2}\sum_{j_2=1}^{2}
\sum_{\stackrel{N_3\leq N_2\leq N_1}{N_2,N_3-{\rm dyadic }}}
\\
\int_{\R\times\Theta}|\widetilde{\Delta}_{N_0,L_0}(v)\Delta_{N_1}(u)\Delta_{N_2}(u)\partial_{r}
\big(\Delta_{N_3}(u)\big)G_{1j_1j_2}^{N_3}(\Delta_{N_3}(u),\tilde{S}_{N_3}(u))|,
\end{multline*}
\begin{multline*}
I_4(L_0,N_0,N_1)=\sum_{j_1=1}^{2}\sum_{j_2=1}^{2}
\sum_{\stackrel{N_3\leq N_2\leq N_1}{N_2,N_3-{\rm dyadic }}}
\\
\int_{\R\times\Theta}|\widetilde{\Delta}_{N_0,L_0}(v)\Delta_{N_1}(u)\Delta_{N_2}(u)\Delta_{N_3}(u)
\partial_{r}\big(G_{1j_1j_2}^{N_3}(\Delta_{N_3}(u),\tilde{S}_{N_3}(u))\big)|,
\end{multline*}
Recall that $G_{1j_1j_2}^{N_3}(z_1,z_2)$ satisfies the bound (\ref{novo}).
If we define $Q$ by (\ref{Q}), 
expanding with respect to the $L$ localizations,
using two times Proposition~\ref{str2} 
to the products $\widetilde{\Delta}_{N_0,L_0}(v)\Delta_{N_2,L_2}(u)$
and $\partial_{r}\big(\Delta_{N_1}(u)\big)\Delta_{N_3,L_3}(u)$
and (\ref{novo}) (together with (\ref{infty}), (\ref{jm}) and (\ref{jm2})) gives 
\begin{multline*}
\sum_{\stackrel{N_1\leq \min(N_0, M/2)}{L_0,N_0,N_1-{\rm dyadic }}}
I_1(L_0,N_0,N_1)
\leq
\\
\sum_{L_0,L_1,L_2,L_3-{\rm dyadic }}\,
\sum_{\stackrel{N_1\geq N_2\geq N_3, N_1\leq N_0}{N_0,N_1,N_2,N_3-{\rm dyadic }}}
L_{0}^{\beta-b'}(L_1L_2L_3)^{\beta-b}\Big(\frac{N_0}{N_1}\Big)^{\sigma-1}\,\,
\frac{N_{3}^{(1-\sigma_1)(\alpha-2)}}{(N_2N_3)^{\sigma_{1}-\varepsilon}}Q.
\end{multline*}
The last expression may estimated exactly as we did for $I_1$, by exchanging the roles of $N_0$ and $N_1$.
Similarly
\begin{multline*}
\sum_{\stackrel{N_1\leq \min(N_0, M/2)}{L_0,N_0,N_1-{\rm dyadic }}}
I_2(L_0,N_0,N_1)
\leq
\\
\sum_{L_0,L_1,L_2,L_3-{\rm dyadic }}\,
\sum_{\stackrel{N_1\geq N_2\geq N_3, N_1\leq N_0}{N_0,N_1,N_2,N_3-{\rm dyadic }}}
L_{0}^{\beta-b'}(L_1L_2L_3)^{\beta-b}\Big(\frac{N_0}{N_1}\Big)^{\sigma}
\Big(\frac{N_2}{N_0}\Big)
\,\,
\frac{N_{3}^{(1-\sigma_1)(\alpha-2)}}{(N_2N_3)^{\sigma_{1}-\varepsilon}}Q.
\end{multline*}
On the other hand on the summation region,
$$
\Big(\frac{N_0}{N_1}\Big)^{\sigma}
\Big(\frac{N_2}{N_0}\Big)\leq 
\Big(\frac{N_0}{N_1}\Big)^{\sigma-1}
$$
and thus, again, we may conclude as in the bound for $I_1$.
The sum
$$
\sum_{\stackrel{N_1\leq \min(N_0, M/2)}{L_0,N_0,N_1-{\rm dyadic }}}
I_3(L_0,N_0,N_1)
$$
can be bounded similarly.
Let us finally estimate the quantity
$$
\sum_{\stackrel{N_1\leq \min(N_0, M/2)}{L_0,N_0,N_1-{\rm dyadic }}}
I_4(L_0,N_0,N_1)\,.
$$
Observe that we can write 
\begin{multline*}
\Big|\partial_{r}\big(G_{1j_1j_2}^{N_3}(\Delta_{N_3}(u),\tilde{S}_{N_3}(u))\big)\Big|\leq
\\
C
\Big(
|\partial_{r}\big(\Delta_{N_3}(u)\big)|+|\partial_{r}\big(\tilde{S}_{N_3}(u)\big)|
\Big)
\Big(1+|\Delta_{N_3}(u)|+|\tilde{S}_{N_3}(u)|\Big)^{\max(\alpha-3,0)}\,.
\end{multline*}
Now using (\ref{infty-bis}), we can write
\begin{eqnarray*}
\|\partial_{r}(\Delta_{N_3}(u))\|_{L^{\infty}(\R\times\Theta)}
+
\|\partial_{r}(\tilde{S}_{N_3}(u))\|_{L^{\infty}(\R\times\Theta)}
& \leq &  
\sum_{\stackrel{N\leq N_3}{N-{\rm dyadic }}}
\|\partial_{r}(\Delta_{N}(u))\|_{L^{\infty}(\R\times\Theta)}
\\
& \leq &
\sum_{\stackrel{L,N\leq N_3}{L,N-{\rm dyadic }}}\|\partial_{r}(\Delta_{N,L}(u))\|_{L^{\infty}(\R\times\Theta)}
\\
& \leq &
\sum_{\stackrel{L,N\leq N_3}{L,N-{\rm dyadic }}}
CN^{2}L^{\frac{1}{2}}\|\Delta_{N,L}(u)\|_{L^{2}(\R\times\Theta)}
\\
& \leq &
CN_{3}^{2-\sigma_1}\|u\|_{X^{\sigma_1,b}_{rad}(\R\times\Theta)}\,.
\end{eqnarray*}
Similarly
$$
\Big(1+|\Delta_{N_3}(u)|+|\tilde{S}_{N_3}(u)|\Big)^{\max(\alpha-3,0)}
\leq 
C(1+N_{3}^{1-\sigma_1}\|u\|_{X^{\sigma_1,b}_{rad}(\R\times\Theta)})^{\max(\alpha-3,0)}\,\,.
$$
Let us suppose that $\alpha\in [3,4[$, the analysis for $\alpha\in [2,3]$ being simpler
(one needs to modify slightly the next several lines by invoking the
assumption $\sigma_1>1/3$).
If we define $Q$ by (\ref{Q}), 
expanding with respect to the $L$ localizations,
using Proposition~\ref{str1} to the product
$\Delta_{N_1,L_1}(u)\Delta_{N_3,L_3}(u)$ and Proposition~\ref{str2}
to the product $\widetilde{\Delta}_{N_0,L_0}(v)\Delta_{N_2,L_2}(u)$, we get
\begin{multline*}
\sum_{\stackrel{N_1\leq \min(N_0, M/2)}{L_0,N_0,N_1-{\rm dyadic }}}
I_4(L_0,N_0,N_1)
\leq
\\
\sum_{L_0,L_1,L_2,L_3-{\rm dyadic }}\,
\sum_{\stackrel{N_1\geq N_2\geq N_3, N_1\leq N_0}{N_0,N_1,N_2,N_3-{\rm dyadic }}}
L_{0}^{\beta-b'}(L_1L_2L_3)^{\beta-b}\Big(\frac{N_0}{N_1}\Big)^{\sigma}\,\,
\frac{1}{N_0}
\frac{N_3N_{3}^{(1-\sigma_1)(\alpha-2)}}{(N_2N_3)^{\sigma_{1}-\varepsilon}}Q.
\end{multline*}
Since on the region of summation
$$
\Big(\frac{N_0}{N_1}\Big)^{\sigma}\frac{1}{N_0}N_3\leq \Big(\frac{N_0}{N_1}\Big)^{\sigma-1}
$$
we may conclude exactly as we did for $I_1$.
This completes the analysis for $I_2$ and thus (\ref{main1}) is established.
Thanks to the multilinear nature of our reasoning (compare with the method of Ginibre-Velo, Kato for treating 
the Cauchy problem for NLS which is not multilinear), 
the proof of (\ref{main2}) is essentially the same as the proof of
(\ref{main1}). 
However one can no longer assume that the frequencies $N_1$ and $N_2$
satisfy $N_1\geq N_2$ but this fact does not affect the analysis since in
contrast with (\ref{main1}) all terms in the right hand-side of (\ref{main2})
have {\it the same} spatial regularity $\sigma$ (this is a
standard feature in the analysis of nonlinear PDE's and not related to the
spaces $X^{\sigma,b}_{rad}$ we work with).
More precisely, we can write
$$
F(u)-F(v)=(u-v)G_{1}(u,v)+(\overline{u}-\overline{v})G_{2}(u,v)
$$
with $G_{j}(z_1,z_2)$, $j=1,2$ satisfying the growth assumption
\begin{equation}\label{rastbis}
\big|
\partial^{k_1}_{z_1}\bar{\partial}^{k_2}_{z_1}
\partial^{l_1}_{z_2}\bar{\partial}^{l_2}_{z_2}
G_j(z_1,z_2)
\big|\leq C_{k_1,k_2,l_1,l_2}(1+|z_1|+|z_2|)^{\alpha-k_1-k_2-l_1-l_2}
\,.
\end{equation}
Since the analysis is very similar to the proof of (\ref{main1}), we shall
only outline the estimate for $(u-v)G_{1}(u,v)$. Again, we can suppose that
$F(u)$ is vanishing at order $3$ at $u=0$ and $\alpha\geq 2$. Let us set
$$
w_1=u-v,\quad w_2=u,\quad w_3=v.
$$
We thus need to prove that
\begin{multline*}
\Big|
\int_{\R\times\Theta}
w_1G_{1}(w_2,w_3)\overline{w_4}
\Big|
\\
\leq
C(1+\|w_2\|_{X^{\sigma,b}_{rad}(\R\times\Theta)}+\|w_3\|_{X^{\sigma,b}_{rad}(\R\times\Theta)})^\alpha
\|w_1\|_{X^{\sigma,b}_{rad}(\R\times\Theta)}\|w_4\|_{X^{-\sigma,b'}_{rad}(\R\times\Theta)}\,.
\end{multline*}
Next, we expand
$$
w_{1}=\sum_{N_1-{\rm dyadic}}\Delta_{N_1}(w_1),\quad w_{4}=\sum_{N_0-{\rm dyadic}}\Delta_{N_0}(w_4)
$$
and
$$
G_{1}(w_2,w_3)=\sum_{N_2-{\rm dyadic}}\Big(G_{1}(\tilde{S}_{2N_2}(w_2),\tilde{S}_{2N_2}(w_3))-
G_{1}(\tilde{S}_{N_2}(w_2),\tilde{S}_{N_2}(w_3))\Big).
$$
Thus, modulo complex conjugations irrelevant in this discussion, one has to
evaluate quantities of type
\begin{multline}\label{derm}
\sum_{N_0,N_1,N_2-{\rm dyadic}}
\Big|
\int_{\R\times\Theta}\overline{\Delta_{N_0}(w_4)}\Delta_{N_1}(w_1)\Delta_{N_2}(w_j)
\\
H_{j}^{N_2}(\Delta_{N_2}(w_2),\tilde{S}_{N_2}(w_2),\Delta_{N_2}(w_3),\tilde{S}_{N_2}(w_3))
\Big|,\quad j=2,3,
\end{multline}
where $H_{j}^{N_2}(z_1,z_2,z_3,z_4)$ are smooth functions satisfying growth
restrictions at infinity coming from (\ref{rast}). In the analysis of
(\ref{derm}), we distinguish two cases for $N_0$, $N_1$, $N_2$ in the sum
defining (\ref{derm}). Since $N_1$ and $N_2$ are not ordered, we need to
compare $N_0$ with $\max(N_1,N_2)$ by performing arguments close in the spirit
to the proof of (\ref{main1}).

{\bf Case 1.}
The first case is when $N_{0}\leq \max(N_1,N_2)$. In
this case, we expand once more $H_{j}^{N_2}$ which introduces a sum over $N_{3}-{\rm
dyadic}$, $N_{3}\leq N_2$ of terms $\Delta_{N_3}(w_{j})$ (or complex conjugate)
times a function which satisfies a decay property coming form (\ref{rast}).
As in the analysis of $I_1$ above, we obtain the bound
\begin{multline}\label{greve}
|(\ref{derm})|\leq
\sum_{L_0,L_1,L_2,L_3-{\rm dyadic }}\,
\sum_{\stackrel{N_1, N_2\geq N_3, \max(N_1,N_2)\geq N_0}{N_0,N_1,N_2,N_3-{\rm dyadic }}}
\\
L_{0}^{\beta-b'}(L_1L_2L_3)^{\beta-b}\Big(\frac{N_0}{\max(N_1,N_2)}\Big)^{\sigma}\,\,
\frac{N_{3}^{(1-\sigma)(\alpha-2)}}{(\min(N_1,N_2)N_3)^{\sigma-\varepsilon}}Q,
\end{multline}
where $Q$ is defined similarly to (\ref{Q}) with the important difference that
$\sigma_1$ is replaced by $\sigma$ and the harmless difference that $u$ is
replaced by a suitable $w_j$, $j=1,2,3$ and $v$ is replaced by $w_4$.
If $\max(N_1,N_2)=N_1$ or $N_1\geq N_3$ then we conclude exactly as in the proof of
(\ref{main1}). 

We can therefore suppose that $\max(N_1,N_2)=N_2$ and $N_1\leq N_3$.
Observe that we can also suppose that $F(u)$ is vanishing at order $5$ at
$u=0$ which allows to expand the non-linearity once more. Indeed, the cubic
term in the Taylor expansion of the non-linearity can be dealt with as in
(\ref{greve}) since in this term $\alpha=2$. Thus in the case $\max(N_1,N_2)=N_2$ and $N_1\leq N_3$, 
we expand once more the non-linearity which introduces a sum over $N_{4}-{\rm
dyadic}$, $N_{4}\leq N_3$ of terms $\Delta_{N_4}(w_{j})$ (or complex conjugate)
times a function which satisfies an appropriate decay property coming form
(\ref{rast}).
We next consider two cases $N_1\geq N_4$ and $N_1\leq N_4$.
Let us suppose first that $N_1\geq N_4$. In this case,
using the bilinear Strichartz estimates as in the analysis of $I_1$ above, we obtain the bound
\begin{multline}\label{greve2}
|(\ref{derm})|\leq
\sum_{L_0,L_1,L_2,L_3,L_4-{\rm dyadic }}\,
\sum_{\stackrel{N_3\geq N_1\geq N_4, N_2\geq N_3\geq N_4, N_2\geq N_0}{N_0,N_1,N_2,N_3,N_4-{\rm dyadic }}}
\\
L_{0}^{\beta-b'}(L_1L_2L_3L_4)^{\beta-b}\Big(\frac{N_0}{N_2}\Big)^{\sigma}\,\,
\frac{N_4^{(1-\sigma)(\alpha-2)}}
{(N_1 N_3)^{\sigma-\varepsilon}}Q,
\end{multline}
where $Q$ is defined similarly to (\ref{Q}) with one additional factor in the
product, i.e. the product runs from $1$ to $4$ instead of $1$ to $3$.
With (\ref{greve2}) at our disposal, we can conclude exactly as in the proof
of (\ref{main1}).
Let us suppose finally that $N_1\leq N_4$. In this case we put the term
involving $\Delta_{N_1}$ in $L^\infty$ and perform the bilinear estimates with
the terms involving $N_0$, $N_2$, $N_3$, $N_4$ to get 
\begin{multline*}
|(\ref{derm})|\leq
\sum_{L_0,L_1,L_2,L_3,L_4-{\rm dyadic }}\,
\sum_{\stackrel{N_1\leq N_4, N_2\geq N_3\geq N_4, N_2\geq N_0}{N_0,N_1,N_2,N_3,N_4-{\rm dyadic }}}
\\
L_{0}^{\beta-b'}(L_1L_2L_3L_4)^{\beta-b}\Big(\frac{N_0}{N_2}\Big)^{\sigma}\,\,
\frac{N_1^{1-\sigma}N_4^{\max(0,(1-\sigma)(\alpha-3))}}
{(N_3 N_4)^{\sigma-\varepsilon}}Q,
\end{multline*}
where $Q$ is defined similarly to (\ref{Q}).
Once again we conclude similarly to the proof
of (\ref{main1}).

{\bf Case 2.}
If $N_{0}\geq \max(N_1,N_2)$, then we integrate by parts by the aid of
$\Delta_{N_0}(w_4)$ and the analysis is very similar to the bound for 
$I_2$ in the proof of (\ref{main1}).
\\

This completes the proof of Proposition~\ref{main}.
\end{proof}
\begin{remarque}
We refer to \cite{Ramona}, where an analysis similar to the proof of Proposition~\ref{main} is performed.
In \cite{Ramona}, one proves bilinear Strichartz estimates for the free evolution and by the transfer principal
of \cite{BGT} these estimates are transformed to estimates involving the projector $\Delta_{N,L}$.
This approach is slightly different from the approach used in \cite{Tz}, based on direct bilinear estimates for
functions enjoying localization properties similar to $\Delta_{N,L}(u)$.
\end{remarque}
\section{Local analysis for NLS and the approximating ODE}
In this section, we state the standard consequence of Proposition~\ref{main} to the local well-posedness of
(\ref{1}) and (\ref{N}).
For $T>0$, we define the restriction spaces
$X^{\sigma,b}_{rad}([-T,T]\times\Theta)$, equipped with the natural norm
$$
\|u\|_{X^{\sigma,b}_{rad}([-T,T]\times\Theta)}=
\inf\{
\|w\|_{X^{\sigma,b}_{rad}(\R\times\Theta)},\quad w\in
X^{\sigma,b}_{rad}(\R\times\Theta)\quad {\rm with}\quad w|_{]-T,T[}=u 
\}.
$$
Similarly, for $I\subset\R$ an interval, we can define the restriction spaces
$X^{\sigma,b}_{rad}(I\times\Theta)$, equipped with the natural norm.
A Sobolev inequality with respect to the time variable yields,
\begin{equation*}
\|u\|_{L^{\infty}([-T,T]\,;\,H^{\sigma}_{rad}(\Theta))}\leq C_{b}
\|u\|_{X^{\sigma,b}_{rad}([-T,T]\times \Theta)},\quad b>\frac{1}{2}.
\end{equation*}
Thus for $b>1/2$ the space $X^{\sigma,b}_{rad}([-T,T]\times \Theta)$ 
is continuously embedded in $C([-T,T]\,;\,H^{\sigma}_{rad}(\Theta))$.
We shall solve (\ref{1}) for short times by applying the Banach contraction mapping principle to the 
``Duhamel formulation'' of (\ref{1})
\begin{equation}\label{venda}
u(t)=e^{it\Delta}u_0-i\int_{0}^{t}e^{i(t-\tau)\Delta}F(u(\tau))d\tau\,,
\end{equation}
where $e^{it\Delta}$ denotes the free propagator.
\begin{remarque}\label{zabelejka}
{\rm
In (\ref{venda}), the operator $e^{it\Delta}$ is defined by the Dirichlet self-adjoint realization of the
Laplacian via the functional calculus of self-adjoint operators. As mentioned before the uniqueness
statements in the well-posedness results in this paper are understood as uniqueness results for
(\ref{venda}).
On the other hand, despite the low regularity situation in this paper, 
the solutions  of (\ref{venda}), we construct here have zero traces on $\R\times\partial\Theta$
(which is a general feature reflecting from the Dirichlet Laplacian, we work with)
and thus the uniqueness issue can be studied in the context of the equation (\ref{1})
subject to zero boundary conditions on $\R\times\partial\Theta$.
If we set $S(t)=e^{it\Delta}$, then $S(t)e_{n}=e^{-itz_n^2}e_{n}$ and the norms in the Bourgain spaces
may be expressed as 
$$
\|u\|_{X^{\sigma,b}_{rad}(\R\times\Theta)}=\|S(-t)u\|_{H^{\sigma,b}_{rad}(\R\times\Theta)},
$$
where $H^{\sigma,b}_{rad}(\R\times\Theta)$ is a classical anisotropic Sobolev space equipped with the norm
$$
\|v\|_{H^{\sigma,b}_{rad}(\R\times\Theta)}^{2}=\sum_{n\geq 1}z_{n}^{2\sigma}\|
\langle\tau\rangle^{b}\widehat{\langle v(t),e_{n}\rangle}(\tau)\|_{L^{2}(\R_{\tau})}^{2}\,,
$$
where again $\langle\cdot,\cdot\rangle$ stays for the $L^2(\Theta)$ pairing and $\widehat{\cdot}$ denotes 
the Fourier transform on $\R$.
Therefore in the context of (\ref{venda}) we are in a situation where the Bourgain approach to well-posedness
of dispersive equations may be applied.
Let us also observe that the solutions of (\ref{venda}) we obtain here 
solve (\ref{1}) in distributional sense (see e.g. \cite[Section~3.2]{BGT} for details on this point).
Let us finally remark that for $\sigma<1/2$ the spaces $H^{\sigma}_{rad}(\Theta)$ are independent of 
the choice of the boundary conditions we work with. In particular, the space $H^s_{rad}(\Theta)$, on 
which the invariant measure $d\rho$ is defined, is independent of the boundary conditions.
On the other hand both the dynamics and the Gibbs measure $d\rho$ do depend on
the choice of the boundary conditions.
}
\end{remarque}
Now we state the following standard consequence of Proposition~\ref{main} (see \cite{Gi} or
\cite[Proposition~6.3]{Tz}).
\begin{proposition}\label{duh}
Let $\max(1/3,1-2/\alpha)<\sigma_1\leq \sigma<1/2$. Then there exist two positive numbers $b,b'$ such
that $b+b'<1$, $b'<1/2<b$, there exists $C>0$ such that for every $T\in]0,1]$, 
every $u,v\in X^{\sigma,b}_{rad}([-T,T]\times\Theta)$, every $u_0\in
H^{\sigma}_{rad}(\Theta)$, 
\begin{equation*}
\big\|e^{it\Delta}u_0 \big\|_{X^{\sigma,b}_{rad}([-T,T]\times\Theta)}\leq 
C\|u_0\|_{H^{\sigma}_{rad}(\Theta)}\, ,
\end{equation*}
\begin{multline*}
\Big\|\int_{0}^{t}e^{i(t-\tau)\Delta}F(u(\tau))d\tau\Big\|_{X^{\sigma,b}_{rad}([-T,T]\times\Theta)}
\leq
\\
\leq CT^{1-b-b'}\Big(1+\|u\|^{\max(2,\alpha)}_{X^{\sigma_{1},b}_{rad}([-T,T]\times\Theta)}\Big)
\|u\|_{X^{\sigma,b}_{rad}([-T,T]\times\Theta)}
\end{multline*}
and
\begin{multline*}
\Big\|\int_{0}^{t}e^{i(t-\tau)\Delta}(F(u(\tau))-F(v(\tau)))d\tau
\Big\|_{X^{\sigma,b}_{rad}([-T,T]\times\Theta)}
\leq
\\
\leq CT^{1-b-b'}\Big(1+\|u\|^{\max(2,\alpha)}_{X^{\sigma,b}_{rad}([-T,T]\times\Theta)}+
\|v\|^{\max(2,\alpha)}_{X^{\sigma,b}_{rad}([-T,T]\times\Theta)}\Big)
\|u-v\|_{X^{\sigma,b}_{rad}([-T,T]\times\Theta)}\,.
\end{multline*}
\end{proposition}
One may also formulate statements in the spirit of Proposition~\ref{duh}, where $[-T,T]$ is replaced by
an interval $I\subset\R$ of size one and $0$ by a point of $I$.
We also remark that the integral terms in  Proposition~\ref{duh} are well-defined thanks to a priori
estimates in the Bourgain spaces (see e.g. \cite{Gi}).

Proposition~\ref{duh} implies (see \cite[Proposition~7.1]{Tz}) 
a local well-posedness result for the Cauchy problem
\begin{equation}\label{1bis}
(i\partial_t+\Delta)u-F(u)=0,\quad u|_{t=0}=u_0.
\end{equation}
\begin{proposition}\label{lwp}
Let us fix $\sigma_1$ and $\sigma$ such that $\max(1/3,1-2/\alpha)<\sigma_{1}\leq\sigma<1/2$.
Then there exist $b>1/2$, $\beta>0$, $C>0$, $\tilde{C}>0$, $c\in]0,1]$ such that for every 
$A>0$ if we set $T=c(1+A)^{-\beta}$ then for every $u_0\in H^{\sigma_{1}}_{rad}(\Theta)$ 
satisfying $\|u_0\|_{H^{\sigma_{1}}}\leq A$ there exists a unique
solution $u$ of (\ref{venda}) in $X^{\sigma_{1},b}_{rad}([-T,T]\times \Theta)$. 
Moreover $u$ solves (\ref{1bis}) and
$$
\|u\|_{L^{\infty}([-T,T];H^{\sigma_{1}}(\Theta))}
\leq 
C\|u\|_{X^{\sigma_{1},b}_{rad}([-T,T]\times \Theta)}
\leq \tilde{C}\|u_0\|_{H^{\sigma_{1}}(\Theta)}\, .
$$ 
If in addition $u_0\in H^{\sigma}_{rad}(\Theta)$ then
$$
\|u\|_{L^{\infty}([-T,T];H^{\sigma}(\Theta))}\leq 
C\|u\|_{X^{\sigma,b}_{rad}([-T,T]\times \Theta)}
\leq \tilde{C}\|u_0\|_{H^{\sigma}(\Theta)}\, .
$$ 
Finally if $u$ and $v$ are two solutions with data $u_0$, $v_0$ respectively,
satisfying 
$$
\|u_0\|_{H^{\sigma_{1}}}\leq A,\quad \|v_0\|_{H^{\sigma_{1}}}\leq A
$$
then 
$$
\|u-v\|_{L^{\infty}([-T,T];H^{\sigma_{1}}(\Theta))}\leq C\|u_0-v_0\|_{H^{\sigma_{1}}(\Theta)}\, .
$$ 
If in addition $u_0,v_0\in H^{\sigma}_{rad}(\Theta)$ then
$$
\|u-v\|_{L^{\infty}([-T,T];H^{\sigma}(\Theta))}\leq C\|u_0-v_0\|_{H^{\sigma}(\Theta)}\, .
$$
\end{proposition}
Since the projector $S_N$ is acting nicely on the Bourgain spaces Proposition~\ref{duh} also implies a 
well-posedness result (the important point is the independence of $N$ of the constants appearing 
in the statement) for the ODE (\ref{N}).
\begin{proposition}\label{lwpbis}
Let us fix $\sigma_1$ and $\sigma$ such that $\max(1/3,1-2/\alpha)<\sigma_{1}\leq\sigma<1/2$.
Then there exist $b>1/2$, $\beta>0$, $C>0$, $\tilde{C}>0$, $c\in]0,1]$ such that for every $A>0$ if we set
$T= c(1+A)^{-\beta}$ then for every $N\geq 1$, every 
$u_0\in E_{N}$ satisfying $\|u_0\|_{H^{\sigma_{1}}}\leq A$ there exists a unique
solution $u=S_{N}(u)$ of 
(\ref{N}) in $X^{\sigma_{1},b}_{rad}([-T,T]\times \Theta)$. Moreover 
\begin{equation*}
\|u\|_{L^{\infty}([-T,T];H^{\sigma_{1}}(\Theta))}
\leq 
C\|u\|_{X^{\sigma_{1},b}_{rad}([-T,T]\times \Theta)}
\leq \tilde{C}\|u_0\|_{H^{\sigma_{1}}(\Theta)}\, .
\end{equation*}
If in addition $u_0\in H^{\sigma}_{rad}(\Theta)$ then
\begin{equation*}
\|u\|_{L^{\infty}([-T,T];H^{\sigma}(\Theta))}
\leq 
C\|u\|_{X^{\sigma,b}_{rad}([-T,T]\times \Theta)}
\leq \tilde{C}\|u_0\|_{H^{\sigma}(\Theta)}\, .
\end{equation*}
Finally if $u$ and $v$ are two solutions with data $u_0$, $v_0$ respectively,
satisfying 
$$
\|u_0\|_{H^{\sigma_{1}}}\leq A,\quad \|v_0\|_{H^{\sigma_{1}}}\leq A
$$
then 
$$
\|u-v\|_{L^{\infty}([-T,T];H^{\sigma_{1}}(\Theta))}\leq C\|u_0-v_0\|_{H^{\sigma_{1}}(\Theta)}\, .
$$ 
If in addition $u_0,v_0\in H^{\sigma}_{rad}(\Theta)$ then
$$
\|u-v\|_{L^{\infty}([-T,T];H^{\sigma}(\Theta))}\leq 
C
\|u_0-v_0\|_{H^{\sigma}(\Theta)}\, .
$$
\end{proposition}
\section{Long time analysis of the approximating ODE}
In this section we study the long time dynamics of 
\begin{equation}\label{Nbis}
(i\partial_t+\Delta)u-S_{N}(F(u))=0,\quad u|_{t=0}\in E_N\,.
\end{equation}
Recall from the introduction that the measure $d\rho_N$ is invariant under the well-defined flow 
of (\ref{Nbis}). Denote this flow by $\Phi_{N}(t):E_{N}\rightarrow E_{N}$, $t\in\R$.
We have the following statement.
\begin{proposition}\label{longtime}
There exists $\Lambda>0$ such that for every integer $i\geq 1$, every $\sigma\in[s,1/2[$, every $N\in\N$, 
there exists a $\rho_{N}$ measurable set $\Sigma_{N,\sigma}^{i}\subset E_{N}$ such that : 
\begin{itemize}
\item
$
\rho_{N}(E_{N}\backslash \Sigma_{N,\sigma}^{i})\leq 2^{-i}\,.
$
\item
$
u\in \Sigma_{N,\sigma}^{i}\,\Rightarrow\, \|u\|_{L^2(\Theta)}\leq \Lambda.
$
\item
There exists $C_{\sigma}$, depending on $\sigma$, such that for every $i\in\N$, 
every $N\in\N$, every $u_0\in \Sigma_{N,\sigma}^{i}$, every $t\in\R$, 
$
\|\Phi_{N}(t)(u_0)\|_{H^{\sigma}(\Theta)}\leq C_{\sigma}(i+\log(1+|t|))^{\frac{1}{2}}\,.
$
\item
For every $\sigma\in ]s,1/2[$, every $\sigma_1\in [s,\sigma[$, every $t\in\R$
there exists $i_{1}$ such that for every integer $i\geq 1$, every $N\geq 1$, if $u_0\in
\Sigma^{i}_{N,\sigma}$ then one has
$
\Phi_{N}(t)(u_0)\in \Sigma^{i+i_1}_{N,\sigma_1}.
$
\end{itemize}
\end{proposition}
\begin{remarque}
One may wish to see the invariance property of the sets $\Sigma^{i}_{N,\sigma}$ displayed 
by the last assertion as a ``weak form of a conservation law''.
\end{remarque}
\begin{proof}
Let $\Lambda>0$ be such that $\chi(x)=0$ for $|x|>\Lambda$.
For $\sigma\in [s,1/2[$, $i,j$ integers $\geq 1$, we set
$$
B_{N,\sigma}^{i,j}(D_{\sigma})\equiv
\Big\{u\in E_{N}\,:\,\|u\|_{H^{\sigma}(\Theta)}\leq D_{\sigma}(i+j)^{\frac{1}{2}},\quad
\|u\|_{L^2(\Theta)}\leq \Lambda\Big\},
$$
where the number $D_{\sigma}\gg 1$ (independent of $i,j,N$) will be fixed later. 
Thanks to Proposition~\ref{lwpbis}, there exist $c>0$, $C>0$, $\beta>0$ only depending on $\sigma$ 
such that if we set $\tau \equiv c D_{\sigma}^{-\beta}(i+j)^{-\beta/2}$ then for every $t\in[-\tau,\tau]$,
\begin{equation}\label{preser}
\Phi_{N}(t)\big(B_{N,\sigma}^{i,j}(D_{\sigma})\big)\subset 
\Big\{u\in E_{N}\,:\,\|u\|_{H^{\sigma}(\Theta)}\leq C\,D_{\sigma}(i+j)^{\frac{1}{2}},\quad
\|u\|_{L^2(\Theta)}\leq \Lambda\Big\}\, .
\end{equation}
Next, following Bourgain \cite{Bo1}, we set
$$
\Sigma_{N,\sigma}^{i,j}(D_{\sigma})\equiv
\bigcap_{k=-[2^{j}/\tau]}^{[2^{j}/\tau]}\Phi_{N}(-k\tau)(B_{N,\sigma}^{i,j}(D_{\sigma}))\, ,
$$
where $[2^{j}/\tau]$ stays for the integer part of $2^{j}/\tau$. 
Using the invariance of the measure $\rho_{N}$ by the flow $\Phi_{N}$, we can write
\begin{eqnarray*}
\rho_{N}(E_{N}\backslash\Sigma_{N,\sigma}^{i,j}(D_{\sigma}))
& = &
\rho_{N}\Big(\bigcup_{k=-[2^{j}/\tau]}^{[2^{j}/\tau]}
\big(E_{N}\backslash\Phi_{N}(-k\tau)\big(B_{N,\sigma}^{i,j}(D_{\sigma})\big)\big)\Big) 
\\
& \leq &
(2[2^{j}/\tau]+1)\rho_{N}(E_N\backslash B_{N,\sigma}^{i,j}(D_{\sigma}))
\\
& \leq &
C2^{j}D_{\sigma}^{\beta}(i+j)^{\beta/2}\rho_{N}(E_N\backslash B_{N,\sigma}^{i,j}(D_{\sigma}))\,.
\end{eqnarray*}
Using the support property of $\chi$, we observe that
set $(u\in E_{N}\,:\, \|u\|_{L^2(\Theta)}>\Lambda)$ is of zero $\rho_N$ measure and therefore 
\begin{eqnarray}
\rho_{N}(E_N\backslash B_{N,\sigma}^{i,j}(D_{\sigma}))  = \tilde{\rho}_{N}
\Big(u\in H^{s}_{rad}(\Theta)\,:\,\|S_{N}(u)\|_{H^{\sigma}(\Theta)}> D_{\sigma}(i+j)^{\frac{1}{2}} \Big).
\end{eqnarray}
Therefore, using Lemma~\ref{gauss-bis}, we can write
\begin{equation}\label{zvez}
\rho_{N}(E_{N}\backslash\Sigma_{N,\sigma}^{i,j}(D_{\sigma}))\leq
C 2^{j}D_{\sigma}^{\beta}(i+j)^{\beta/2}e^{-cD_{\sigma}^2(i+j)}\leq 2^{-(i+j)},
\end{equation}
provided $D_{\sigma}\gg 1$, depending on $\sigma$ but independent of $i,j,N$.
Thanks to (\ref{preser}), we obtain that for
$u_0\in\Sigma_{N,\sigma}^{i,j}(D_\sigma)$, the solution of (\ref{Nbis}) with data $u_0$ satisfies
\begin{equation}\label{jjj1}
\|\Phi_{N}(t)(u_0)\|_{H^{\sigma}(\Theta)}\leq CD_{\sigma}(i+j)^{\frac{1}{2}},\quad |t|\leq 2^{j}\,.
\end{equation}
Indeed, for $|t|\leq 2^{j}$, we may find an integer $k\in [-[2^{j}/\tau],[2^{j}/\tau]]$ and 
$\tau_1\in [-\tau,\tau]$ so that $t=k\tau+\tau_1$ and thus 
$u(t)=\Phi_{N}(\tau_1)\big(\Phi_{N}(k\tau)(u_0)\big)$.
Since $u_0\in\Sigma_{N,\sigma}^{i,j}(D_\sigma)$ implies that $\Phi_{N}(k\tau)(u_0)\in 
B_{N,\sigma}^{i,j}(D_\sigma)$, we may apply (\ref{preser}) and arrive at (\ref{jjj1}).
Next, we set
$$
\Sigma_{N,\sigma}^{i}=\bigcap_{j= 1}^{\infty}\Sigma_{N,\sigma}^{i,j}(D_{\sigma})\,.
$$
Thanks to (\ref{zvez}),
\begin{equation*}
\rho_{N}(E_{N}\backslash \Sigma_{N,\sigma}^{i})\leq 2^{-i}\,.
\end{equation*}
In addition, using (\ref{jjj1}), we get that there exists $C_{\sigma}$ such that for every $i$, every
$N$, every $u_0\in \Sigma_{N,\sigma}^{i}$, every $t\in \R$,
$$
\|\Phi_{N}(t)(u_0)\|_{H^{\sigma}(\Theta)}\leq C_{\sigma}(i+\log(1+|t|))^{\frac{1}{2}}\,.
$$
Indeed for $t\in \R$ there exists $j\in\N$ such that $2^{j-1}\leq 1+|t|\leq 2^j$ and we apply 
(\ref{jjj1}) with this $j$.

Let us now turn to the proof of the last assertion.
Fix $t\in\R$ and $u_{0}\in \Sigma_{N,\sigma}^{i}$. 
Since $u_{0}\in \Sigma_{N,\sigma}^{i}$, for every integer $j\geq 1$, we have the bound
$$
\|\Phi_{N}(t_1)(u_0)\|_{H^{\sigma}(\Theta)}\leq C_{\sigma}(i+j)^{\frac{1}{2}},\quad  |t_1|\leq 2^j.
$$
Let $i_1\in\N$ (depending on $t$) be such that for every $j\geq 1$, $2^{j}+|t|\leq 2^{j+i_1}$.
Therefore, if we set $u(t)\equiv\Phi_{N}(t)(u_0)$, we have that
$$
\|\Phi_{N}(t_1)(u(t))\|_{H^{\sigma}(\Theta)}
=
\|\Phi_{N}(t+t_1)(u_0)\|_{H^{\sigma}(\Theta)}
\leq
C_{\sigma}(i+j+i_1)^{\frac{1}{2}},\quad  |t_1|\leq 2^j.
$$
Thanks to  the $L^2$ conservation law, for $u_0\in\Sigma^{i}_{N,\sigma}$ one has
$$
\|\Phi_{N}(t_1)(u(t))\|_{L^{2}(\Theta)}=\|u_0\|_{L^2(\Theta)}\leq \Lambda.
$$
Therefore
\begin{eqnarray*}
\|\Phi_{N}(t_1)(u(t))\|_{H^{\sigma_1}(\Theta)}
& \leq &
\|\Phi_{N}(t_1)(u(t))\|_{H^{\sigma}(\Theta)}^{\frac{\sigma_1}{\sigma}}
\|\Phi_{N}(t_1)(u(t))\|_{L^{2}(\Theta)}^{\frac{\sigma-\sigma_1}{\sigma}}
\\
& \leq &
[\Lambda]^{\frac{\sigma-\sigma_1}{\sigma}}\Big[C_{\sigma}(i+j+i_1)\Big]^{\frac{\sigma_1}{2\sigma}}\,.
\end{eqnarray*}
Let us fix $i_1\geq 1$ such that in addition to the property
$$
2^{j}+|t|\leq 2^{j+i_1},\quad \forall\,\, j\geq 1,
$$
we also have that for every $i,j\geq 1$,
$$
[\Lambda]^{\frac{1-\sigma_1}{\sigma}}\Big[C_{\sigma}(i+j+i_1)\Big]^{\frac{\sigma_1}{2\sigma}}
\leq D_{\sigma_1}(i+j+i_1)^{\frac{1}{2}}\,.
$$
Thus
$$
\|\Phi_{N}(t_1)(u(t))\|_{H^{\sigma_1}(\Theta)}
\leq
D_{\sigma_1}(i+j+i_1)^{\frac{1}{2}}, \quad  |t_1|\leq 2^j,
$$
i.e. for every $|t_1|\leq 2^j$ one has 
$\Phi_{N}(t_1)(u(t))\in B_{N,\sigma_1}^{i+i_1,j}(D_{\sigma})$.
We can therefore conclude that $u(t)\in\Sigma_{N,\sigma_1}^{i+i_1,j}(D_{\sigma})$ for every $j\geq
1$. Hence 
$
u(t)\in \Sigma_{N,\sigma_1}^{i+i_1}
$
and the restriction on $i_1$ depends only on $\sigma$, $\sigma_1$ and $t$.
This completes the proof of Proposition~\ref{longtime}.
\end{proof}
\section{Construction of the statistical ensemble (long time analysis for NLS)}
Let us set for  integers $i\geq 1$, $N \geq 1$ and $\sigma\in[s,1/2[$, 
$$
\tilde{\Sigma}_{N,\sigma}^{i}\equiv
\big\{
u\in H^{\sigma}_{rad}(\Theta)\,:\, S_{N}(u)\in \Sigma_{N,\sigma}^{i}
\big\}.
$$
Next, for an integer $i\geq 1$ and $\sigma\in[s,1/2[$, we set
$$
\Sigma_{\sigma}^{i}\equiv 
\big\{
u\in H^{\sigma}_{rad}(\Theta)\,:\,\exists\, N_k\rightarrow\infty, N_k\in\N,\,
\exists\, u_{N_k}\in \Sigma_{N_k,\sigma}^{i},\,u_{N_k}\rightarrow u\,\, {\rm in}\, H^{\sigma}_{rad}(\Theta)
\big\}.
$$
We have the following statement.
\begin{lemme}\label{nicolas}
The set $\Sigma_{\sigma}^{i}$ is a closed set in $H^{\sigma}_{rad}(\Theta)$
(in particular $\rho$ measurable).
\end{lemme}
\begin{proof}
Let $(u_{m})_{m\in\N}$ be a sequence of $\Sigma_{\sigma}^{i}$ which converges to 
$u$ in $H^{\sigma}_{rad}(\Theta)$. Our goal is to show that $u\in \Sigma_{\sigma}^{i}$.
Since $u_m\in \Sigma_{\sigma}^{i}$ there exist a sequence of integers 
$N_{m,k}\rightarrow\infty$ as $k\rightarrow \infty$
and a sequence $(u_{N_{m,k}})_{k\in\N}$ of $\Sigma_{N_{m,k},\sigma}^{i}$ such that
\begin{equation}\label{noel}
\lim_{k\rightarrow\infty} u_{N_{m,k}}=u_{m}\quad \, {\rm in}\,\,\,\, H^{\sigma}(\Theta)\,.
\end{equation}
For every $j\in\N$, we can find $m_j\in \N$ such that
$$
\|u-u_{m_j}\|_{H^{\sigma}(\Theta)}<\frac{1}{2j}\,.
$$
Then, thanks to (\ref{noel}) (with $m=m_j$), for every $j\in\N$, we can find $N_{m_j,k_j}\in \N$ 
and $u_{N_{m_j,k_j}}\in\Sigma^{i}_{N_{m_j,k_j},\sigma}$ such that
$$
N_{m_j,k_j}>j,\quad \|u_{m_j}-u_{N_{m_j,k_j}}\|_{H^{\sigma}(\Theta)}<\frac{1}{2j}\,.
$$
Therefore, if we set
$v_{j}\equiv u_{N_{m_j,k_j}}$ and $M_{j}\equiv N_{m_j,k_j}$ then
$M_{j}\rightarrow \infty$ as $j\rightarrow\infty$, $v_{j}\in \Sigma^{i}_{M_j,\sigma}$
and $v_j\rightarrow u$ as $j\rightarrow\infty$ in $H^{\sigma}_{rad}(\Theta)$.
Consequently $u\in \Sigma_{\sigma}^{i}$.
This completes the proof of Lemma~\ref{nicolas}.
\end{proof}
We have the inclusion
$$
\limsup_{N\rightarrow\infty}\tilde{\Sigma}_{N,\sigma}^{i}
\equiv\bigcap_{N=
  1}^{\infty}\bigcup_{N_1=N}^{\infty}\tilde{\Sigma}_{N_1,\sigma}^{i}
\subset \Sigma_{\sigma}^{i}.
$$
Indeed, if $u\in \limsup_{N\rightarrow\infty}\tilde{\Sigma}_{N,\sigma}^{i}$
then there exists a sequence of integers $(N_k)$ tending to infinity as
$k\rightarrow\infty$ such that $u\in \tilde{\Sigma}_{N_k,\sigma}^{i}$, i.e.
$S_{N_k}(u)\in \Sigma_{N_k,\sigma}^{i}$. Thus $u\in  \Sigma_{\sigma}^{i}$
since $S_{N_k}(u)$ tends to $u$ in $H^{\sigma}_{rad}(\Theta)$. 
Therefore 
\begin{equation}\label{kr1}
\rho(\Sigma_{\sigma}^{i})   \geq   \rho(\limsup_{N\rightarrow\infty}\tilde{\Sigma}_{N,\sigma}^{i})\,.
\end{equation}
Let us next show that
\begin{equation}\label{kr2}
\rho(\limsup_{N\rightarrow\infty}\tilde{\Sigma}_{N,\sigma}^{i})
 \geq 
\limsup_{N\rightarrow \infty}\rho(\tilde{\Sigma}_{N,\sigma}^{i})\,.
\end{equation}
Indeed, if we set $A_{N}\equiv \tilde{\Sigma}_{N,\sigma}^{i}$ and 
$B_{N}\equiv H^s_{rad}(\Theta)\backslash A_{N}$ then
\begin{equation}\label{kr3}
\limsup_{N\rightarrow \infty}\rho(A_{N})=\limsup_{N\rightarrow \infty}\Big(\rho(H^s_{rad}(\Theta))-
\rho(B_{N})\Big)
=\rho(H^s_{rad}(\Theta))-\liminf_{N\rightarrow \infty}\rho(B_{N})\,.
\end{equation}
Using Fatou's lemma, we can obtain
\begin{equation*}
-\liminf_{N\rightarrow \infty}\rho(B_{N})\leq -\rho\Big(\liminf_{N\rightarrow \infty}B_{N}\Big),
\end{equation*}
where
$$
\liminf_{N\rightarrow \infty}B_{N}\equiv
\bigcup_{N= 1}^{\infty}\bigcap_{N_1=N}^{\infty}B_{N_1}\,.
$$
Therefore, coming back to (\ref{kr3}), we get
\begin{equation*}
\limsup_{N\rightarrow \infty}\rho(A_{N})\leq \rho\Big(
H^s_{rad}(\Theta)\backslash
\liminf_{N\rightarrow \infty}B_{N}
\Big)
=\rho\Big(\limsup_{N\rightarrow \infty}A_{N}\Big).
\end{equation*}
Therefore (\ref{kr2}) holds. Since
$$
\rho(\tilde{\Sigma}_{N,\sigma}^{i})=\int_{\tilde{\Sigma}_{N,\sigma}^{i}}f(u)d\mu(u)
$$
and
$$
\rho_{N}(\Sigma_{N,\sigma}^{i})=\int_{\Sigma_{N,\sigma}^{i}}f_{N}(u)d\mu_{N}(u)=\int_{\tilde{\Sigma}_{N,\sigma}^{i}}f_{N}(u)d\mu(u)
$$
thanks to Lemma~\ref{lem3}, we get
$$
\lim_{N\rightarrow \infty}\big((\rho(\tilde{\Sigma}_{N,\sigma}^{i})-\rho_{N}(\Sigma_{N,\sigma}^{i})\big)=0\,.
$$
Thus, using  Proposition~\ref{longtime} and Theorem~\ref{thm2}, we obtain
\begin{equation}\label{kr4}
\limsup_{N\rightarrow \infty}\rho(\tilde{\Sigma}_{N,\sigma}^{i})
=
\limsup_{N\rightarrow \infty}\rho_{N}(\Sigma_{N,\sigma}^{i})
\geq
\limsup_{N\rightarrow \infty}\big(\rho_{N}(E_{N})-2^{-i}\big)
=
\rho\big(H^s_{rad}(\Theta)\big)-2^{-i}.
\end{equation}
Collecting (\ref{kr1}), (\ref{kr2}) and (\ref{kr4}), we arrive at
$$
\rho(\Sigma_{\sigma}^{i})  \geq \rho\big(H^s_{rad}(\Theta)\big)-2^{-i}.
$$
Now, we set
$$
\Sigma_{\sigma}\equiv\bigcup_{i\geq 1}\Sigma_{\sigma}^{i}\,.
$$
Thus $\Sigma_{\sigma}$ is of full $\rho$ measure.
It turns out that one has global existence for $u_0\in \Sigma_{\sigma}^{i}$. 
\begin{proposition}\label{global_existence}
Let us fix $\sigma\in [s,1/2[$, $\sigma_1\in ]\max(1/3,1-2/\alpha),\sigma[$ and $i\in\N$.
Then for every $u_0\in \Sigma_{\sigma}^{i}$, the local solution $u$ 
of (\ref{1bis}) given by Proposition~\ref{lwp} is globally defined.
In addition there exists $C>0$ such that for every $u_0\in \Sigma_{\sigma}^{i}$,
\begin{equation}\label{growth}
\|u(t)\|_{H^{\sigma_1}(\Theta)}\leq C(i+\log(1+|t|))^{\frac{1}{2}}\,.
\end{equation}
Moreover, if $(u_{0,k})_{k\in\N}$, $u_{0,k}\in \Sigma^{i}_{\sigma,N_k}$,
$N_k\rightarrow\infty$ converges to $u_0$ as $k\rightarrow\infty$ in $H^{\sigma}_{rad}(\Theta)$ then  
\begin{equation}\label{limit}
\lim_{k\rightarrow\infty}\|u(t)-\Phi_{N_k}(t)(u_{0,k})\|_{H^{\sigma_1}(\Theta)}=0\,.
\end{equation}
\end{proposition}
\begin{proof}
Let $u_0\in \Sigma_{\sigma}^{i}$ and $(u_{0,k})$ $u_{0,k}\in \Sigma^{i}_{\sigma,N_k}$,
$N_k\rightarrow\infty$
a sequence
tending to $u_0$ in $H^{\sigma}_{rad}(\Theta)$. Let us fix $T>0$. Our aim so to extend the solution 
of (\ref{1bis}) given by Proposition~\ref{lwp} to the interval $[-T,T]$.
Using Proposition~\ref{longtime}, we have that there exists a constant $C$ such that for every $k\in\N$,
every $t\in\R$,
\begin{equation}\label{ant}
\|\Phi_{N_k}(t)(u_{0,k})\|_{H^{\sigma}(\Theta)}\leq C(i+\log(1+|t|))^{\frac{1}{2}}\,.
\end{equation}
Therefore, if we set $u_{N_k}(t)\equiv  \Phi_{N_k}(t)(u_{0,k})$ and 
$\Lambda\equiv C(i+\log(1+T))^{\frac{1}{2}}$, we have the bound
\begin{equation}\label{david1}
\|u_{N_k}(t)\|_{H^{\sigma}}\leq \Lambda,\quad \forall\,|t|\leq T,\quad \forall\, k\in\N.
\end{equation}
In particular $\|u_0\|_{H^{\sigma}}\leq\Lambda$ (apply (\ref{david1}) with $t=0$ and let $k\rightarrow\infty$).
Let $\tau>0$ be the local existence time for (\ref{1bis}), provided by
Proposition~\ref{lwp} for $\sigma_1$, $\sigma$ and $A=\Lambda+1$. 
Recall that we can assume $\tau=c(2+\Lambda)^{-\beta}$
for some $c>0$, $\beta>0$ depending only on $\sigma$ and $\sigma_1$. 
We can of course assume that $T>\tau$.
Denote by $u(t)$ the solution of (\ref{1bis}) with data $u_0$ on the time interval $[-\tau,\tau]$. Then
$v_{N_k}\equiv u-u_{N_k}$ solves the equation
\begin{equation}\label{eqnv}
(i\partial_{t}+\Delta) v_{N_k} =F(u)-S_{N_k}(F(u_{N_k})), \quad v_{N_k}|_{t=0}=u_0-u_{0,k} \, .
\end{equation}
Next, we write
$$
F(u)-S_{N_k}(F(u_{N_k}))=S_{N_k}\big(F(u)-F(u_{N_k})\big)+(1-S_{N_k})F(u).
$$
Therefore
\begin{multline*}
v_{N_k}(t)=e^{it\Delta}(u_0-u_{0,k})
\\
-i\int_{0}^{t}e^{i(t-\tau)\Delta}S_{N_k}\big(F(u(\tau))-F(u_{N_k}(\tau))\big)d\tau
-i\int_{0}^{t}e^{i(t-\tau)\Delta}(1-S_{N_k})F(u(\tau))d\tau\,.
\end{multline*}
Let us observe that for $\sigma_1<\sigma$ the map $1-S_{N}$ sends $H^{\sigma}_{rad}(\Theta)$ to
$H^{\sigma_{1}}_{rad}(\Theta)$ with norm $\leq CN^{\sigma_{1}-\sigma}$.
Similarly, for $I\subset\R$ an interval, the map
$1-S_{N}$ sends $X^{\sigma,b}_{rad}(I\times\Theta)$ to
$X^{\sigma_{1},b}_{rad}(I\times\Theta)$ with norm $\leq CN^{\sigma_{1}-\sigma}$.
Moreover $S_{N}$ acts as a bounded operator (with norm $\leq 1$) on the Bourgain spaces
$X^{\sigma,b}_{rad}$. Therefore, using Proposition~\ref{duh}, we obtain that there exist $C>0$,
$b>1/2$ and $\theta>0$ (depending only on $\sigma$, $\sigma_1$) such that one has the bound
\begin{eqnarray*}
\|v_{N_k}\|_{X^{\sigma_1,b}_{rad}([-\tau,\tau]\times\Theta)}
& \leq &
C\Big(\|u_0-u_{0,k}\|_{H^{\sigma_1}(\Theta)}
\\
&  & 
+ \tau^{\theta}\|v_{N_k}\|_{X^{\sigma_1,b}_{rad}([-\tau,\tau]\times\Theta)}
\big(1+\|u\|_{X^{\sigma_1,b}_{rad}([-\tau,\tau]\times\Theta)}^{\alpha_1}+
\|u_{N_k}\|_{X^{\sigma_1,b}_{rad}([-\tau,\tau]\times\Theta)}^{\alpha_1}
\big)
\\
& & + \tau^{\theta}N_k^{\sigma_1-\sigma}\|u\|_{X^{\sigma,b}_{rad}([-\tau,\tau]\times\Theta)}
\big(
1+\|u\|_{X^{\sigma_1,b}_{rad}([-\tau,\tau]\times\Theta)}^{\alpha_1}
\big)\Big),
\end{eqnarray*}
where $\alpha_1\equiv\max(2,\alpha)$. 
A use of Proposition~\ref{lwp} and Proposition~\ref{lwpbis} yields
\begin{eqnarray*}
\|v_{N_k}\|_{X^{\sigma_1,b}_{rad}([-\tau,\tau]\times\Theta)}
& \leq &
C\|u_0-u_{0,k}\|_{H^{\sigma_1}(\Theta)}
\\
& &
+
C\tau^{\theta}\|v_{N_k}\|_{X^{\sigma_1,b}_{rad}([-\tau,\tau]\times\Theta)}
(1+C\|u_{0}\|_{H^{\sigma_1}(\Theta)}^{\alpha_1}+C\|u_{0,k}\|_{H^{\sigma_1}(\Theta)}^{\alpha_1})
\\
& & +
C\tau^{\theta}N_k^{\sigma_1-\sigma}\|u_0\|_{H^{\sigma}(\Theta)}
(1+C\|u_{0}\|_{H^{\sigma_1}(\Theta)}^{\alpha_1})
\\
& \leq &
C\|u_0-u_{0,k}\|_{H^{\sigma_1}(\Theta)}
+
C\tau^{\theta}(1+\Lambda)^{\alpha_1}\|v_{N_k}\|_{X^{\sigma_1,b}_{rad}([-\tau,\tau]\times\Theta)}
\\
& &
+
C\tau^{\theta}(1+\Lambda)^{\alpha_1}N_k^{\sigma_1-\sigma}\|u_0\|_{H^{\sigma}(\Theta)}\,.
\end{eqnarray*}
Recall that $\tau=c(2+\Lambda)^{-\beta}$, where $c>0$ and $\beta>0$ are depending only on 
$\sigma$ and $\sigma_1$. In the last estimate the constants $C$ and $\theta$ also depend only on 
$\sigma_1$ and $\sigma$. Therefore, if we assume that $\beta>\alpha_1/\theta$ then the restriction on $\beta$
remains to depend only on $\sigma_1$ and $\sigma$. Similarly, if we assume that $c$ is so small that
$$
C\tau^{\theta}(1+\Lambda)^{\alpha_1}\leq
Cc^{\theta}(2+\Lambda)^{-\beta\theta}(1+\Lambda)^{\alpha_1}
\leq Cc^{\theta}
<1/2
$$
then the smallness restriction on $c$ remains to depend only on $\sigma_1$ and $\sigma$. 
Therefore, we have that after possibly slightly modifying the values of $c$ and $\beta$ 
(keeping $c$ and $\beta$ only depending on $\sigma$ and $\sigma_1$ and independent of $N_k$)
in the definition of $\tau$ that
\begin{eqnarray*}
\|v_{N_k}\|_{X^{\sigma_{1},b}_{rad}([-\tau,\tau]\times\Theta)}& \leq & C\|u_0-u_{0,k}\|_{H^{\sigma_1}(\Theta)}+
\frac{1}{2}N_k^{\sigma_1-\sigma}\|u_0\|_{H^{\sigma}(\Theta)}
\\
& = &
C\|v_{N_k}(0)\|_{H^{\sigma_1}(\Theta)}+
\frac{1}{2}N_k^{\sigma_1-\sigma}\|u(0)\|_{H^{\sigma}(\Theta)}\,.
\end{eqnarray*}
Since $b>1/2$, the last inequality implies
\begin{equation}\label{vvvv}
\|v_{N_k}(t)\|_{H^{\sigma_1}(\Theta)}\leq C\|v_{N_k}(0)\|_{H^{\sigma_1}(\Theta)}+
CN_k^{\sigma_1-\sigma}\|u(0)\|_{H^{\sigma}(\Theta)},\quad |t|\leq \tau= c(1+\Lambda)^{-\beta}\,,
\end{equation}
where the constants $c$, $C$ and $\beta$ depend only $\sigma_1$ and $\sigma$.
Therefore, using that
$$
\lim_{k\rightarrow\infty}\|v_{N_k}(0)\|_{H^{\sigma_1}(\Theta)}=0,
$$
we obtain that
\begin{equation}\label{rub1}
\lim_{k\rightarrow\infty}\|v_{N_k}(t)\|_{H^{\sigma_1}(\Theta)}=0,\quad |t|\leq \tau\,.
\end{equation}
Thus by taking  $N_k$ large enough in (\ref{vvvv}) one has via a use of
the triangle inequality,
\begin{equation}\label{david2}
\|u(t)\|_{H^{\sigma_1}(\Theta)}\leq\|u_{N_k}(t)\|_{H^{\sigma_1}(\Theta)}
+\|v_{N_k}(t)\|_{H^{\sigma_1}(\Theta)}\leq\Lambda+1,\quad |t|\leq\tau.
\end{equation}
Let us define the function $g_k(t)$ by
$$
g_k(t)\equiv \|v_{N_k}(t)\|_{H^{\sigma_1}(\Theta)}+N_k^{\sigma_1-\sigma}\|u(t)\|_{H^{\sigma}(\Theta)}\,.
$$
The function $g_k(t)$ is a priori defined only on $[-\tau,\tau]$. Our goal is to extend it on $[-T,T]$.
Using (\ref{vvvv}) and the bound 
$$
\|u(t)\|_{H^{\sigma}(\Theta)}\leq C\|u(0)\|_{H^{\sigma}(\Theta)}, \quad |t|\leq \tau,
$$ 
provided from Proposition~\ref{lwp}, we obtain that there exists a constant $C(\sigma,\sigma_1)$ depending 
only on $\sigma_1$ and $\sigma$ such that
$$
g_k(t)\leq C(\sigma,\sigma_1)g_k(0),\quad \forall\, t\in [-\tau,\tau]\,.
$$
We now repeat the argument for obtaining (\ref{vvvv}) on 
$[\tau,2\tau]$ and thanks to the bounds (\ref{david1}) and (\ref{david2}), we obtain
that $v_{N_k}(t)$ and $u$ exist on $[\tau,2\tau]$ and one has the bound
\begin{equation*}
\|v_{N_k}(t)\|_{H^{\sigma_1}(\Theta)}\leq C\|v_{N_k}(\tau)\|_{H^{\sigma_1}(\Theta)}+
CN_k^{\sigma_1-\sigma}\|u(\tau)\|_{H^{\sigma}(\Theta)},\quad t\in[\tau,2\tau]\,.
\end{equation*}
Therefore, thanks to (\ref{rub1}) (with $t=\tau$)
$$
\lim_{k\rightarrow\infty}\|v_{N_k}(t)\|_{H^{\sigma_1}(\Theta)}=0,\quad \tau\leq t\leq 2\tau\,.
$$
By  taking $N_k\gg 1$, we get via a use of the triangle inequality
\begin{equation*}
\|u(t)\|_{H^{\sigma_{1}}(\Theta)}\leq
\|u_{N_k}(t)\|_{H^{\sigma_1}(\Theta)}
+
\|v_{N_k}(t)\|_{H^{\sigma_1}(\Theta)}
\leq
\Lambda+1,\quad \tau \leq t\leq 2\tau.
\end{equation*}
Using (\ref{vvvv}) and the bound 
$$
\|u(t)\|_{H^{\sigma}(\Theta)}\leq C\|u(\tau)\|_{H^{\sigma}(\Theta)},\quad  \tau\leq t\leq 2\tau,
$$
provided from Proposition~\ref{lwp}, we obtain that 
$$
g_k(t)\leq C(\sigma,\sigma_1)g_k(\tau),\quad \forall\, t\in [\tau,2\tau].
$$
Then, we can continue by covering the interval $[-T,T]$ with intervals of size $\tau$, which yields
the existence of $u(t)$ on $[-T,T]$
(the point is that at each step the $H^{\sigma}$ norm of $u$ remains bounded by $\Lambda+1$
and the limit as $k\rightarrow\infty$ of the $H^{\sigma}$ norm of $v_{N_k}$ is zero).
Since $T>0$ was chosen arbitrary, we obtain that for every
$u_0\in\Sigma_{\sigma}^{i}$ the local solution of (\ref{1bis}) is globally defined.
Moreover 
$$
\|u(t)\|_{H^{\sigma_{1}}(\Theta)}\leq\Lambda+1,\quad |t|\leq T
$$
which by recalling the definition of $\Lambda$ implies the bound (\ref{growth}).
In addition, by iterating the bounds on $g_k$ we get at each step, 
we obtain the existence of a constant $C$ depending 
only on $\sigma$ and $\sigma_1$ such that
$$
g_k(t)\leq e^{C(1+|t|)}g_k(0)
$$
which implies that there exists a constant $C$ depending only on $\sigma_1$ and $\sigma$ such that 
$v_{N_k}$ enjoys the bound 
\begin{equation*}
\|v_{N_k}(t)\|_{H^{\sigma_{1}}(\Theta)}\leq 
C^{1+T}\Big(N_k^{\sigma_1-\sigma}\|u_0\|_{H^{\sigma}(\Theta)}
+\|u_0-u_{0,k}\|_{H^{\sigma_1}(\Theta)}\Big),\quad |t|\leq T.
\end{equation*}
Therefore for every $\varepsilon>0$ there exists $N^{\star}$ such that for $N_k\geq N^\star$ 
one has the inequality 
\begin{equation*}
\sup_{|t|\leq T}\|u(t)-\Phi_{N_k}(t)(u_{0,k})\|_{H^{\sigma_1}(\Theta)}<\varepsilon\,.
\end{equation*}
Hence we have (\ref{limit}).
This completes the proof of Proposition~\ref{global_existence}.
\end{proof}
By the Proposition~\ref{global_existence}, we can define a flow $\Phi$ acting on
$\Sigma_{\sigma}$, $\sigma\in [s,1/2[$ and defining the global dynamics of (\ref{1bis}) for 
$u_0\in \Sigma_{\sigma}$.
Let us now turn to the construction of a set invariant under $\Phi$.
Let $ l=(l_{j})_{j\in\N}$ be a increasing sequence of real numbers
such that $l_0=s$, $l_{j}<1/2$ and 
$
\lim_{j\rightarrow\infty}l_{j}=1/2.
$
Then, we set
\begin{equation*}
\Sigma=\bigcap_{\sigma\in l}\Sigma_{\sigma}\,.
\end{equation*}
The set $\Sigma$ is of full $\rho$ measure. It is the one involved in the statement of Theorem~\ref{thm3}.
Using the invariance property of $\Sigma^{i}_{N,\sigma}$, we now obtain that the set $\Sigma$  
is invariant under $\Phi$.
\begin{proposition}\label{compare}
For every $t\in\R$, $\Phi(t)(\Sigma)=\Sigma$.
In addition for every $\sigma\in l$, 
$\Phi(t)$ is continuous with respect to the induced by $H^\sigma_{rad}(\Theta)$ on $\Sigma$ topology.
In particular, the map
$\Phi(t):\Sigma\rightarrow\Sigma$ is a measurable map with respect to $\rho$.
\end{proposition}
\begin{proof}
Since the flow is time reversible, it suffices to show that
\begin{equation}\label{inclu}
\Phi(t)(\Sigma)\subset\Sigma,\quad \forall t\in\R\,.
\end{equation}
Indeed, if we suppose that (\ref{inclu}) holds true then for $u\in\Sigma$ and $t\in\R$, 
we have that thanks to (\ref{inclu}) $u_0\equiv\Phi(-t)u\in\Sigma$ (recall that $\Phi$ is well-defined 
on $\Sigma$ by Proposition~\ref{global_existence}) and thus
$u=\Phi(t)u_{0}$, i.e. $\Sigma\subset\Phi(t)(\Sigma)$.
Hence $\Phi(t)(\Sigma)=\Sigma$ is a consequence of (\ref{inclu}). 

Let us now prove (\ref{inclu}).
Fix $u_0\in\Sigma$ and $t\in\R$. It suffices to show that for every $\sigma_1\in l$, 
we have 
$$
\Phi(t)(u_0)\in \Sigma_{\sigma_1}\,.
$$
Let us take $\sigma\in ]\sigma_1,1/2[$, $\sigma\in l$. Since $u_0\in\Sigma$, we have that
$u_0\in \Sigma_{\sigma}$. 
Therefore there exists $i$ such that 
$
u_0\in \Sigma^{i}_{\sigma}.
$
Let $u_{0,k}\in\Sigma^{i}_{N_k,\sigma}$, $N_k\rightarrow\infty$ be a sequence which 
tends to $u_0$ in $H^{\sigma}(\Theta)$.
Thanks to Proposition~\ref{longtime} there exists $i_1$ such that
$$
\Phi_{N_k}(t)(u_{0,k})\in\Sigma^{i+i_1}_{N_k,\sigma_1}\,,\quad \forall\, k\in\N.
$$
Therefore using (\ref{limit}) of Proposition~\ref{global_existence}, we obtain that
$$
\Phi(t)(u_0)\in \Sigma^{i+i_1}_{\sigma_1}.
$$
Thus $\Phi(t)(u_0)\in\Sigma_{\sigma_1}$ which proves (\ref{inclu}).

Let us finally prove the continuity of $\Phi(t)$ on $\Sigma$ with respect to the $H^\sigma_{rad}(\Theta)$ 
topology.
Let $u\in\Sigma$ and $u_n\in\Sigma$ be a sequence such that $u_n\rightarrow u$ in $H^\sigma_{rad}(\Theta)$. 
We need to prove that for every $t\in\R$,
$
\Phi(t)(u_n)\rightarrow \Phi(t)(u)
$
in $H^\sigma_{rad}(\Theta)$. Let us fix $t\in\R$. 
Since $u\in\Sigma$ (and thus in all $\Sigma_{\sigma}$, $\sigma\in l$), 
using Proposition~\ref{global_existence}, we obtain that there exists $C>0$ such that
\begin{equation}\label{kyoto}
\sup_{|\tau|\leq |t|}\|\Phi(\tau)(u)\|_{H^\sigma(\Theta)}\leq C(\log(2+|t|))^{\frac{1}{2}}\equiv \Lambda.
\end{equation}
Let us denote by $\tau_0$ the local existence time in Proposition~\ref{lwp}, associated to $\sigma$ and
$A=2\Lambda $. Then, by the continuity of the flow given by Proposition~\ref{lwp}, we have
$\Phi(\tau_0)(u_n)\rightarrow \Phi(\tau_0)(u)$ in $H^\sigma_{rad}(\Theta)$. Next, we cover
the interval $[0,t]$ by intervals of size $\tau_0$ and we apply the continuity
of the flow established in  Proposition~\ref{lwp} at each step. 
The applicability of Proposition~\ref{lwp} is possible thanks to the bound (\ref{kyoto}).
Therefore, we obtain that $\Phi(t)(u_n)\rightarrow \Phi(t)(u)$ in $H^\sigma_{rad}(\Theta)$.
This completes the proof of Proposition~\ref{compare}.
\end{proof}
\section{Proof of the measure invariance}
Fix $\sigma\in]s,1/2[$, $\sigma\in l$.
Thanks to the invariance by $\Phi$ of the set $\Sigma$,
using the regularity of the measure $\mu$ (which is a finite Borel measure) and Remark~\ref{rem}, 
we deduce that it suffices to prove the measure invariance for subsets $K$ of $\Sigma$
which are compacts of $H^s_{rad}(\Theta)$ and which are bounded in $H^\sigma_{rad}(\Theta)$.
Let us fix $t\in \R$ and a compact $K$ of $H^s_{rad}(\Theta)$ which is a bounded set in 
$H^\sigma_{rad}(\Theta)$. Our aim is to show that 
$
\rho(\Phi(t)(K))=\rho(K).
$
By the time reversibility of the flow, we may suppose that $t>0$.
Since $K$ is bounded in $H^\sigma_{rad}(\Theta)$ and a compact in $H^s_{rad}(\Theta)$, 
using the continuity property displayed
by Proposition~\ref{compare} and Proposition~\ref{lwp}, we infer that there exists $R>0$ such that
\begin{equation}\label{krum}
\{\Phi(\tau)(K),\,\, 0\leq \tau\leq t\}\subset
\{u\in H^\sigma_{rad}(\Theta)\,:\, \|u\|_{H^\sigma(\Theta)}\leq R\}\equiv B_{R}\,.
\end{equation}
Indeed, the left hand-side of (\ref{krum}) is included in a sufficiently large
$H^s_{rad}(\Theta)$ ball thanks to the continuity property of the flow on
$H^s_{rad}(\Theta)$ shown in Proposition~\ref{compare} and the compactness of
$K$. Then, by iterating the propagation of regularity statement of
Proposition~\ref{lwp},
applied with $A$ such that the $H^s_{rad}(\Theta)$ ball centered at the origin of
radius $A$ contains the left hand-side of (\ref{krum}), we arrive at
(\ref{krum}) (observe that we only have the poor bound $R\sim e^{Ct}$).
\\

Let $c$ and $\beta$ (depending only on $s$ and $\sigma$)
be fixed by an application of Proposition~\ref{lwp} with $s=\sigma_1$ and $\sigma=\sigma$.
Next, we set
$$
\tau_0\equiv c_0(1+R)^{-\beta_0},
$$
where $0<c_0\leq c$, $\beta_0\geq \beta$, depending only on $s$ and $\sigma$,
are to be fixed in the next lemma which allows to compare $\Phi$ and $\Phi_N$ for data in $B_{R}$.

\begin{lemme}\label{nedelia}
There exist $c_0$ and $\beta_0$ depending only on $s$ and $\sigma$ such that
for every $\varepsilon>0$ there exists $N_0\geq 1$ such that for every $N\geq N_0$, every $u_0\in B_{R}$,
every $\tau\in [0,\tau_0]$,
$$
\|\Phi(\tau)(u_0)-\Phi_{N}(\tau)(S_{N}(u_0))\|_{H^s(\Theta)}<\varepsilon\,.
$$
\end{lemme}
\begin{proof}
For $u_0\in B_{R}$, we denote by $u$ the solution of (\ref{1bis}) with data $u_0$
and by $u_{N}$ the solution of (\ref{Nbis}) with data $S_{N}(u_0)$, defined on $[0,\tau_0]$. 
Next, we set $v_{N}\equiv u-u_N$. Then $v_{N}$ solves
\begin{equation}\label{eqnvpak}
(i\partial_t+\Delta)v_N= F(u)-S_{N}(F(u_{N})), \quad v_{N}(0)=(1-S_{N})u_0\,.
\end{equation}
By writing
$$
F(u)-S_{N}(F(u_{N}))=S_{N}\big(F(u)-F(u_{N})\big)+(1-S_{N})F(u)
$$
and using Proposition~\ref{duh}, we obtain that there exists $b>1/2$ and $\theta>0$ 
depending only on $s$ and $\sigma$ such that one has
\begin{eqnarray*}
\|v_{N}\|_{X^{s,b}_{rad}([0,\tau_0]\times\Theta)}
& \leq &
CN^{s-\sigma}\|u_0\|_{H^{\sigma}(\Theta)}
\\
&  & + C\tau_0^{\theta}\|v_N\|_{X^{s,b}_{rad}([0,\tau_0]\times\Theta)}
\big(1+\|u\|_{X^{s,b}_{rad}([0,\tau_0]\times\Theta)}^{\max(2,\alpha)}
+\|u_{N}\|_{X^{s,b}_{rad}([0,\tau_0]\times\Theta)}^{\max(2,\alpha)}\big)
\\
& & 
+ C\tau_0^{\theta}N^{s-\sigma}\|u\|_{X^{\sigma,b}_{rad}([0,\tau_0]\times\Theta)}
\big(1+\|u\|_{X^{s,b}_{rad}([0,\tau_0]\times\Theta)}^{\max(2,\alpha)}\big).
\end{eqnarray*}
Using Proposition~\ref{lwp} and Proposition~\ref{lwpbis}, we get
\begin{eqnarray*}
\|v_N\|_{X^{s,b}_{rad}([0,\tau_0]\times\Theta)}
& \leq &
CN^{s-\sigma}\|u_0\|_{H^{\sigma}(\Theta)}
\\
& &
+C\tau_0^{\theta}\|v_N\|_{X^{s,b}_{rad}([0,\tau_0]\times\Theta)}(1+C\|u_{0}\|_{H^{s}(\Theta)}^{\max(2,\alpha)})
\\
& & 
+C\tau_0^{\theta}N^{s-\sigma}\|u_0\|_{H^{\sigma}(\Theta)}(1+C\|u_{0}\|_{H^{s}(\Theta)}^{\max(2,\alpha)})\,.
\end{eqnarray*}
Coming back to the definition of $\tau_0$
we can choose $c_0$ small enough and $\beta_0$ large enough,
but keeping their dependence only on $s$ and $\sigma$, to infer that
$$
\|v_N\|_{X^{s,b}_{rad}([0,\tau_0]\times\Theta)}\leq CN^{s-\sigma}\|u_0\|_{H^{\sigma}(\Theta)}.
$$
Since $b>1/2$, by the Sobolev embedding, the space $X^{s,b}_{rad}([0,\tau_0]\times\Theta)$ is continuously 
embedded in $L^{\infty}([0,\tau_0];H^s_{rad}(\Theta))$ and thus there exists $C$ depending only on $s$, 
$\sigma$ such that
\begin{equation*}
\|v_N(t)\|_{H^{s}(\Theta)}\leq CRN^{s-\sigma},\quad t\in [0,\tau_0].
\end{equation*}
This completes the proof of Lemma~\ref{nedelia}.
\end{proof}
It suffices to prove that
\begin{equation}\label{reduction}
\rho(\Phi(\tau)(K))=\rho(K),\quad \tau\in [0,\tau_0].
\end{equation}
Indeed, it suffices to cover $[0,t]$ by intervals of size $\tau_0$ and apply
(\ref{reduction}) at each step. Such an iteration is possible since 
by the continuity property of $\Phi(t)$ at each
step the image remains a compact of $H^s_{rad}(\Theta)$ included in the ball $B_{R}$. 
Let us now prove (\ref{reduction}).
Let $B_{\varepsilon}$ be the open ball in $H^{s}_{rad}(\Theta)$ centered at
the origin and of radius $\varepsilon$.
We have that $\Phi(\tau)(K)$ is a closed set of $H^{s}_{rad}(\Theta)$
contained in $\Sigma$. Therefore, by Theorem~\ref{thm2}, we can write
$$
\rho\Big(\Phi(\tau)(K)+\overline{B_{2\varepsilon}}\Big) \geq \limsup_{N\rightarrow \infty}
\rho_{N}\Big(\big(\Phi(\tau)(K)+\overline{B_{2\varepsilon}}\big)\cap E_{N}\Big)\,,
$$
where $\overline{B_{2\varepsilon}}$ is the closed ball in
$H^{s}_{rad}(\Theta)$, centered at the origin and of radius $2\varepsilon$.
Using Lemma~\ref{nedelia}, we obtain that for every $\varepsilon>0$, if we take $N$ large enough, we
have
$$
\big(\Phi_{N}(\tau)(S_{N}(K))+B_{\varepsilon}\big)\cap E_{N}
\subset
\big(\Phi(\tau)(K)+\overline{B_{2\varepsilon}}\big)\cap E_{N}
$$
and therefore
$$
\limsup_{N\rightarrow \infty}
\rho_{N}\Big(\big(\Phi(\tau)(K)+\overline{B_{2\varepsilon}}\big)\cap E_{N}\Big)
\geq 
\limsup_{N\rightarrow \infty}
\rho_{N}\Big(
\big(\Phi_{N}(\tau)(S_{N}(K))+B_{\varepsilon}\big)\cap E_{N}
\Big).
$$
Next, using the Lipschitz continuity of the flow $\Phi_N$ (see Proposition~\ref{lwpbis}), we obtain that
there exists $c\in]0,1[$, independent of $\varepsilon$ such that for $N$ large enough, we have
$$
\Phi_{N}(\tau)\big((K+B_{c\varepsilon})\cap E_{N}\big)
\subset
\big(\Phi_{N}(\tau)(S_{N}(K))+B_{\varepsilon}\big)\cap E_{N},
$$
where $B_{c\varepsilon}$ is the open ball in $H^{s}_{rad}(\Theta)$ centered at the origin and 
of radius $c\varepsilon$. Therefore
$$
\limsup_{N\rightarrow \infty}
\rho_{N}\Big(\big(\Phi_{N}(\tau)(S_{N}(K))+B_{\varepsilon}\big)\cap E_{N}\Big) 
\geq 
\limsup_{N\rightarrow \infty}\rho_{N}\Big(\Phi_{N}(\tau)\big(
(K+B_{c\varepsilon})\cap E_{N}\big)\Big)\, .
$$
Further, using the invariance of $\rho_N$ under $\Phi_N$, we obtain that
$$
\rho_{N}\Big(\Phi_{N}(\tau)\big((K+B_{c\varepsilon})\cap E_{N}\big)\Big)
=
\rho_{N}\Big((K+B_{c\varepsilon})\cap E_{N}\Big)
$$
and thus
$$
\limsup_{N\rightarrow \infty}\rho_{N}\Big(\Phi_{N}(\tau)\big((K+B_{c\varepsilon})\cap E_{N}\big)\Big)
\geq
\liminf_{N\rightarrow \infty}\rho_{N}\Big((K+B_{c\varepsilon})\cap E_{N}\Big).
$$
Finally, invoking once again  Theorem~\ref{thm2}, we can write
$$
\liminf_{N\rightarrow \infty}\rho_{N}\Big((K+B_{c\varepsilon})\cap E_{N}\Big)
\geq
\rho(K+B_{c\varepsilon})\geq \rho(K).
$$
Therefore, we have the inequality
$$
\rho\Big(\Phi(\tau)(K)+\overline{B_{2\varepsilon}}\Big) \geq \rho(K).
$$
By letting $\varepsilon\rightarrow 0$, thanks to the dominated convergence, 
we obtain that 
$$
\rho(\Phi(\tau)(K))\geq \rho(K).
$$
By the time reversibility of the flow we get $\rho(\Phi(\tau)(K))= \rho(K)$
and thus the measure invariance. 
\\

This completes the proof of Theorem~\ref{thm3}.\qed
\section{Concerning the three dimensional case}
\subsection{General discussion}
The extension of the result to the 3d case is an interesting problem.
In this case one can still prove the measure existence. The Cauchy problem issue 
is much more challenging. Despite the fact that the Cauchy problem for $H^{\sigma}$, $\sigma<1/2$ data
is ill-posed, in the sense of failure of continuity of the flow map (see the work of 
Christ-Colliander-Tao \cite{CCT}, or the appendix of \cite{BGT}), we may hope that estimates on Wiener chaos
can help us to resolve globally (with uniqueness) the Cauchy problem a.s. on a suitable statistical 
ensemble $\Sigma$ (which is included in the intersection of $H^{\sigma}$, $\sigma<1/2$ and misses $H^{1/2}$).
This would be an example showing the possibility to get strong solutions of a dispersive equation,
a.s. with respect to a measure, beyond the Hadamard well-posedness threshold.
In this section, we prove an estimate which shows that one has a control on the second Picard iteration, 
in all $H^{\sigma}$, $\sigma<1/2$, a.s. with respect to the measure. 
We will consider zonal solutions of the cubic defocusing NLS on the sphere $S^3$. The analysis of 
this model has a lot of similarities with the analysis on the ball of $\R^3$
(which is the three dimensional analogue of (\ref{1})).
There are however some simplifications because of the absence of boundary on $S^3$ and a nice formula 
for the products of zonal eigenfunctions.
In this section, we will benefit from some computations of the unpublished manuscript \cite{BGT-zonal}.
\subsection{Zonal functions on $S^3$ }
Let $S^3$ be the unit sphere in $\R^{4}$. If we consider functions on $S^3$ depending only on the geodesic 
distance to the north pole, we obtain the zonal functions on $S^3$.
The zonal functions can be expressed in terms of zonal spherical harmonics which in their turn can be 
expressed in terms the classical Jacobi polynomials. 
Let $\theta\in[0,\pi]$ be a local parameter measuring the geodesic distance to the north pole of  $S^3$.
Define the space $L^{2}_{rad}(S^3)$ to be equipped with the following norm
$$
\|f\|_{L^{2}_{rad}(S^3)}=\Big(\int_{0}^{\pi}|f(\theta)|^{2}(\sin\theta)^{2} d\theta\Big)^{\frac{1}{2}},
$$
where $f$ is a zonal function on $S^3$ and $(\sin\theta)^{2} d\theta$ is the surface measure on $S^3$. 
One can define similarly other functional spaces of zonal functions, 
for example $L^{p}_{rad}(S^3)$, $H^{s}_{rad}(S^3)$ etc. The Laplace-Beltrami 
operator on $L^2(S^3)$ can be restricted to
$L^{2}_{rad}(S^3)$ and in the coordinate $\theta$ it reads
$$
\frac{\partial^{2}}{\partial \theta^2}+\frac{2}{{\rm tg }\, \theta}\frac{\partial}{\partial\theta}
$$
since using the parametrization of $S^3$ in terms of $\theta$ and $S^{2}$, one can write,
$$
\Delta_{S^3}=\frac{\partial^{2}}{\partial \theta^2}+\frac{2}{{\rm tg }\, \theta}
\frac{\partial}{\partial\theta}+\frac{1}{\sin^{2}\theta}\Delta_{S^{2}}\,.
$$
It follows from the Sturm-Liouville theory (see also e.g. \cite{SW}) that an orthonormal basis of 
$L^{2}_{rad}(S^3)$ can be build by the functions 
$$
P_{n}(\theta)=\sqrt{\frac{2}{\pi}}\,\frac{\sin n\theta}{\sin\theta},\quad\theta\in[0,\pi],\quad n\geq 1,
$$
where $\theta$ connotes the geodesic distance to the north pole of $S^3$. 
The functions $P_n$ are eigenfunctions of $-\Delta_{S^3}$ 
with corresponding eigenvalue $\lambda_{n}=n^2-1$.
We next define the function $\gamma : \N^{4}\longrightarrow \R$ by
$$
\gamma(n,n_1,n_2,n_3)\equiv\int_{S^3}
P_{n}P_{n_1}P_{n_2}P_{n_3}.
$$
Then clearly
$$
P_{n_1}P_{n_2}P_{n_3}
=
\sum_{n=1}^{\infty}\gamma(n,n_1,n_2,n_3)P_{n}
$$
and thus the behaviour of $\gamma$ would be of importance when analysing cubic expressions on $S^3$.
In the next lemma we give a bound for $\gamma(n,n_1,n_2,n_3)$.
\begin{lemme}\label{l1} 
One has the bound
$
0\leq \gamma(n,n_1,n_2,n_3)\leq (2/\pi)\min(n,n_1,n_2,n_3).
$
\end{lemme}
\begin{proof}
Using the explicit formula for $P_n$ and some trigonometric considerations, we obtain the relation 
\begin{equation}\label{basic}
P_{k}P_{l}=
\sqrt{\frac{2}{\pi}}\,
\sum_{j=1}^{\min(k,l)}
P_{|k-l|+2j-1},\quad k\geq 1,\, l\geq 1.
\end{equation}
By symmetry we can suppose that $n_{1}\geq n_{2}\geq n_{3}\geq n$. Then due to (\ref{basic}) we obtain 
that $P_{n}P_{n_3}$ can be
expressed as a sum of $n$ terms while the sum corresponding to $P_{n_1}P_{n_2}$ contains $n_2$ terms. 
Since for $k\neq l$ one has
$\int_{S^3}P_{k}P_{l}=0$, we obtain that the contribution to $\gamma(n,n_1,n_2,n_3)$ of any of the term 
of the sum for $P_{n}P_{n_3}$
is not more than $2/\pi$ and therefore $\gamma(n,n_1,n_2,n_3)\leq (2/\pi)n$. This completes the proof of Lemma~\ref{l1}.
\end{proof}
We shall also make use of the following property of 
$\gamma(n,n_1,n_2,n_3)$.
\begin{lemme}\label{l2}
Let $n>n_{1}+n_{2}+n_{3}$. Then $\gamma(n,n_1,n_2,n_3)=0.$
\end{lemme}
\begin{proof}
One needs simply to observe that in the spectral decomposition of $P_{n_1}P_{n_2}P_{n_3}$ there 
are only spherical
harmonics of degree $\leq n_{1}+n_2+n_3$ and therefore $P_{n_1}P_{n_2}P_{n_3}$ is orthogonal to $P_n$. 
This completes the proof of Lemma~\ref{l2}.
\end{proof}
\begin{remarque}
Let us observe that $(\pi/2)\gamma(n,n_1,n_2,n_3)\in\Z$. This fact is however not of importance for the sequel. 
\end{remarque}
\subsection{The cubic defocusing NLS on $S^3$}
Consider the cubic defocusing nonlinear Schr\"odinger equation, posed on $S^3$,
\begin{equation}\label{5}
(i\partial_{t}+\Delta_{S^3})u-|u|^{2}u=0,
\end{equation}
where $u:\R\times S^3\longrightarrow \C$. By the variable change $u\rightarrow e^{it}u$, 
we can reduce (\ref{5}) to 
\begin{equation}\label{6}
(i\partial_{t}+\Delta_{S^3}-1)u-|u|^{2}u=0.
\end{equation}
We will perform our analysis to the equation (\ref{6}).
The Hamiltonian associated to (\ref{6}) is 
$$
H(u,\bar{u})=\int_{S^3}|\nabla u|^{2}+\int_{S^3}|u|^{2}+\frac{1}{2}\int_{S^3}|u|^{4},
$$
where $\nabla$ denotes the riemannian gradient on $S^3$.
We will study zonal solutions of (\ref{6}), i.e. solutions such that $u(t,\cdot)$ is a zonal function on $S^3$.
Let us fix $s<1/2$.
The free measure, denoted by $\mu$, associated to (\ref{6}) is the distribution of the $H^s_{rad}(S^3)$ 
random variable
\begin{equation*}
\varphi(\omega,\theta)=
\sum_{n= 1}^{\infty}\frac{g_n(\omega)}{n}P_{n}(\theta)\, ,
\end{equation*}
where $g_n(\omega)$ is a sequence of centered, normalised, independent identically distributed (i.i.d.)
complex Gaussian random variables, defined in a probability space $(\Omega,{\mathcal F},p)$.
Using Lemma~\ref{l1}, we obtain that
$$
\|P_{n}\|_{L^4(S^3)}\leq n^{\frac{1}{4}}\,.
$$
and therefore using Lemma~\ref{lem1}, as in the proof of Theorem~\ref{thm1},
we get
$$
\|\varphi(\omega,\theta)\|^{2}_{L^4(\Omega\times S^3)}\leq\sum_{n= 1}^{\infty}
\frac{C}{n^2}\|P_{n}\|_{L^4(S^3)}^{2}\leq C \sum_{n= 1}^{\infty}\frac{n^{\frac{1}{2}}}{n^2}<\infty\,.
$$
Hence the image measure on $H^s_{rad}(S^3)$ under the map
$$
\omega\longmapsto \sum_{n= 1}^{\infty}\frac{g_n(\omega)}{n}P_{n}(\theta)\, ,
$$
of 
$$
\exp\Big(-\frac{1}{2}\|\varphi(\omega,\cdot)\|_{L^4(S^3)}^{4}\Big)dp(\omega)
$$
is a nontrivial measure which could be expected to be invariant under a flow of (\ref{6}).
For that purpose one should define global dynamics of (\ref{6}) on a set of full $\mu$ measure, 
i.e. solutions of (\ref{6}) with data $\varphi(\omega,\theta)$ for 
typical $\omega$'s. Using for instance the Fernique integrability theorem one has that 
$\|\varphi(\omega,\cdot)\|_{H^{1/2}(S^3)}=\infty$ $\mu$ a.s. Thus one needs to establish a 
well-define (and stable in a suitable sense) dynamics for data of Sobolev regularity $<1/2$. 
There is a major problem if one tries to solve this problem for individual 
$\omega$'s since the result of \cite{CCT} (see also the appendix of \cite{BGT}) shows that (\ref{6}) is 
in fact ill-posed for data of Sobolev regularity $<1/2$ and 
the data giving the counterexample can be chosen to be a zonal function since the analysis uses only 
point concentrations.
Therefore, it is possible that solving (\ref{6}) with data $\varphi(\omega,\theta)$, for typical $\omega$'s, 
would require a probabilistic argument in the spirit 
of the definition of the stochastic integration. 
Below, we present an estimate which gives a control on the second Picard iteration with data 
$\varphi(\omega,\theta)$.
\\

Let us consider the integral equation (Duhamel form) corresponding to (\ref{6}) with data 
$\varphi(\omega,\theta)$
\begin{equation}\label{Duhamel}
u(t)=S(t)(\varphi(\omega,\cdot))-i\int_{0}^{t}S(t-\tau)(|u(\tau)|^{2}u(\tau))d\tau,
\end{equation}
where $S(t)=\exp(it(\Delta_{S^3}-1))$ is the unitary group generated by the free evolution. 
The operator $S(t)$ acts as an isometry on $H^{s}(S^3)$ which can be easily seen by expressing $S(t)$ 
in terms of the spectral decomposition. 
One can show (see \cite{BGT}) that for $s>1/2$, the Picard iteration applied in the context of 
(\ref{Duhamel}) converges, if we replace $\varphi(\omega,\cdot)$ 
in (\ref{Duhamel}) by data in $u_0\in H^{s}(S^3)$, in the Bourgain spaces $X^{s,b}([-T,T]\times S^3)$, 
where $b>1/2$ is close to $1/2$, $T\sim (1+\|u_0\|_{H^s(S^3)})^{-\beta}$ 
(for some $\beta>0$ depending on $b$ and $s$).
For the definition the Bourgain spaces $X^{s,b}([-T,T]\times S^3)$ associated
to $\Delta_{S^3}$, we refer to \cite{BGT} (see also (\ref{7}) below).
The modification for $\Delta_{S^3}-1$ is then direct. 
Let us set (the first Picard iteration)
$$
u_{1}(\omega,t,\theta)\equiv S(t)(\varphi(\omega,\cdot))=
\sum_{n=1}^{\infty}\frac{g_{n}(\omega)}{n}P_{n}(\theta)e^{-itn^2}  \,.
$$
The random variable $u_{1}$ represents the free evolution. Notice that again 
$$
\|u_1(\omega,t,\cdot)\|_{H^{1/2}(S^3)}=\infty,\quad {\rm a.s.}
$$   
but for every $\sigma<1/2$,
$$
\|u_1(\omega,t,\cdot)\|_{H^{\sigma}(S^3)}<\infty,\quad {\rm a.s.}
$$
Let us consider the second Picard iteration
$$
u_{2}(\omega,t,\theta)\equiv S(t)(\varphi(\omega,\cdot))-
i\int_{0}^{t}S(t-\tau)(|u_1(\omega,\tau)|^{2}u_1(\omega,\tau))d\tau\,.
$$
Set
$$
v_{2}(\omega,t,\theta)\equiv\int_{0}^{t}S(t-\tau)(|u_1(\omega,\tau)|^{2}u_1(\omega,\tau))d\tau\,.
$$
Thanks to the ``dispersive effect'', $v_{2}$ is again a.s. in all $H^{\sigma}(S^3)$ for $\sigma<1/2$.
\begin{proposition}\label{thm4}
Let us fix $\sigma<1/2$. Then for $b>1/2$ close to $1/2$ and every $T>0$,
$$
\|v_{2}(\omega,t,\theta)\|_{L^2(\Omega\,;\,X^{\sigma,b}([-T,T]\times S^3))}<\infty.
$$
In particular
$$
\|v_{2}(\omega,t,\theta)\|_{L^2(\Omega\,;\,L^{\infty}([-T,T]\,;\,H^{\sigma}(S^3)))}<\infty
$$
and thus $\|v_{2}(\omega,\cdot,\cdot)\|_{L^{\infty}([-T,T]\,;\,H^{\sigma}(S^3))}$ is a.s. finite
which implies that the second Picard iteration for (\ref{Duhamel}) is a.s. in $H^{\sigma}$.
\end{proposition}
\begin{remarque}
Using estimates on the third order Wiener chaos, we might show that higher moments and Orlitch norms 
with respect to $\omega$ are finite.
\end{remarque}
\begin{proof}[Proof of Proposition~\ref{thm4}]
Let $\psi\in C_{0}^{\infty}(\R;\R)$ be a bump function localizing in $[-T,T]$. 
Let $\psi_1\in C_{0}^{\infty}(\R;\R)$ be a bump function which equals one on the support of $\psi$.
Set 
$$
w_{1}(t)\equiv \psi_1(t)u_{1}(t).
$$
Then using \cite{Gi}, for $b>1/2$ (close to $1/2$),
\begin{eqnarray*}
\|v_{2}(\omega,\cdot)\|_{X^{\sigma,b}([-T,T]\times S^3)}
& \leq &
\|\psi\, v_{2}(\omega,\cdot)\|_{X^{\sigma,b}(\R\times S^3)}
\\
& \leq & 
C\||w_{1}(\omega,\cdot)|^{2}w_{1}(\omega,\cdot)\|_{X^{\sigma,b-1}(\R\times S^3)}.
\end{eqnarray*}
Set
$$
w(\omega,t,\theta)\equiv |w_{1}(\omega,t,\theta)|^{2}w_{1}(\omega,t,\theta).
$$
We need to show that the $L^2(\Omega)$ of $\|w(\omega,\cdot)\|_{X^{\sigma,b-1}(\R\times S^3)}$ is finite. 
If
$$
w(\omega,t,\theta)=\sum_{n=1}^{\infty}c(\omega,n,t)P_{n}(\theta)
$$
then we have 
\begin{equation}\label{7}
\|w(\omega,\cdot)\|_{X^{\sigma,b-1}(\R\times S^3)}^{2}
=
\sum_{n=1}^{\infty}\int_{-\infty}^{\infty}
\langle \tau+n^{2}\rangle^{2(b-1)}n^{2\sigma}|\widehat{c}(\omega,n,\tau)|^{2}d\tau,
\end{equation}
where $\widehat{c}(\omega,n,\tau)$ denotes the Fourier transform with respect to $t$ of $c(\omega,n,t)$.
Let us next compute $c(\omega,n,t)$. This will of course make appeal to the function $\gamma$ introduced 
in the previous section.
We have that
$$
w(\omega,t,\theta)=\psi_{1}^{3}(t)
\sum_{(n_1,n_2,n_3)\in\N^3}
\frac{g_{n_1}(\omega) \overline{g_{n_2}(\omega)} g_{n_3}(\omega) }
{n_1 n_2 n_3}P_{n_1}(\theta)P_{n_2}(\theta)P_{n_3}(\theta)
e^{-it(n_1^2-n_2^2+n_3^2)}
$$
and therefore
$$
c(\omega,n,t)=\psi_{1}^{3}(t)\sum_{(n_1,n_2,n_3)\in\N^3}\gamma(n,n_1,n_2,n_3)
\frac{g_{n_1}(\omega) \overline{g_{n_2}(\omega)} g_{n_3}(\omega) }
{n_1 n_2 n_3}e^{-it(n_1^2-n_2^2+n_3^2)}\,.
$$
If we denote $\psi_2=\psi_1^3$ then
$$
\widehat{c}(\omega,n,\tau)=\sum_{(n_1,n_2,n_3)\in\N^3}\gamma(n,n_1,n_2,n_3)
\frac{g_{n_1}(\omega) \overline{g_{n_2}(\omega)} g_{n_3}(\omega) }
{n_1 n_2 n_3}\widehat{\psi_2}\big(\tau+n_1^2-n_2^2+n_3^2\big)\,.
$$
Let us observe that thanks to the independence of $(g_{n})_{n\in\N}$ we have that there are essentially 
two different situations when the expression
\begin{equation}\label{combin}
\int_{\Omega}
g_{n_1}(\omega) \overline{g_{n_2}(\omega)} g_{n_3}(\omega)
\overline{g_{m_1}(\omega)} g_{m_2}(\omega)\overline{g_{m_3}(\omega)}
dp(\omega)
\end{equation}
is different from zero. Namely 
\begin{itemize}
\item
$n_1=m_1$, $n_2=m_2$, $n_3=m_3$, 
\item
$n_1=n_2$, $n_3=m_1$, $m_2=m_3$.
\end{itemize}
Indeed, the complex gaussians $g_n$ satisfy
$$
\int_{\Omega}g_{n}(\omega)dp(\omega)=\int_{\Omega}g^2_{n}(\omega)dp(\omega)
=\int_{\Omega}g^3_{n}(\omega)dp(\omega)=\int_{\Omega}|g_{n}(\omega)|^2g_{n}(\omega)dp(\omega)=0
$$
and thus in order to have a nonzero contribution of (\ref{combin}) each gaussian without a bar 
in the integral (\ref{combin}) should be coupled with another gaussian
having a bar and the same index.

Therefore coming back to (\ref{7}), we get
$$
\int_{\Omega}\|w(\omega,\cdot)\|_{X^{\sigma,b-1}(\R\times S^3)}^{2}\,dp(\omega)\leq C(I_1+I_2),
$$
where
\begin{equation}\label{8}
I_1=\sum_{(n,n_1,n_2,n_3)\in\N^4}\int_{-\infty}^{\infty}\frac{n^{2\sigma}}{\langle \tau+n^{2}\rangle^{\beta}}
\frac{\gamma^{2}(n,n_1,n_2,n_3)}{(n_1 n_2 n_3)^{2}}|\widehat{\psi_2}|^{2}\big(\tau+n_1^2-n_2^2+n_3^2\big)d\tau,
\end{equation}
with $\beta\equiv 2(1-b)$ and
\begin{equation}\label{8bis}
I_2=
\sum_{(n,n_1,n_2,n_3)\in\N^4}\int_{-\infty}^{\infty}\frac{n^{2\sigma}}{\langle \tau+n^{2}\rangle^{\beta}}
\frac{\gamma(n,n_1,n_1,n_2)\gamma(n,n_2,n_3,n_3)}{(n_1 n_1 n_2)(n_2 n_3 n_3)}
|\widehat{\psi_2}|^{2}\big(\tau+n_2^2\big)d\tau\,.
\end{equation}
Notice that $\beta<1$ is close to $1$ when $b>1/2$ is close to $1/2$.
Thus our goal is to show the convergence of (\ref{8}) and (\ref{8bis}).
For that purpose we make appeal to the following lemma.
\begin{lemme}\label{ihp}
For every $\sigma\in]0,1/2[$ there exist $\beta<1$ and $C>0$ such that for every $\alpha\in\R$,
\begin{equation}\label{vlak}
\sum_{n=1}^{\infty}\frac{n^{2\sigma}}{(1+|n^2-\alpha|)^{\beta}}\leq
C(1+|\alpha|)^{\sigma}\,.
\end{equation}
\end{lemme}
\begin{proof}
Let $\beta<1$ be such that $2\beta-2\sigma>1$, i.e.
$
1/2+\sigma<\beta<1.
$
We prove (\ref{vlak}) for such values of $\beta$. The contribution of the region $\frac{1}{4}n^2\geq|\alpha|$ 
to the left hand-side of (\ref{vlak}) can be bounded by
$$
\sum_{n=1}^{\infty}\frac{n^{2\sigma}}{(1+\frac{3}{4}n^2)^{\beta}}\leq C_{\sigma}\leq 
C_{\sigma}(1+|\alpha|)^{\sigma}
$$
thanks to the assumption $2\beta-2\sigma>1$ and since for $\frac{1}{4}n^2\geq|\alpha|$ one has
$|n^2-\alpha|\geq \frac{3}{4}n^2$. We next estimate the contribution of the region 
$\frac{1}{4}n^2\leq|\alpha|$ (if it is not empty) by
$$
(4|\alpha|)^{\sigma}\sum_{n=1}^{\infty}\frac{1}{(1+|n^2-\alpha|)^{\beta}}
\leq C_{\sigma}|\alpha|^{\sigma}\,.
$$
This completes the proof of Lemma~\ref{ihp}.
\end{proof}
Let us now show the convergence of (\ref{8}).
Using the rapid decay of $|\widehat{\psi_2}|^{2}$, we can eliminate the $\tau$ integration and arrive at
\begin{equation}\label{9}
(\ref{8})\leq
C\sum_{(n,n_1,n_2,n_3)\in\N^4}
\frac{n^{2\sigma}\gamma^{2}(n,n_1,n_2,n_3)}{(1+|n^2-n_1^2+n_2^2-n_3^2|)^{\beta}(n_1 n_2 n_3)^{2}}\,.
\end{equation}
Using Lemma~\ref{l1} ans Lemma~\ref{ihp}, we obtain that with a suitable choice of $\beta<1$ one has
\begin{eqnarray*}
(\ref{8})& \leq &
C\sum_{(n_1,n_2,n_3)\in\N^3}
\frac{(1+|n_1^2-n_2^2+n_3^2|)^{\sigma}(\min(n_1,n_2,n_3))^{2}}
{(n_1 n_2 n_3)^{2}}
\\
& \leq &
C\sum_{(n_1,n_2,n_3)\in\N^3}
\frac{(\max(n_1,n_2,n_3))^{2\sigma}(\min(n_1,n_2,n_3))^{2}}
{(n_1 n_2 n_3)^{2}}
\\
& \leq &
C\sum_{n_3\leq n_2\leq n_1}
\frac{n_1^{2\sigma}n_2 n_3}{(n_1 n_2 n_3)^{2}}
\leq
C\sum_{n_1=1}^{\infty}\frac{n_1^{2\sigma}(\log(1+n_1))^{2}}{n_1^2}<\infty.
\end{eqnarray*}
Let us next analyse (\ref{8bis}).
Using the rapid decay of $|\widehat{\psi_2}|^{2}$, we can eliminate the $\tau$ integration and arrive at
\begin{equation*}
(\ref{8bis})\leq
C\sum_{(n,n_1,n_2,n_3)\in\N^4}\frac{n^{2\sigma}}{(1+|n_2^{2}-n^{2}|)^{\beta}}
\frac{\gamma(n,n_1,n_1,n_2)\gamma(n,n_2,n_3,n_3)}{(n_1 n_1 n_2)(n_2 n_3 n_3)}\,.
\end{equation*}
Using Lemma~\ref{l1} ans Lemma~\ref{ihp}, we obtain that with a suitable choice of $\beta<1$ one has
\begin{equation*}
(\ref{8bis})\leq
C\sum_{(n_1,n_2,n_3)\in\N^3}
\frac{n_2^{2\sigma} \min(n_1,n_2)\min (n_3,n_2) }{(n_1 n_2 n_3)^{2}}\,.
\end{equation*}
Let us fix $\varepsilon>0$ such that $\sigma+\varepsilon<1/2$.
Therefore, we can write
$$
(\ref{8bis})\leq
C\sum_{(n_1,n_2,n_3)\in\N^3}\frac{n_2^{2\sigma+2\varepsilon}(n_1 n_3)^{1-\varepsilon}}{(n_1 n_2 n_3)^{2}}
<\infty\,.
$$
This completes the proof of Proposition~\ref{thm4}.
\end{proof}


\begin{thebibliography}{10}
\bibitem{Ramona} R. Anton, {\it Cubic nonlinear Schr\"odinger equation on three dimensional balls 
with radial data}, Preprint~2006.
%
\bibitem{AT} A. Ayache, N. Tzvetkov, {\it $L^p$ properties of Gaussian random series},
to appear in Trans. AMS.
%
\bibitem{Bo1}
J. Bourgain, {\it Periodic nonlinear Schr\"odinger equation and invariant measures}, 
Comm. Math. Phys. 166 (1994) 1-26.
%
\bibitem{Bo2}
J. Bourgain, {\it Invariant measures for the 2D-defocusing nonlinear Schr\"odinger equation}, 
Comm. Math. Phys. 176 (1996) 421-445.
%
\bibitem{BGT} N. Burq, P. G\'erard, N. Tzvetkov, 
{\it Multilinear eigenfunction estimates and global existence for the three dimensional nonlinear 
Schr\"odinger equations}, Ann. ENS, 38 (2005) 255-301.
%
\bibitem{BGT-zonal} N. Burq, P. G\'erard, N. Tzvetkov, {\it Zonal low regularity solutions of the 
nonlinear Schr\"odinger equation on $S^d$}, Unpublished manuscript, Summer~2002.
%
\bibitem{CCT} M. Christ, J. Colliander, T. Tao, {\it Ill-posedness for nonlinear Schr\"odinger 
and wave equations }, Preprint~2003. 
%
\bibitem{Gi} J. Ginibre,
{\it Le probl\`eme de Cauchy pour des EDP semi-lin\'eaires p\'eriodiques en variables d'espace 
(d'apr\`es Bourgain)}, S\'eminaire Bourbaki, Exp. 796, Ast\'erisque 237 (1996) 163-187.
%
\bibitem{LRS} J. Lebowitz, R. Rose, E. Speer, {\it Statistical dynamics of the
nonlinear Schr\"odinger equation}, J. Stat. Physics V 50 (1988) 657-687.
%
\bibitem{KS} S. Kuksin, A. Shirikyan,
{\it Randomly forced CGL equation : stationary measures and the inviscid limit}, J. Phys A 37 (2004) 1-18.
%
\bibitem{SW} E. Stein and G. Weiss,
{\it Introduction to Fourier analysis on euclidean spaces}, 
Princeton University Press, Princeton, N.J., 1971. Princeton Mathematical Series, No. 32.
%
\bibitem{Tz} N.~Tzvetkov, {\it Invariant measures for the Nonlinear Schr\"odinger equation on the disc},
Dynamics of PDE 3 (2006) 111-160.
%
%
\bibitem{Zh} P. Zhidkov, {\it KdV and nonlinear Schr\"odinger equations :
Qualitative theory}, Lecture Notes in Mathematics 1756, Springer 2001.
\end{thebibliography}
\end{document}